\documentclass[a4paper,11pt,oneside]{article}

\usepackage{amsmath,amssymb}
\usepackage{epsfig}
\usepackage[utf8]{inputenc}
\usepackage{textcomp}
\usepackage{algorithm}
\usepackage{algpseudocode}
\usepackage{algcompatible}
\usepackage{color}
\usepackage{graphicx}
\usepackage{caption} 
\usepackage{subcaption}
\usepackage{bm}
\usepackage{amsthm}
\usepackage{mathrsfs}
\usepackage{dsfont}
\usepackage{url}

\usepackage{geometry}
\newtheorem{thm}{Theorem}[section]

\newtheorem{definition}[thm]{Definition}

\newtheorem{proper}[thm]{Property}

\newtheorem{remark}[thm]{Remark}

\DeclareMathOperator*{\argmin}{arg\,min}

\newcommand{\Prob}{\mathbb{P}}

\newcommand{\E}{\mathbb{E}}

\newcommand{\g}[1]
{
\mathbf{ #1}
}

\newcommand{\bs}[1]
{
\boldsymbol{ #1}
}

\algnewcommand\INPUT{\item[\textbf{Input:}]}%
\algnewcommand\OUTPUT{\item[\textbf{Output:}]}%

\title{\textbf{\textsc{Physics-informed cluster analysis and a priori efficiency criterion for the construction of local reduced-order bases}}}
\author{Thomas DANIEL\footnotemark[1]  \footnotemark[2] \footnotemark[3] , Fabien CASENAVE\footnotemark[1] , Nissrine AKKARI\footnotemark[1] , Ali KETATA\footnotemark[2] ,\\David RYCKELYNCK\footnotemark[2]}
\date{December 3, 2021}

\begin{document}
\maketitle

\renewcommand{\thefootnote}{\fnsymbol{footnote}}
\footnotetext[1]{ \ SafranTech, Rue des Jeunes Bois, Ch\^ateaufort, 78114 Magny-les-Hameaux (France).}
\footnotetext[2]{ \ MINES ParisTech, PSL University, Centre des mat\'{e}riaux (CMAT), CNRS UMR 7633, BP 87, 91003 Evry (France).}
\footnotetext[3]{ \ Correspondence: thomas.daniel@mines-paristech.fr }
\renewcommand{\thefootnote}{\arabic{footnote}}

\begin{abstract}
Nonlinear model order reduction has opened the door to parameter optimization and uncertainty quantification in complex physics problems governed by nonlinear equations. In particular, the computational cost of solving these equations can be reduced by means of local reduced-order bases. This article examines the benefits of a physics-informed cluster analysis for the construction of cluster-specific reduced-order bases. We illustrate that the choice of the dissimilarity measure for clustering is fundamental and highly affects the performances of the local reduced-order bases. It is shown that clustering with an angle-based dissimilarity on simulation data efficiently decreases the intra-cluster Kolmogorov $N$-width. Additionally, an a priori efficiency criterion is introduced to assess the relevance of a ROM-net, a methodology for the reduction of nonlinear physics problems introduced in our previous work in [T.~Daniel, F.~Casenave, N.~Akkari, D.~Ryckelynck, Model order reduction assisted by deep neural networks (ROM-net), \textit{Advanced Modeling and Simulation in Engineering Sciences} 7 (16), 2020]. This criterion also provides engineers with a very practical method for ROM-nets' hyperparameters calibration under constrained computational costs for the training phase. On five different physics problems, our physics-informed clustering strategy significantly outperforms classic strategies for the construction of local reduced-order bases in terms of projection errors.

\end{abstract}

\section{Introduction}

Differential equations are widely used for the mathematical modeling of physical phenomena. These differential equations involve boundary conditions, initial conditions, constants and source terms that can be considered as parameters of the physics problem. For well-posed parametrized differential equations, any point of the parameter space is associated to one single solution. Under the third Hadamard well-posedness condition, the solution is a continuous function of the parameters. Therefore, any connected set in the parameter space defines a connected set in the solution space, called \textit{solution manifold}. This concept can be extended to any quantity of interest, be it an internal variable or a function of the solution. The solution manifold can be interpreted as the support of a probability density function {for the solution} in uncertainty propagation, when a probabilistic model is given to describe uncertainties on the parameters. The concept of solution manifold also appears in design optimization where the objective is to minimize a solution-dependent cost function by modifying parameters such as material constants, microstructural properties, or geometrical characteristics, for instance. Both uncertainty quantification and design optimization are many-query problems since they require solving the parametrized differential equations for a large number of points in the parameter space, which is sometimes prohibitive. To mitigate the computational cost related to such applications, many \textit{model order reduction} methods~\cite{10.5555/2568435, keiper2018reduced} have been developed, including methods based on tensor decompositions (e.g. the \textit{Proper Generalized Decomposition}~\cite{PGDreview, PGDbook}) and projection-based methods (e.g. the \textit{Reduced Basis method}~\cite{RBmethodPrudhomme, RozzaReducedBasis} and the \textit{POD Galerkin method}~\cite{cordier:hal-00417819, RowleyPODGalerkin}). \textit{Projection-based model order reduction} consists in computing an approximate solution in a low-dimensional subspace of the solution space, which can give accurate predictions provided that the solution manifold is embedded in a low-dimensional space.

The potential of such numerical methods is related to the maximum distance between a point of the solution manifold and its orthogonal projection onto the approximation space. This distance is used to define the \textit{Kolmogorov $N$-width} measuring the worst-case error for the best $N$-dimensional approximation space. The behavior of the Kolmogorov width when increasing the dimension $N$ provides information about the reducibility of a given physics problem: slowly decaying Kolmogorov widths indicate that increasing the dimension of the {linear approximation subspace} does not significantly improve the quality of the approximate solution. The asymptotic behavior of the Kolmogorov width is studied in~\cite{10.1093/imanum/dru066, GREIF2019216}. Problems combining large values of the Kolmogorov widths with low decay rates are not reducible, which means that one cannot compute accurate approximate solutions at low computational costs. This is the case in particular when considering wave propagation~\cite{GREIF2019216} and {advection-dominated problems~\cite{CagniartBook, nonino2019overcoming, doi:10.1137/110823158, dahmen_plesken_welper_2014, rim2020manifold, taddei2020spacetime}}. For such problems, in the projection-based model order reduction community, \cite{localROB, localROB2}~suggested using multiple \textit{local} approximation spaces dedicated to different subsets of the solution manifold. These local approximation spaces are spanned by local reduced-order bases (ROBs), usually computed with the Proper Orthogonal Decomposition (POD~\cite{Lumley, cordier:hal-00417819}). ROM interpolation techniques have also been extensively studied in~\cite{Interpolation0, Interpolation1, Interpolation2, Interpolation3, Interpolation4, Interpolation5, Interpolation6, doi:10.2514/1.J050233, AMSALLEM2016373, Interpolation7, Interpolation8, CHOI2020109787}. More recently, \cite{LEE2020108973}~introduced two nonlinear model order reduction techniques called \textit{manifold Galerkin method} and \textit{manifold least-squares Petrov-Galerkin method}, where a deep convolutional autoencoder is used for nonlinear dimensionality reduction. In~\cite{KimChoi2020}, this idea is extended to a hyper-reduction framework, using a shallow masked autoencoder with fully-connected layers.

The present article focuses on the use of multiple local ROBs to tackle the slow decrease of the Kolmogorov $N$-width, because of the compatibility of this approach with the Galerkin method and its classical extension to hyper-reduction methods~\cite{EIM, Ryckelynck2005, ECSW, ECM, Grimberg2020}. Splitting a non-reducible problem into multiple reducible ones can be achieved with \textit{cluster analysis}. {Cluster analysis belongs to unsupervised learning tasks and can be defined as the search of groups (or \textit{clusters}) of similar objects in a database. The choice of the clustering algorithm depends on the underlying motivation and thus on the clusters' topological properties that are expected. Representative-based algorithms~\cite{kmeans, aggarwal2015data, kMedoidsPAM} use a dissimilarity measure to assign each object of the database to the cluster corresponding to the closest representative object, leading to compact clusters. We refer the reader to the books~\cite{aggarwal2013data}, \cite{aggarwal2015data}~(chapters 6 and 7), \cite{Hastie2005TheEO}~(section 14.3), or articles~\cite{10.1145/331499.331504, Saxena2017} for more details about clustering algorithms.}
{In cases where efficient error bounds for reduced representations are available, ideal cases being for linear problems with affine parameter dependencies, a powerful framework has been proposed in~\cite{balabanov2021randomized}. Random sketching is used to improve the orthogonal greedy algorithm, and then used in a greedy algorithm where the solution is approximated in a subspace spanned by vectors selected online in a dictionary of candidate basis vectors.
}

The present work studies \textit{physics-informed} clustering strategies for the construction of dictionaries of local ROBs, {for complex nonlinear problems parametrized by a field}. Physics-informed cluster analysis consists in clustering the parameter space by means of a dissimilarity measure which involves physical quantities obtained when solving the physics problem. In other words, clusters in the parameter space are implicitly defined as the preimages of clusters found in a database of numerical simulation results. {In practice, this means that a clustering algorithm is applied in the solution space, as proposed for the first time in~\cite{localROB, localROB2}. The focus is on finding a clustering strategy that is appropriate for model order reduction purposes. In~\cite{AmsallemHaasdonk}, it was noticed that clustering based on the Euclidean distance (or $L^2$ distance) in the solution space was not adapted for the construction of local reduced-order models (ROMs), which led to the definition of \textit{projection-error based local ROMs} (PEBL-ROM) where the solution space is hierarchically partitioned using the projection error as a dissimilarity criterion.} {We propose here to work with the sine dissimilarity related to \textit{relative} projection errors instead, giving a symmetric dissimilarity measure that can be plugged into a representative-based clustering method. We show that this dissimilarity measure corresponds to the Hilbert-Schmidt distance between the projections onto the snapshots' approximation spaces, which is of particular interest when using Galerkin projection to solve the governing equations.}

{Our main contributions are the development of a physics-informed clustering strategy based on the sine dissimilarity and k-medoids clustering, an automatic snapshot selection procedure for the construction of POD bases, and an \textit{a priori} efficiency criterion enabling hyperparameters calibration for dictionary-based ROM-nets~\cite{ROM-net} where a classifier is trained to automatically recommend the best ROM in the dictionary for a given point in the parameter space without computing the dissimilarity online.} {Section~\ref{SectionDicoROMs} gives an overview of model order reduction methods and of the techniques that have been developed to deal with non-reducible problems, including dictionaries of local ROMs. Section~\ref{ClusteringBackground} presents representative-based clustering algorithms, and in particular k-medoids clustering. Our physics-informed clustering strategy and a priori efficiency criterion are introduced in Sections~\ref{SectionPhyInformedClus} and~\ref{SectionROMnetEffCrit} respectively. Applications to various physics problems are developed in Section~\ref{SectionApplication} and show the importance of choosing an appropriate dissimilarity measure for clustering.}

\section{Model order reduction background}
\label{SectionDicoROMs}

Let us consider a physics problem described by the following parametrized differential equation:\begin{equation}
\mathcal{D}(u;x) = 0,
\label{PDE}
\end{equation} where $u$ is the primal variable belonging to a Hilbert space $\mathcal{H}$ whose inner product is denoted by $\langle .,. \rangle_{\mathcal{H}}$, $x$ denotes the parameters of the problem, and $\mathcal{D}$ is an operator involving a differential operator and operators for initial conditions and/or boundary conditions. Equation~\eqref{PDE} can be a system of ordinary differential equations or partial differential equations depending on the physics problem. Let us assume that this physics problem is well-posed in the sense of Hadamard, that is to say that there exists a unique solution $u(x)$ for any parameter $x$, and that this solution changes continuously with $x$. Let us introduce the set $\mathcal{X}$ of all the possible parameters $x$. The \textit{solution manifold} $\mathcal{M}$ is defined by:\begin{equation}
\mathcal{M} := u(\mathcal{X}) = \{ u(x) \ | \ x \in \mathcal{X} \}.
\label{solutionManifold}
\end{equation}

\subsection{Introduction to model order reduction}

Model order reduction~\cite{10.5555/2568435, keiper2018reduced} is a discipline in numerical analysis consisting in replacing a computationally expensive high-fidelity model by a fast reduced-order model (ROM) to calculate approximate solutions of some complex physics equations. A ROM can be either a data-driven metamodel (or surrogate model) calibrated with a regression algorithm, or a physics-based model obtained by numerical methods such as the Proper Generalized Decomposition~\cite{PGDreview, PGDbook}, the Reduced Basis method~\cite{RBmethodPrudhomme, RozzaReducedBasis}, and the POD Galerkin method~\cite{cordier:hal-00417819, RowleyPODGalerkin}, among others. It is generally used for parametrized equations whose solution must be known for different points in the parameter space. As in machine learning, a model order reduction procedure starts by a \textit{training phase} (or \textit{offline stage}) where the ROM is built from some training data. The ROM is then used on test data in an \textit{exploitation phase} (or \textit{online stage}). In the training phase, high-fidelity solutions, called \textit{snapshots}, are computed with the high-fidelity model for different points of the parameter space to get a sampled representation of the solution manifold. The model order reduction algorithm analyzes these snapshots to learn how the solution is affected by parameter variations. Contrary to usual machine learning problems, the amount of training data is limited because the high-fidelity model giving snapshots is time-consuming and costly. The selection of relevant points in the parameter space can be optimized to ensure that the snapshots are representative of the behavior of the solution, like in the greedy approach of the Reduced Basis method where an a posteriori error estimator is used to select snapshots. Given the cost of computing snapshots in the training phase, a ROM is profitable only if it is extensively used in the exploitation phase. This paper addresses issues that are specific to \textit{projection-based model order reduction} (e.g. POD Galerkin, Reduced Basis method) where the approximate solution is obtained by solving the physics equations with the Galerkin method on a well-chosen reduced-order basis (ROB).

\subsection{The Proper Orthogonal Decomposition (POD)}

It is now assumed that Equation~\eqref{PDE} defines a parametrized partial differential equation whose solution $u(x)$ for a given point $x\in\mathcal{X}$ in the parameter space is a function of space and time defined on $\Omega \times [0;t_f]$, with $\Omega \subset \mathbb{R}^{\alpha}$, $\alpha=1$, $2$ or $3$ and $t_f \in \mathbb{R}_{+}^*$. Most of the time, the Hilbert space $\mathcal{H}$ is a subspace of the Lebesgue space $L^{2}(\Omega \times [0;t_f])$ of square-integrable functions. However, parameters $x$ and time $t$ can be considered together in a variable $\chi$ called \textit{generalized parameters} living in the set $\tilde{\mathcal{X}} = \mathcal{X} \times [0;t_f ]$. Therefore, the solution $u(\chi)$ belongs to the space $L^{2}(\Omega)$. The Stiefel manifold $V(N,L^{2}(\Omega))$ represents the set of all orthonormal $N$-frames in $L^{2}(\Omega)$. For two square-integrable functions $f$ and $g$, the notation $\langle f,g \rangle_{L^{2}(\Omega)}$ stands for the $L^{2}(\Omega)$ inner product $\int_{\Omega} f(\bs{\xi})g(\bs{\xi})d\bs{\xi}$. The following definition gives a theoretical continuous definition of the proper orthogonal decomposition (POD~\cite{Lumley, cordier:hal-00417819}), also known as the Karhunen-Lo\`eve decomposition or principal component analysis:

\begin{definition}
[POD basis] Let $u : \tilde{\mathcal{X}} \rightarrow L^{2}(\Omega)$. A \textit{POD basis} $\{ \psi_{k}^{*} \}_{1 \leq k \leq N} \in V(N,L^{2}(\Omega))$ of order $N \in \mathbb{N}^*$ of $u$ is a solution of the following optimization problem:\begin{equation}
\{ \psi_{k}^{*} \}_{1 \leq k \leq N} \in \underset{\{ \psi_{k} \}_{1 \leq k \leq N} \in V(N,L^{2}(\Omega))}{\argmin} \ \int_{\chi \in \tilde{\mathcal{X}}} || u(\chi) - \sum_{k=1}^N \langle u(\chi), \psi_k \rangle_{L^{2}(\Omega)} \psi_k   ||_{L^{2}(\Omega)}^2 \ d\chi .
\label{ContinuousPOD}
\end{equation}
\end{definition}

The sum in Equation~\eqref{ContinuousPOD} is the proper orthogonal decomposition of order $N$ of $u$. When $N \rightarrow + \infty$, the approximation error given by the minimum of the cost function in Equation~\eqref{ContinuousPOD} tends towards zero (Theorem 4 in~\cite{HenriPOD}).

Let $\mathcal{H}$ be a Hilbert space with an orthonormal basis $\{ e_i \}_{i \in \mathbb{N}^*}$, and $A: \mathcal{H} \rightarrow \mathcal{H}$ a linear operator. We define the Hilbert-Schmidt function $\Lambda_{HS(\mathcal{H})}$ as:\begin{equation}
\Lambda_{HS(\mathcal{H})}(A) := \displaystyle \sqrt{\sum_{i=1}^{\infty} || A(e_i) ||_{\mathcal{H}}^2},
\end{equation} which can potentially take infinite values. A linear operator $A$ on a Hilbert space $\mathcal{H}$ is a \textit{Hilbert-Schmidt operator} if $\Lambda_{HS(\mathcal{H})}(A)$ is finite. As shown in~\cite{gohberg1990classes} (Chapter VIII, Theorem 2.3), the set $HS(\mathcal{H})$ of all Hilbert-Schmidt operators on $\mathcal{H}$ is a Hilbert space with respect to the following inner product:\begin{equation}
\langle A , B  \rangle_{HS(\mathcal{H})} := \sum_{i=1}^{\infty} \langle A(e_i), B(e_i)  \rangle_{\mathcal{H}}.
\label{HSinnerProd}
\end{equation} The Hilbert-Schmidt function $\Lambda_{HS(\mathcal{H})}$ is actually the norm induced by this inner product, and corresponds to the Frobenius norm for matrices when the vector space $\mathcal{H}$ is finite-dimensional. We now use the more conventional notation $|| A ||_{HS(\mathcal{H})} := \Lambda_{HS(\mathcal{H})}(A)$ for the Hilbert-Schmidt norm. The Hilbert-Schmidt inner product and norm are independent of the choice of the basis $\{ e_i \}_{i \in \mathbb{N}^*}$, see Proposition~9.18 in Chapter~9 of~\cite{cheverry:cel-01587623}, which will be useful for proofs of some properties of the dissimilarity measure introduced in this paper. The POD is highly related to the theory of Hilbert-Schmidt operators. In~\cite{Djouadi2009, Djouadi2012}, it is shown that the POD optimization problem is equivalent to finding the optimal approximation of a Hilbert-Schmidt operator related to $u$ by a finite rank operator in the Hilbert-Schmidt norm. The POD basis functions can also be obtained from the eigenfunctions of the Hilbert-Schmidt integral operator~\cite{cordier:hal-00417819} $\mathcal{R}_{u}$:\begin{equation}
L^2(\tilde{\mathcal{X}})\ni\varphi\mapsto \mathcal{R}_u(\varphi)\in L^2(\tilde{\mathcal{X}}),
\end{equation}
where $\mathcal{R}_u(\varphi)$ is defined by:
\begin{equation}
\tilde{X}\ni\chi\mapsto \mathcal{R}_u(\varphi)(\chi):=\int_{\chi'\in\tilde{\mathcal{X}}}\left<u(\chi),u(\chi')\right>_{L^2\left(\Omega\right)}\varphi(\chi')d\chi'\in \mathbb{R}.
\end{equation}

In this work, we keep the explicit distinction between the time $t$ and the parameters $x \in \mathcal{X}$ rather than working on the generalized parameters $\chi \in \tilde{\mathcal{X}}$, because we do not consider the time as a clustering variable. Nonetheless, spatio-temporal functions $f\in L^{2}(\Omega \times [0;t_f])$ are considered as trajectories $(f(.,t))_{t \in [0;t_f]}$ in the Hilbert space $L^{2}(\Omega)$. In other words, such functions are seen as functions defined on $\Omega$ and parametrized by the time. For this reason, the manifold $\mathcal{M}$ is rather defined by:\begin{equation}
\mathcal{M} := \{ u(x)(.,t) \ | \ x\in\mathcal{X}, \ t \in [0;t_f] \},
\end{equation} and the approximation spaces are subspaces of $L^{2}(\Omega)$, leading to an approximate solution expressed  as a time-dependent linear combination of basis functions defined on $\Omega$.

In practice, we are given a finite set of $m$ points $x_i$ of the parameter space $\mathcal{X}$, for which high-fidelity solutions $u(x_i)$ are computed in a high-dimensional approximation space whose dimension is denoted by $\mathcal{N}$. These solutions, called snapshots, provide information about the behavior of the physical system and give a sampled version of the solution manifold. The POD is applied as a linear dimensionality reduction technique, processing this information to build a ROB that can be used to accelerate future numerical simulations for new parameters.

\begin{definition}
[POD basis construction] Given an integer $N\leq \mathcal{N}$, a POD basis $\{ \psi_{k}^{*} \}_{1 \leq k \leq N} \in V(N,L^{2}(\Omega))$ is computed from the snapshots $\{ u(x_i)\}_{1\leq i \leq m}$ as a solution of the following optimization problem:\begin{equation}
\{ \psi_{k}^{*} \}_{1 \leq k \leq N} \in \underset{\{ \psi_{k} \}_{1 \leq k \leq N} \in V(N,L^{2}(\Omega))}{\argmin} \ \sum_{i=1}^{m} || u(x_i) - \sum_{k=1}^N \langle u(x_i), \psi_k \rangle_{L^{2}(\Omega)} \psi_k   ||_{L^{2}(\Omega \times [0;t_f])}^2 .
\label{PODdef}
\end{equation}
\end{definition}

The uniqueness of the POD basis is obtained by specifying a construction algorithm, such as the Snapshot POD~\cite{Sirovich, Chatterjee} or the singular value decomposition (SVD) for instance. By construction, the subspace spanned by the ROB minimizes the projection errors of the snapshots $u(x_i)$. The optimality of the POD basis is discussed and illustrated in~\cite{Meyer2003}. In practice, when using a numerical procedure to solve Equation~\eqref{PDE}, for example the finite-element method with a time-stepping scheme, the coordinates of the snapshots in the finite-element basis are stored in columns in a matrix $\g{Q} \in \mathbb{R}^{\mathcal{N}\times m n_t}$ called snapshots matrix, with $n_t$ being the number of time steps. The coordinates of the POD modes $\psi_k$ are given in the $N$ first columns of the matrix $\g{M}^{-1/2}\g{V}$, where $\g{M} \in \mathbb{R}^{\mathcal{N} \times \mathcal{N}}$ is the finite-element mass matrix and $\g{V}\in\mathbb{R}^{\mathcal{N}\times \textrm{rank}(\g{Q})}$ is the matrix of left singular vectors in the SVD of the snapshots matrix $\g{Q}$, when indexing the singular values in decreasing order. The decay rate of the singular values of the snapshots matrix $\g{Q}$ is related to the behavior of the sequence of Kolmogorov widths. It enables evaluating the reducibility of the physics problem. When computing a POD basis for a variable defined at integration points rather than the finite-element mesh nodes, for the purpose of applying Gappy-POD after hyper-reduced simulations, the POD modes are simply given by the $N$ first left singular vectors in the SVD of the corresponding snapshots matrix.

\subsection{Non-reducible problems}

{Approximate solutions $\tilde{u}(x)$ of Equation~\eqref{PDE} can be obtained by solving the PDEs on a finite-dimensional subspace $\mathcal{H}_{N} \in \textrm{Gr}(N,\mathcal{H})$ spanned by a ROB, where the Grassmannian $\textrm{Gr}(N,\mathcal{H})$ is the set of all $N$-dimensional subspaces of $\mathcal{H}$. The best approximation of the solution in $\mathcal{H}_N$ for a given parameter $x$ is the orthogonal projection $\pi_{\mathcal{H}_N}(u(x))$ of the theoretical solution $u(x)$ onto the approximation space:\begin{equation}
\pi_{\mathcal{H}_N}(u(x)) = \underset{v \in \mathcal{H}_N}{\argmin} \ || u(x) - v ||_{\mathcal{H}},
\end{equation} with $|| . ||_{\mathcal{H}}$ denoting the norm induced by the inner product of the Hilbert space $\mathcal{H}$.}

\noindent {The Kolmogorov $N$-width is defined by:\begin{equation}
d_{N}(\mathcal{M}) := \underset{\mathcal{H}_{N} \in \textrm{Gr}(N,\mathcal{H})}{\inf} \  \underset{u \in \mathcal{M}}{\sup}  \ \underset{v \in \mathcal{H}_N}{\inf} || u - v ||_{\mathcal{H}} = \underset{\mathcal{H}_{N} \in \textrm{Gr}(N,\mathcal{H})}{\inf} \  \underset{u \in \mathcal{M}}{\sup}  \  || u - \pi_{\mathcal{H}_N}(u) ||_{\mathcal{H}},
\label{KolmogorovWidth}
\end{equation} and quantifies how well the solution manifold $\mathcal{M}$ can be approximated by searching approximate solutions in a $N$-dimensional subspace of $\mathcal{H}$. The Kolmogorov $N$-width corresponds to the worst projection error on the best $N$-dimensional approximation space. For a fixed solution manifold $\mathcal{M}$, the sequence $\left( d_{N}(\mathcal{M})\right)_{N\in\mathbb{N}}$ is decreasing, which means that approximation errors get lower when increasing the dimension of the approximation space.}

For some problems, the Kolmogorov width slowly decays when increasing the dimension $N$ of the approximation space. For these \textit{non-reducible} problems, the dimension $N$ of the {linear approximation space} giving a sufficiently small Kolmogorov width is generally too high to enable the fast computation of approximate solutions. Qualitatively, the solution manifold $\mathcal{M}$ covers too many independent directions to be embedded in a low-dimensional subspace. To address this issue, several techniques have been developed:
\begin{itemize}
\item Problem-specific methods tackle the difficulties of some specific physics problems that are known to be non-reducible, such as advection-dominated problems which have been largely investigated, {for instance in~\cite{doi:10.1137/110823158, dahmen_plesken_welper_2014, rim2020manifold, taddei2020spacetime, PhysRevE.89.022923, doi:10.1137/17M1140571, MadayConvection}}.
\item Online-adaptive model reduction methods update the ROM in the exploitation phase by collecting new information online as explained in~\cite{Zimmermann2017}, in order to limit extrapolation errors when solving the parametrized governing equations in a region of the parameter space that was not explored in the training phase. The ROM can be updated for example by querying the high-fidelity model when necessary for basis enrichment~\cite{Ryckelynck2005, 10.1145/1618452.1618469, doi:10.1137/151003660, CasenaveAkkari, he2020insitu}. Other methods propose enrichment procedures that do not require solving the equations with the high-fidelity model, whose complexity scales linearly with~(\cite{doi:10.1137/140989169, doi:10.1137/19M1257275}) or is independent of~(\cite{ETTER2020112931}) the dimension of the high-fidelity model.
\item ROM interpolation methods~\cite{Interpolation0, Interpolation1, Interpolation2, Interpolation3, Interpolation4, Interpolation5, Interpolation6, doi:10.2514/1.J050233, AMSALLEM2016373, Interpolation7, Interpolation8, CHOI2020109787} use interpolation techniques on Grassmann manifolds or matrix manifolds to adapt the ROM to the parameters considered in the exploitation phase by interpolating between two precomputed ROMs.
\item Dictionaries of basis vector candidates enable building a parameter-adapted ROM in the exploitation phase by selecting a few basis vectors. This technique is presented in~\cite{doi:10.1137/120873868, Kaulmann2012ONLINEGR} for the Reduced Basis method.
\item Dictionaries of ROMs rely on the construction of several local ROMs adapted to different regions of the solution manifold. These local ROMs can be obtained by partitioning the time interval~\cite{Drohmann2010, Dihlmann2011}, the parameter space~\cite{Drohmann2010, Eftang2010, Haasdonk2011, LDEIM2014, he2020insitu, doi:10.2514/6.2020-0418, kapteyn2020physicsbased}, or the {solution space~\cite{localROB, localROB2, LDEIM2014, Amsallem2015localHROM, AmsallemHaasdonk, RyckelynckComputerVision, ROM-net, Grimberg2020}.} Local ROMs have been used both with the Reduced Basis method and the POD Galerkin method. In the same vein as online-adaptive model reduction methods, local ROBs can be adapted online using for example a low-rank SVD update method, as in~\cite{localROB2, Amsallem2015localHROM} when switching from one local ROB to another or in~\cite{he2020insitu} when an error indicator detects extrapolation errors. This concept of local ROMs should not be confused with another type of local (or localized) ROMs described in~\cite{BuhrIapichinoSmetana+2020+245+306}, where the ROMs are associated to subdomains of the computational domain, in the spirit of domain decomposition techniques.
\item Nonlinear manifold ROM methods~\cite{LEE2020108973, lee2020DeepConservation, KimChoi2020} learn a nonlinear embedding and project the governing equations onto the corresponding approximation manifold, by means of a nonlinear function mapping a low-dimensional latent space to the solution space. This function is the decoder of an undercomplete autoencoder trained with the mean squared error loss to compress the snapshots and reconstruct them from their compressed representations. In this way, the nonlinear manifold is approximated with one single nonlinear ROM. Classical linear ROMs are obtained when the autoencoder has only one hidden-layer with linear activation functions. In this case, the decoder simply returns a linear combination of the POD modes.
\end{itemize}

\subsection{Dictionaries of local reduced-order models}

This paper focuses on dictionaries of ROMs, where the solution manifold is partitioned to get a collection of subsets $\mathcal{M}_k \subset \mathcal{M}$ that can be covered by a dictionary of low-dimensional subspaces, enabling the use of linear ROMs. If $\{ \mathcal{M}_k \}_{k \in [\![ 1;K ]\!]}$ is a partition of $\mathcal{M}$, then:\begin{equation}
\forall k \in [\![ 1;K ]\!], \ \forall N \in \mathbb{N}^{*}, \quad d_{N}(\mathcal{M}_k) \leq d_{N}(\mathcal{M}).
\end{equation} For a given number $K$ of subsets, two partitions can be compared on the basis of the ratios $d_{N}(\mathcal{M}_k) / d_{N}(\mathcal{M})$. The idea of clustering training data to define local ROBs traces back to the work of D. Amsallem, K. Washabaug, M.J. Zahr and C. Farhat in 2012, published in~\cite{localROB, localROB2} and validated on nonlinear problems in computational fluid dynamics and fluid-structure-electric interactions. In these papers, the set of snapshots is partitioned with k-means clustering to define $K$ clusters represented by their means $\{\overline{u}_k \}_{1\leq k \leq K}$. One local ROB is computed for each cluster using the POD. In the exploitation phase, given the solution at the $i$-th time increment, one looks for the closest mean $\overline{u}_k$ in terms of the norm $|| . ||_{L^{2}(\Omega)}$ and computes the state of the solution at the $i+1$-th time increment with the corresponding local ROB. This technique has been used more recently in a hyper-reduction framework in~\cite{Amsallem2015localHROM, Grimberg2020}.

When using a clustering algorithm to partition the solution manifold, the quality of the partition is related to the choice of the clustering method and the dissimilarity measure $\delta$ used to group similar solutions on the manifold. Among physics-informed clustering strategies, \textit{i.e.} strategies incorporating simulation data to compute dissimilarities, \cite{localROB, localROB2, Amsallem2015localHROM, RyckelynckComputerVision, Grimberg2020} used k-means with Euclidean distances in the solution space or in a subspace of the solution space found by PCA, \cite{ROM-net} used k-medoids with the Grassmann distance between subspaces spanned by the trajectories of the solutions, {\cite{AmsallemHaasdonk} applied a hierarchical partitioning based on a binary tree structure with the projection error as dissimilarity criterion,} and~\cite{LDEIM2014} proposed working on the governing equations' nonlinear term, using either a variant of k-means with the DEIM~\cite{DEIM} residual as clustering criterion or k-means on a low-dimensional representation of the governing equations' nonlinear term obtained by a DEIM-based feature selection. It is recalled that k-means is a representative-based clustering algorithm equipped with the Euclidean distance, and that changing this distance leads to other clustering methods. The Local Decomposition Method~\cite{LocalPODGPR} also relies on a physics-informed clustering strategy even though no dissimilarity measure is used, because a Gaussian mixture model is applied to shock sensors computed from the field of a quantity of interest, which enables separating subsonic and transonic flows in computational fluid dynamics.

\subsection{Dictionary-based ROM-nets}

\begin{figure}[!h]
\centering
\includegraphics[scale=0.33]{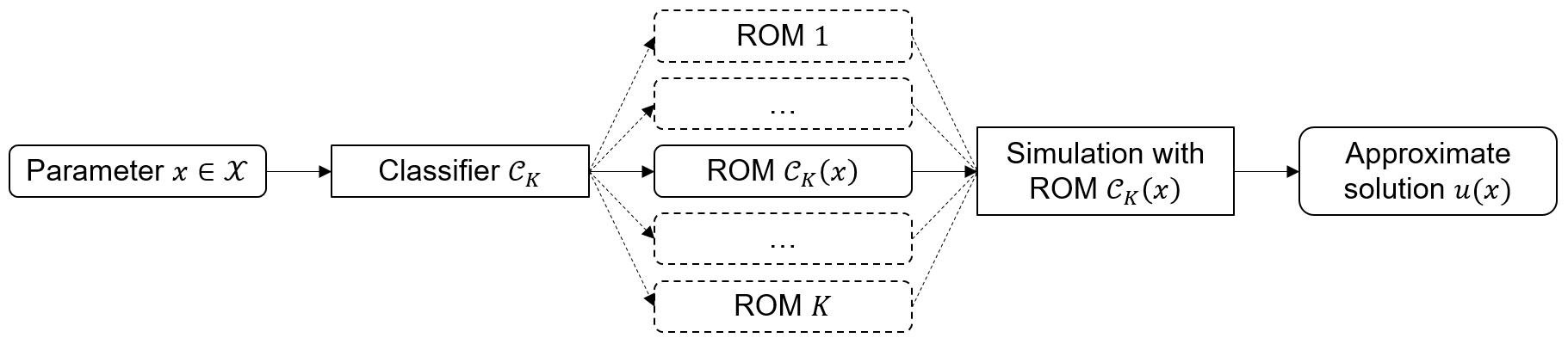}
\caption{Exploitation phase of a dictionary of $K$ local ROMs combined with a classifier $\mathcal{C}_K$ for automatic model recommendation. The parameter-based LDEIM uses a nearest neighbor classifier and the DEIM for the ROMs. The dictionary-based ROM-net uses artificial neural networks for the classification task.}
\label{ROM-net_online}
\end{figure}

The use of ROM dictionaries introduces the need for a model selection method that identifies the most suitable model in the dictionary. In~\cite{localROB, localROB2, Amsallem2015localHROM, Grimberg2020}, the local ROM is selected by finding the closest cluster representative from the current state of the solution with the Euclidean distance. When model selection is not straightforward and slows down the simulation process, one can use a classifier to learn the model selection task and enable fast model recommendation in the exploitation phase. {In~\cite{mainini2015surrogate}, global POD-bases for inputs and outputs of a black-box simulation model are constructed, classifiers are trained in the form of self-organizing maps in the POD coefficients space, and local surrogate models are trained -- while keeping a reduced representation on a global POD basis.} Dictionaries of ROMs with automatic model recommendation made by a classifier can be found in~\cite{RyckelynckComputerVision, LDEIM2014, ROM-net, doi:10.2514/6.2020-0418, kapteyn2020physicsbased}. To our knowledge, the idea of combining physics-informed clustering for the definition of local ROMs with a classifier for model recommendation came from the pioneering works of Peherstorfer, Butnaru, Willcox, and Bungartz on the Localized Discrete Empirical Interpolation Method (LDEIM~\cite{LDEIM2014}, 2014) and of Nguyen, Barhli, Mu\~noz and Ryckelynck on computer vision~\cite{RyckelynckComputerVision} in 2018. 
When the classification task is performed by deep neural networks and takes the parameters as inputs, this methodology is known as \textit{dictionary-based ROM-net}~\cite{ROM-net}, see Figure~\ref{ROM-net_online}.

\begin{figure}[!h]
\centering
\includegraphics[scale=0.33]{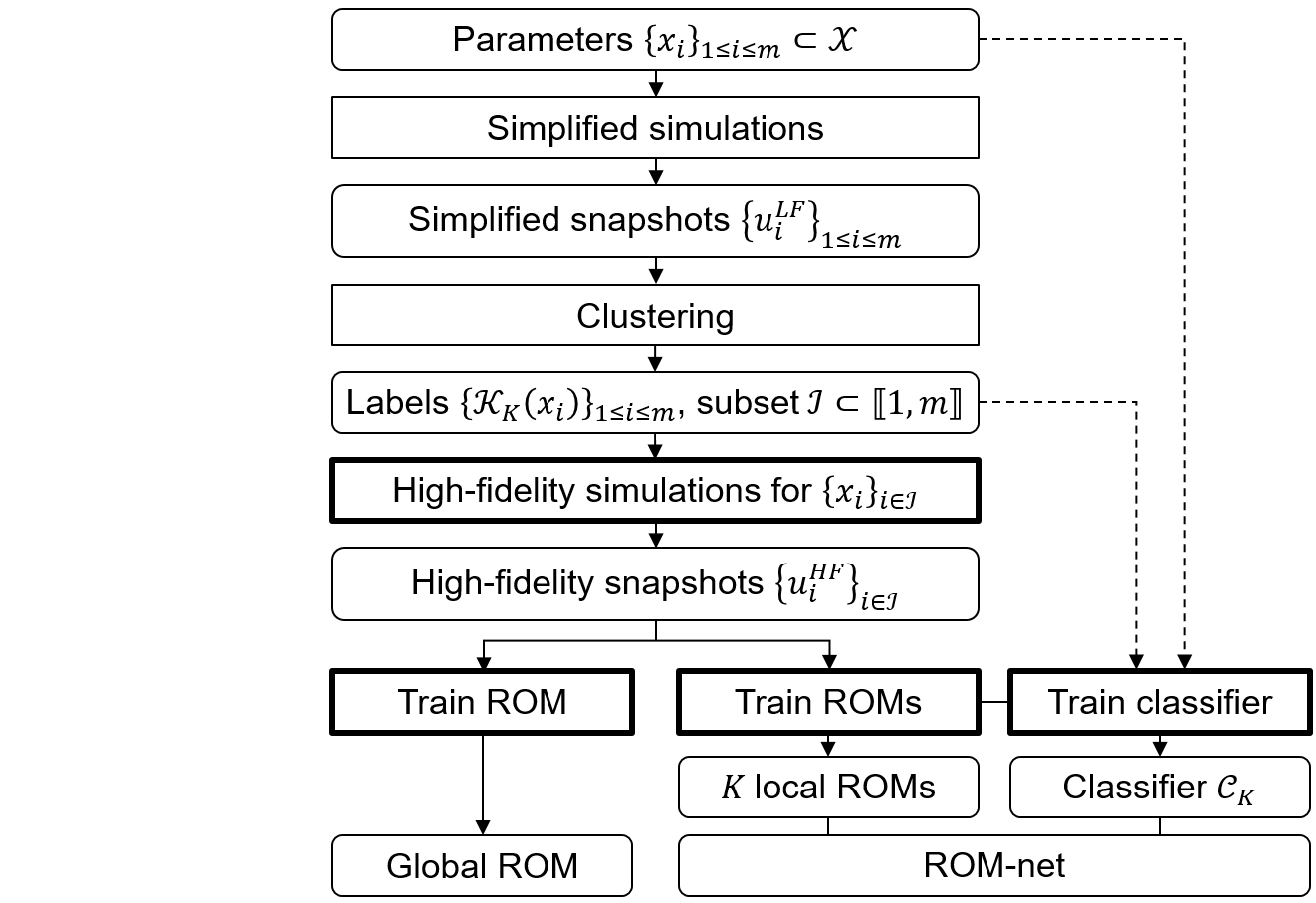}
\caption{Training phases of a dictionary-based ROM-net and a global ROM with a physics-informed clustering strategy.}
\label{ROM-net_offline}
\end{figure}

Dictionary-based ROM-nets consist in a dictionary of local ROMs and a classifier acting as a model selector, which enables the automatic adaptation of the ROM to the state and the environment of the physical system. The ROM-net's classifier (\textit{real classifier} denoted by $\mathcal{C}_K$ where $K$ is the number of local ROMs) approximates the theoretical \textit{perfect classifier} $\mathcal{K}_K$ returning the index of the best local ROM for a given point in the parameter space. For simplicity, the term \textit{dictionary-based ROM-net} (or simply \textit{ROM-net}) is used throughout this paper to refer to the methodology described by Figures~\ref{ROM-net_online} and~\ref{ROM-net_offline} {and Algorithm~\ref{alg:ROM-net_train}}, even when the classifier is not an artificial neural network.

Figure~\ref{ROM-net_offline} gives the main steps of the training phase of a dictionary-based ROM-net and draws a comparison with the construction of a global ROM benefiting from the ROM-net's physics-informed cluster analysis. The dictionary of local ROMs is built from clusters given by a physics-informed clustering procedure. First, a simplified version of the physics problem is solved for each input example of the training database. The simplified physics problem must be less computationally demanding than the target problem. In particular, it can be solved with a coarse mesh to reduce the dimension of the approximation space. The simplified simulations provide what we call \textit{simplified snapshots}: these snapshots cannot be exploited to build ROMs, but they give information about how the physical system reacts to parameter changes. The clustering algorithm finds clusters from the information contained in these simplified snapshots. In light of the clustering results, one must identify a few relevant training examples for which the target problem is solved to get \textit{high-fidelity snapshots}, that is, snapshots that well represent the solution manifold and can then be used for the construction of the local ROMs. The different needs in terms of training data for reduced-order modeling and machine learning can be seen through these two families of snapshots: the clustering and classification algorithms use information related to the simplified snapshots to get a sufficiently large training set, while the ROMs use a limited number of high-fidelity snapshots in order to learn to make predictions in a physics problem. This distinction between these two types of simulation data is essential when considering complex problems with many degrees of freedom. The three major differences between our work and the seminal works of~\cite{RyckelynckComputerVision} and~\cite{LDEIM2014} are the use of simplified simulations, the clustering strategy with a new ROM-oriented dissimilarity measure, and an a priori efficiency criterion introduced hereinafter. It is noteworthy that, among the four variants of the LDEIM, the parameter-based LDEIM with clustering of snapshots (section 4.2. of~\cite{LDEIM2014}) is the one that shares the more similarities with our work, since it applies clustering on simulation data and uses the parameters as inputs for the classifier. {The training algorithm is given in Algorithm~\ref{alg:ROM-net_train}.}

\noindent When using a dictionary-based ROM-net, two natural questions arise:
\begin{itemize}
\item Which dissimilarity measure and clustering algorithm should be used for model order reduction purposes?
\item Can the high-fidelity snapshots be automatically selected and how?
\end{itemize}

After clustering, the training phase of the dictionary-based ROM-net still includes expensive steps corresponding to boxes with thick lines in Figure~\ref{ROM-net_offline}, namely the computation of the high-fidelity snapshots, the construction of the local ROMs (which can involve a hyper-reduction algorithm), and the training of a classifier for automatic model recommendation. Therefore, an evaluation criterion is needed in order to assess the quality of the clusters before continuing the ROM-net's training phase, \textit{i.e.} right after the clustering step, see Figure~\ref{ROM-net_offline}. This criterion should enable the evaluation of the profitability of the ROM-net and the tuning of clustering hyperparameters, using the simplified snapshots only. Put briefly, in addition to the two aforementioned questions, this work must also address the following issues:
\begin{itemize}
\item Is it possible to define a simple practical method to select good hyperparameters (number of clusters, number of POD modes, number of high-fidelity snapshots)?
\item {Can one define an efficiency criterion computable after the clustering step to evaluate the expected performances of the ROM-net with respect to a single global ROM?}
\end{itemize} 

\begin{remark}
Dictionary-based ROM-nets use machine learning to assist model order reduction procedures in the training and exploitation phases. It does not replace physics models by regression models, since numerical predictions made by a ROM are obtained by solving physics equations.
\end{remark}

\begin{remark}
{This paper focuses on the choice of the dissimilarity measure in the clustering task for the construction of the dictionary of ROMs, and on hyperparameters calibration. Our article~\cite{mca26010017} gives more details about the other important component of a dictionary-based ROM-net, namely the classifier designed for automatic model recommendation.}
\end{remark}

\begin{algorithm}[H]
    \caption{Training algorithm for dictionary-based ROM-nets}
    \label{alg:ROM-net_train}
  \begin{algorithmic}[1]
    \INPUT Points $\{ x_{i} \}_{1 \leq i\leq m}$ in the parameter space, number of clusters $K$, number $n_s$ of high-fidelity snapshots per cluster. 
    \OUTPUT Trained dictionary-based ROM-net made of $K$ ROMs and one classifier for automatic model recommendation.
    \STATE \textbf{Stage 1 (simplified simulation data generation):}
    \STATE Call the high-fidelity solver to compute solutions of a simplified version of the target physics problem.
    \STATE Simplified snapshots $\{ u_{i}^{LF} \}_{1 \leq i\leq m} := $ solutions for $\{ x_{i} \}_{1 \leq i\leq m}$.
    \STATE \textbf{Stage 2 (clustering based on simulation data)}
    \STATE Compute sine dissimilarities between simplified snapshots.
    \STATE Run a k-medoids clustering algorithm using sine dissimilarities as clustering distance, to get $K$ clusters.
    \STATE $\mathcal{K}_K (x_i) :=$ label/cluster of $u_{i}^{LF}$.
    \STATE \textbf{Stage 3 (snapshots selection)}
    \FOR{$k \in [\![ 1;K ]\!]$}
      \STATE Run a k-medoids clustering algorithm on the $k$-th cluster to get $n_s$ subclusters.
    \ENDFOR
    \STATE $\mathcal{I} := $ indices of the subclusters' medoids.
    \STATE Call the high-fidelity solver to compute solutions of the target physics problem for $\{ x_{i} \}_{i\in \mathcal{I}}$.
    \STATE High-fidelity snapshots $\{ u_{i}^{HF} \}_{i\in \mathcal{I}} := $ solutions for $\{ x_{i} \}_{i\in \mathcal{I}}$.
    \STATE \textbf{Stage 4 (ROM dictionary construction):}
    \FOR{$k \in [\![ 1;K ]\!]$}
      \STATE Build a local (hyper-)reduced-order model for the $k$-th cluster using its $n_s$ high-fidelity snapshots.
    \ENDFOR
    \STATE \textbf{Stage 5 (classification for automatic model recommendation)}
    \STATE Train a classifier $\mathcal{C}_K$ on the labeled dataset $\{ (x_{i}, \mathcal{K}_K (x_i) \}_{1 \leq i\leq m}$.
  \end{algorithmic}
\end{algorithm}

\section{Clustering background}
\label{ClusteringBackground}

\subsection{Representative-based clustering}

In representative-based clustering algorithms, each cluster is associated to a partitioning representative, \textit{i.e.} a reference point that well represents the cluster's members. Clusters representatives are useful in our case, because they can be used to select the high-fidelity snapshots. Representative-based algorithms generally define the clusters thanks to the Voronoi diagram generated by the representatives, which gives clusters with high cohesion.

\begin{definition}
[Representative-based clustering] Let us consider a finite set $\{ x_i \}_{1 \leq i \leq m}$ of elements of a topological space $\mathcal{T}$ endowed with a dissimilarity measure $\delta$. For a given integer $K\in [\![2;m]\!]$, \textit{representative-based clustering} consists in finding $K$ representatives $\{ \tilde{x}_k \}_{1 \leq k \leq K} \subset \mathcal{T}$ minimizing the objective function:\begin{equation}
\sum_{i=1}^{m} \underset{k \in [\![1;K]\!]}{\min} \delta(x_i, \tilde{x}_k)^2 .
\label{RepBasedClusCostFct}
\end{equation} The clusters $C_k$ are given by:\begin{equation}
C_k := \{ x_i \ | \ \delta(x_i , \tilde{x}_k) \leq \delta(x_i , \tilde{x}_l) \ \forall l \in [\![1;K]\!] \}.
\label{ClusterDefinition}
\end{equation}
\end{definition}

When the dissimilarity measure is the Euclidean distance, the optimal representatives are the clusters' means or centroids (see~\cite{aggarwal2015data}, p.162). This problem corresponds to k-means clustering~\cite{kmeans}, where the cost function in Equation~\eqref{RepBasedClusCostFct} corresponds to the within-cluster variance and is related to clusters inertia.

\subsection{K-medoids clustering}

In k-medoids, the representatives must be taken among the elements of the dataset. This restriction is particularly useful when functions of the training examples (such as mean and median) do not make sense or cannot be easily computed. It enables working with any type of data with any dissimilarity measure. The next definitions introduce the k-medoids optimization problem:

\begin{definition}
[Binary matrices] A \textit{binary matrix} is a matrix whose coefficients are either $0$ or $1$. The set of binary matrices of size $m\times n$ is denoted by $\mathcal{B}_{m,n}$.
\end{definition}

\begin{definition}
[K-medoids clustering] Let us consider a finite set $\{ x_i \}_{1 \leq i \leq m}$ of elements of a topological space $\mathcal{T}$ endowed with a dissimilarity measure $\delta$. For a given integer $K\in [\![2;m]\!]$, let us introduce the set $\mathcal{Z}_{m,K}$:\begin{equation}
\mathcal{Z}_{m,K} := \left\lbrace \g{Z} \in \mathcal{B}_{m,K} \ | \ \sum_{k=1}^K z_{ik} = 1 \ \forall i \in [\![1;m]\!] \ \textrm{and} \ \sum_{i=1}^m z_{ik} \geq 1 \ \forall k \in [\![1;K]\!] \right\rbrace .
\label{DefinitionSetZnK}
\end{equation} \textit{K-medoids clustering} consists in solving the following optimization problem:\begin{equation}
\g{Z}^{*} := \underset{\g{Z}\in \mathcal{Z}_{m,K}}{\argmin} \sum_{k=1}^{K} \sum_{i=1}^{m} z_{ik} \delta(x_i, \tilde{x}_k)^2 ,
\end{equation} where the medoids $\tilde{x}_k$ are given by: \begin{equation}
\tilde{x}_k := \underset{x_j \in \{ x_l \}_{1 \leq l \leq m}}{\argmin} \sum_{i=1}^{m} z_{ik} \delta(x_i, x_j)^2 .
\label{DefMedoid}
\end{equation}
\end{definition}

This formulation of the k-medoids problem has similarities with the k-means formulation proposed in~\cite{10.2307/2529785}. With this formulation, the definition of the clusters $C_k$ in Equation~\eqref{ClusterDefinition} is equivalent to:\begin{equation}
C_k = \{ x_i \ | \ z^{*}_{ik} = 1 \}.
\end{equation} Equation~\eqref{DefinitionSetZnK} defining the set $\mathcal{Z}_{m,K}$ ensures that each point is assigned to one single cluster, and that each cluster contains at least one element. Equation~\eqref{DefMedoid} defines the medoid of a cluster as its most central member. K-medoids is a combinatorial optimization problem, for which several heuristic approaches have been proposed to find a suboptimal solution at lower cost. The Partitioning Around Medoids (PAM~\cite{kMedoidsPAM} and Chap.~2 of~\cite{kaufman1990finding}) is the most known algorithm. It iteratively looks for the best swap between nonmedoid points and medoids. Clustering Large Applications (CLARA~\cite{kMedoidsCLARA} and Chap.~3 of~\cite{kaufman1990finding}) applies PAM on different subsamples to reduce the computational complexity of PAM. As explained in Section~11.2.1 of~\cite{aggarwal2013data}, both PAM and CLARA algorithms can be interpreted as graph-searching problems: PAM explores the entire graph of clustering solutions, while CLARA explores a subgraph only. Clustering Large Applications based on Randomized Sampling (CLARANS~\cite{kMedoidsCLARANS2, kMedoidsCLARANS}) only considers a sample of the neighbors of the current graph node at each iteration, which enables searching over the entire graph as in PAM but at lower cost. These three algorithms have been improved recently in~\cite{kMedoids2019} in terms of computational complexity. Apart from these approaches, a simple and fast k-medoids algorithm has been proposed in~\cite{kMedoidsParkJun2009} following the standard implementation of k-means, \textit{i.e.} alternating between a cluster assignment step and updating the medoids with Equation~\eqref{DefMedoid}. However, as explained in~\cite{kMedoids2019}, this algorithm does not explore as many configurations as PAM does. For all these algorithms, the dissimilarities $\delta(x_i, x_j)$ are precomputed before looking for clusters.

\section{Proposed local ROM approach}
\label{SectionPhyInformedClus}

\subsection{Motivations}

Instead of considering the absolute projection error when defining the Kolmogorov $N$-width, one can use the relative projection error, which leads to the following definition:

\begin{definition}
[Normalized Kolmogorov $N$-width] Let $N \in \mathbb{N}^*$. If $\mathcal{M}$ contains at least one nonzero element, the \textit{normalized Kolmogorov $N$-width} of the manifold $\mathcal{M}$ in the ambient Hilbert space $\mathcal{H}$ is defined by:\begin{equation}
\tilde{d}_{N}(\mathcal{M}) := \underset{\mathcal{H}_{N} \in \emph{\textrm{Gr}}(N,\mathcal{H})}{\inf} \  \underset{u \in \mathcal{M} \setminus \{ 0 \}}{\sup}  \ \underset{v \in \mathcal{H}_N}{\inf} \frac{|| u - v ||_{\mathcal{H}}}{|| u ||_{\mathcal{H}}}.
\end{equation}
\end{definition}

\noindent Let $\measuredangle_{\mathcal{H}} (u,v) \in [0;\pi /2]$ denote the angle between two nonzero elements $u$ and $v$ of $\mathcal{H}$:\begin{equation}
\measuredangle_{\mathcal{H}} (u,v) := \arccos \left( \frac{ | \langle u,v \rangle_{\mathcal{H}} |}{|| u ||_{\mathcal{H}} || v ||_{\mathcal{H}}} \right),
\label{AngleDefinition1}
\end{equation} and let $\measuredangle_{\mathcal{H}} \left(u,\mathcal{V}\right) \in [0;\pi /2]$ be the angle between $u$ and a subspace $\mathcal{V} \subset \mathcal{H}$:\begin{equation}
\measuredangle_{\mathcal{H}} \left(u,\mathcal{V}\right) := \underset{v \in \mathcal{V}}{\inf} \ \measuredangle_{\mathcal{H}} (u,v).
\end{equation} The normalized Kolmogorov $N$-width is related to the largest angle between elements of the solution manifold and the approximation space:

\begin{proper}
Let $N \in \mathbb{N}^*$, and suppose that $\mathcal{M}$ contains at least one nonzero element. Then:\begin{equation}
\tilde{d}_{N}(\mathcal{M}) = \underset{\mathcal{H}_{N} \in \emph{\textrm{Gr}}(N,\mathcal{H})}{\inf} \  \underset{u \in \mathcal{M} \setminus \{ 0 \}}{\sup}  \ \sin \measuredangle_{\mathcal{H}} \left(u,\mathcal{H}_N \right).
\label{KolmogorovWidthSineAngle}
\end{equation}
\label{LinkKolmogorovSineAngle}
\end{proper}

\begin{figure}[!h]
\centering
\includegraphics[scale=0.22]{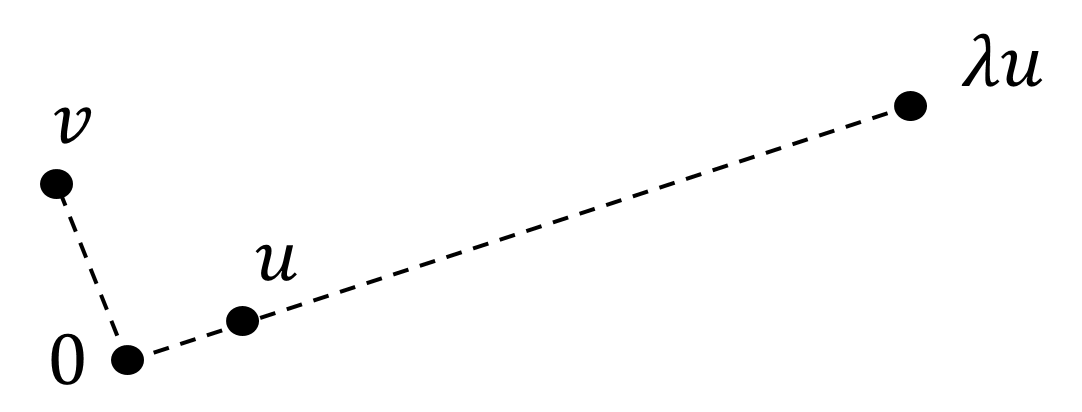}
\caption{Clustering with the Euclidean distance would assign $u$ and $v$ to the same cluster and $\lambda u$ to another, whereas $u$ and $\lambda u$ could be computed with the same 1D approximation space.}
\label{euclid_dist_vs_angle}
\end{figure}

{The proof of this property is given in Appendix A, with another property linking the normalized Kolmogorov width $\tilde{d}_{N}(\mathcal{M})$ with the absolute Kolmogorov width $d_{N}(\mathcal{M})$ via an inequality. As suggested by Equation~\eqref{KolmogorovWidthSineAngle} of Property~\ref{LinkKolmogorovSineAngle}, the dissimilarity measure should be defined as a function of the angle between elements of the solution manifold in order to focus on the shape of the fields $u\in\mathcal{M}$ rather than their intensities.} In this way, clustering would efficiently decrease projection errors by limiting the maximum angular deviations within clusters. The Euclidean distance $|| u-v ||_{\mathcal{H}}$ used in~\cite{localROB, localROB2, Amsallem2015localHROM, RyckelynckComputerVision, Grimberg2020} does not always ensure the reduction of projection errors. Indeed, the solution manifold can contain solutions that are relatively close in terms of the Euclidean distance but distributed in many different directions of the space $\mathcal{H}$. On the other hand, having a subset $\mathcal{M}_k$ with a large diameter in terms of the Euclidean distance is not a problem if it is embedded in a low-dimensional space, as indicated by Property~\ref{LinkKolmogorovSineAngle}. Let us suppose that the solution manifold contains two elements $u$ and $v$ having disjoint supports $\textrm{supp}(u)$ and $\textrm{supp}(v)$ and such that there exists a large real number $\lambda$ such that $\lambda u$ is still in the solution manifold (see Figure~\ref{euclid_dist_vs_angle}). The elements $u$ and $\lambda u$ are aligned in the same direction and could then be obtained with the same $1$-dimensional approximation space. However, if $\lambda$ is large enough, the distance $|| u-\lambda u ||_{\mathcal{H}}$ can be very large with respect to $|| u-v ||_{\mathcal{H}}$. In this case, it is possible to assign $u$ and $v$ to the same cluster while assigning $\lambda u$ to another, whereas $u$ and $\lambda u$ are aligned along a direction that is orthogonal to $v$. For these reasons, the Euclidean distance does not seem to be adapted, except if the number $K$ of clusters is large enough to get very local subsets $\mathcal{M}_k$ with restricted angular deviations.

A more natural and straightforward approach would consist in clustering the parameter space $\mathcal{X}$ to define the subsets $\mathcal{M}_k = u(\mathcal{X}_k)$ for each cluster $\mathcal{X}_k$. Note that the subsets $\mathcal{M}_k$ no longer form a partition of $\mathcal{M}$, although their union still equals to $\mathcal{M}$. This strategy may not be appropriate when $u$ is a nonlinear function of the parameters $x\in\mathcal{X}$. The physics of the underlying problem can also generate situations where small changes of the parameters in some directions of the parameter space totally modifies the shape of the solution in a nonlinear way, while large variations in other directions of the parameter space only imply linear variations. An example is given in~\cite{ROM-net}, where it is shown that clusters identified in the parameter space give subsets~$\mathcal{M}_k$ spreading all over the solution manifold~$\mathcal{M}$. To avoid this issue, it is preferable to apply a physics-informed clustering strategy by partitioning the solution manifold directly with an appropriate dissimilarity measure $\delta$.

\subsection{Proposed dissimilarity measure}

This section introduces the dissimilarity measure used in this paper for clustering and gives some of its properties. {It is important to stress that this dissimilarity is computed from the simplified snapshots $u^{LF}$ given by the simplified simulations. Hence, in this section, the notation $u^{LF} \in L^{2}(\Omega \times [0;t_f])$ represents a simplified snapshot.}

\begin{definition}
[Principal angles between subspaces] Let $\mathcal{V}_1$ and $\mathcal{V}_2$ be two subspaces of $L^2 (\Omega)$. The \textit{principal angles} or \textit{canonical angles} $\theta_k (\mathcal{V}_1 , \mathcal{V}_2 ) \in [0;\pi/2]$ between $\mathcal{V}_1$ and $\mathcal{V}_2$ are defined by:\begin{equation}
\forall k \in \mathbb{N}^{*}, \quad \theta_k (\mathcal{V}_1 , \mathcal{V}_2 ) := \measuredangle (v_{1}^k ,v_{2}^k),
\end{equation} where the angle $\measuredangle := \measuredangle_{L^{2}(\Omega)}$ is measured in $L^{2}(\Omega)$ (see Equation~\eqref{AngleDefinition1}), and where the vectors $v_{1}^k \in \mathcal{V}_1$ and $v_{2}^k \in \mathcal{V}_2$ are given by the following sequence of optimization problems:\begin{equation}
\left\lbrace \begin{array}{rcl}
(v_{1}^{1},v_{2}^{1}) & \in & \underset{(v_1, v_2) \in \mathcal{V}_1 \times \mathcal{V}_2}{\argmin} \ \measuredangle (v_1,v_2) \\
(v_{1}^{k+1},v_{2}^{k+1}) & \in & \argmin \ \left\lbrace \measuredangle (v_1,v_2) \ \vline \ v_j \in \mathcal{V}_{j} \cap \left( \emph{\textrm{span}} (\{ v_{j}^i \}_{1 \leq i \leq k}) \right)^{\perp}, \ j \in \{1;2\} \right\rbrace
\end{array}\right.
\end{equation} with the notation $\mathcal{V}^{\perp}$ denoting the orthogonal complement of $\mathcal{V}\subset L^{2}(\Omega)$ in $L^{2}(\Omega)$.
\end{definition}

In practice, when the spaces $\mathcal{V}_1$ and $\mathcal{V}_2$ are finite-dimensional, it can be shown (see Theorem~1 of~\cite{Bjorck73}) that the principal angles are given by:\begin{equation}
\forall k \in [\![ 1; \min(\dim(\mathcal{V}_1), \dim(\mathcal{V}_2)) ]\!], \quad \theta_{k}\left(\mathcal{V}_1, \mathcal{V}_2 \right) = \arccos \sigma_k ,
\end{equation} with $\sigma_1 \geq \sigma_2 \geq ... \geq \sigma_{\min(\dim(\mathcal{V}_1), \dim(\mathcal{V}_2))}$ being the singular values of the matrix $\g{C}(\mathcal{V}_1, \mathcal{V}_2)\in \mathbb{R}^{\dim(\mathcal{V}_1) \times \dim(\mathcal{V}_2)}$ defined by:\begin{equation}
C_{ij}(\mathcal{V}_1, \mathcal{V}_2) := \langle \psi_{i}^{(1)}, \psi_{j}^{(2)} \rangle_{L^{2}(\Omega)} ,
\end{equation} where the functions $\psi_{i}^{(1)}$ (resp. $\psi_{j}^{(2)}$) form an orthonormal basis of $\mathcal{V}_1$ (resp. $\mathcal{V}_2$). The vector $\bs{\theta}(\mathcal{V}_1 , \mathcal{V}_2)$ denotes the vector containing the principal angles between the spaces $\mathcal{V}_1$ and $\mathcal{V}_2$.

\begin{definition}
[$n$-dimensional elementary basis] Let $u \in L^{2}(\Omega \times [0;t_f])$ and $n \in [\![ 1 ; \mathcal{N} ]\!]$. The \textit{$n$-dimensional elementary basis} associated to $u$ is the orthonormal $n$-frame $\Psi_{n}(u) \in V(n,L^{2}(\Omega))$ obtained by solving the POD minimization problem given in Equation~\eqref{PODdef} with the Snapshot POD algorithm, using the trajectory of $u$ over time as a snapshot.
\label{ElemIdealBasis}
\end{definition}

\begin{definition}
[$n$-dimensional elementary approximation space] Let $u \in L^{2}(\Omega \times [0;t_f])$ and $n \in [\![ 1 ; \mathcal{N} ]\!]$. The \textit{$n$-dimensional elementary approximation space} $\mathcal{V}_{n}(u) \in \emph{\textrm{Gr}}(n,L^{2}(\Omega))$ is the subspace spanned by $\Psi_{n}(u)$.
\end{definition}

In Definition~\ref{ElemIdealBasis}, the POD basis $\Psi_{n}(u)$ is used for clustering only, it is not supposed to be used for numerical simulations since it is computed from simplified snapshots. Qualitatively, the subspace $\mathcal{V}_{n}(u)$ spanned by this POD basis is the best $n$-dimensional approximation space for the trajectory of $u(.,t)$ in $L^{2}(\Omega)$, that is to say:\begin{equation}
\mathcal{V}_{n}(u) = \underset{\mathcal{V}_{n} \in \textrm{Gr}(n,L^{2}(\Omega))}{\argmin} \ \int_{0}^{t_f} \underset{v \in \mathcal{V}_{n}}{\inf} \ || u(.,t) - v ||_{L^{2}(\Omega)}^2 \ dt .
\end{equation}

\begin{definition}
[Chordal distance between subspaces~\cite{chordalDist}, p.~140, Section~2] Let $\mathcal{H}$ be a Hilbert space, and $n, m$ be two integers with $n \leq m$. The chordal distance between subspaces $\mathcal{V}_1 \in \emph{\textrm{Gr}}(n,\mathcal{H})$ and $\mathcal{V}_2 \in \emph{\textrm{Gr}}(m,\mathcal{H})$ is defined by:\begin{equation}
d_{c}(\mathcal{V}_1 , \mathcal{V}_2) := || \sin \bs{\theta}(\mathcal{V}_1 , \mathcal{V}_2) ||_2 = \left( \sum_{k=1}^n \sin^2 \theta_k (\mathcal{V}_1, \mathcal{V}_2) \right)^{1/2} .
\end{equation}
\end{definition}

\begin{definition}
[Sine dissimilarity between functions] Given $n \in [\![ 1 ; \mathcal{N} ]\!]$, the sine dissimilarity $\tilde{\delta}_n $ between functions $u$ and $v$ in $L^{2}(\Omega \times [0;t_f])$ is defined by:\begin{equation}
\tilde{\delta}_n (u,v) := d_{c}(\mathcal{V}_{n}(u) , \mathcal{V}_{n}(v)) .
\end{equation}
\label{DefSineDissimilarity}
\end{definition}

Let us now recall the definition of the orthogonal projection $\pi_{\mathcal{H}_{n}}: L^{2}(\Omega) \rightarrow L^{2}(\Omega)$ on a $n$-dimensional subspace $\mathcal{H}_{n}$ of $L^{2}(\Omega)$, with an orthonormal basis $\{ \psi_k \}_{1\leq k \leq n}$:\begin{equation}
\forall u \in L^{2}(\Omega), \quad \pi_{\mathcal{H}_{n}}(u) = \sum_{k=1}^n \langle u, \psi_{k} \rangle_{L^{2}(\Omega)} \psi_{k} .
\label{DefOrthogProj}
\end{equation} {The following properties give interesting interpretations of the sine dissimilarity that motivate its use for the construction of dictionaries of local ROMs. The proofs of these properties are given in Appendix B.}

\begin{proper}
[Sine dissimilarity and $L^2$ projection errors] For all $n \in [\![ 1 ; \mathcal{N} ]\!]$, the sine dissimilarity is symmetric and satisfies:\begin{equation}
\forall (u,v) \in L^{2}(\Omega \times [0;t_f])^2, \quad \tilde{\delta}_n (u,v) = \left( \sum_{i=1}^n || \psi_{i}(u) - \pi_{\mathcal{V}_{n}(v)}(\psi_{i}(u)) ||_{L^{2}(\Omega)}^2 \right)^{1/2} ,
\label{firstFormulaDeltaTilde}
\end{equation} with $\pi_{\mathcal{V}_{n}(v)}$ denoting the orthogonal projection on $\mathcal{V}_{n}(v)$ and where the functions $\psi_{i}(u) \in L^{2}(\Omega)$ for $i \in [\![ 1 ; n ]\!]$ are the vectors of the elementary basis $\Psi_{n}(u)$.
\label{SineDissAndProjErrorProp}
\end{proper}

\begin{proper}
[Sine dissimilarity and Hilbert-Schmidt distance] For all $n \in [\![ 1 ; \mathcal{N} ]\!]$, for all $(u,v) \in L^{2}(\Omega \times [0;t_f])^2$, the sine dissimilarity satisfies:\begin{equation}
\tilde{\delta}_n (u,v) = \frac{1}{\sqrt{2}} || \pi_{\mathcal{V}_{n}(u)} - \pi_{\mathcal{V}_{n}(v)} ||_{HS(L^{2}(\Omega))} .
\label{symFormulaDeltaTilde}
\end{equation}
\label{PropHSdistance}
\end{proper}

{These properties show that the sine dissimilarity can be interpreted either in terms of projection errors, or as the Hilbert-Schmidt distance between the projections onto the elementary approximation spaces. The sine dissimilarity is therefore relevant for model order reduction methods using Galerkin projection for the computation of approximate solutions in low-dimensional approximation spaces. In addition to the proofs of the two aforementioned properties, two other mathematical results on the sine dissimilarity are given in Appendix B: we show that this dissimilarity is a pseudometric on $L^{2}(\Omega \times [0;t_f])$, and that it is asymptotically equivalent to the Grassmann dissimilarity used in our previous paper~\cite{ROM-net} for small angles. Finally, the next definition introduces the ROM-oriented dissimilarity between parameters as the sine dissimilarity between the corresponding solutions:}

\begin{definition}
[ROM-oriented dissimilarity between parameters] Given $n \in [\![ 1 ; \mathcal{N} ]\!]$, the \textit{ROM-oriented dissimilarity} between parameters $x$ and $x'$ in $\mathcal{X}$ is defined by:\begin{equation}
\delta_{n}\left(x, x' \right) := \tilde{\delta}_n \left(u^{LF} \left(x\right), u^{LF} \left(x'\right) \right) ,
\label{DefROMdissimilarity}
\end{equation} where {$u^{LF}:\mathcal{X} \rightarrow L^{2}(\Omega \times [0;t_f])$} is either the primal variable (\textit{i.e.} the solution of the physics problem) or a dual variable (\textit{i.e.} an internal variable) defining a quantity of interest.
\label{DefDissimilarity}
\end{definition}

It is recalled that this dissimilarity is computed from simplified snapshots. Property~\ref{pseudometric} {given in Appendix B} implies that the ROM-oriented dissimilarity is a pseudometric on $\mathcal{X}$. Several variants of this dissimilarity can be obtained according to the definition of the variable {$u^{LF}$}. Using the primal variable should improve the quality of the POD-Galerkin approximation, since the data would be clustered according to the angles between the subspaces spanned by the trajectories of the primal solution. This would give a \textit{method-oriented} dissimilarity, that is, a dissimilarity favoring the accuracy of the numerical method (namely model order reduction) used for numerical simulations. Using a dual variable instead would improve the quality of the Gappy-POD~\cite{GappyPOD} reconstruction for the quantity of interest when hyper-reduction is used. This would define a \textit{goal-oriented} method favoring the accuracy of numerical predictions of a quantity of interest. Of course, one could mix both strategies by taking a weighted average of these two variants of the ROM-oriented dissimilarity.

\subsection{Choice of the clustering method}

As an unsupervised learning task, clustering has no indisputable evaluation criterion. This is the reason why there is no hierarchy in the large variety of clustering algorithms. The algorithm must be selected according to the purpose. For model order reduction purposes, we have seen that the Kolmogorov $N$-width relates the physics problem's reducibility to projection errors on the approximation space, which makes the projection error a good candidate for an evaluation criterion:

\begin{definition}
[Relative projection error] Let $u \in L^{2}(\Omega)$ be a nonzero square-integrable function, and $\Psi = \{ \psi_k \}_{1 \leq 1 \leq N} \in V(N,L^{2}(\Omega))$ be an orthonormal reduced-order basis of dimension $N\in\mathbb{N}^*$ in $L^{2}(\Omega)$. The \textit{relative projection error} $\eta(u,\Psi)$ of $u$ on $\emph{\textrm{span}}(\Psi)$ is given by:\begin{equation}
\eta(u,\Psi) := \frac{|| u - \sum_{k=1}^N \langle u, \psi_k \rangle_{L^{2}(\Omega)} \psi_k ||_{L^{2}(\Omega)}}{|| u ||_{L^{2}(\Omega)}} .
\label{ProjectionError}
\end{equation}
\end{definition}

\begin{remark}
{The relative projection error does not depend on the choice of the orthonormal basis used to represent the subspace $\emph{\textrm{span}}(\Psi_N)$. Therefore, the notations $\eta(u,\Psi_N)$ and $\eta(u,\emph{\textrm{span}}(\Psi_N))$ can be used interchangeably.}
\end{remark}

As shown in Equation~\eqref{KolmogorovWidthSineAngle} in Property~\ref{LinkKolmogorovSineAngle}, Kolmogorov widths can be decreased by limiting the angular deviation within the clusters. Having defined a dissimilarity measure $\delta_n$ based on angles in Definitions~\ref{DefSineDissimilarity} and~\ref{DefDissimilarity}, one must look for compact-shaped clusters in terms of the dissimilarity $\delta_n$. {Therefore, we use PAM k-medoids clustering algorithm to reduce the intra-cluster maximum angular deviations as much as possible. Our physics-informed clustering method consists in running simplified simulations and applying PAM to simulation data using the ROM-oriented dissimilarity.}

{When computed in $L^2 (\Omega)$ and therefore with $n=1$, the sine dissimilarity has a simple formula:\begin{equation}
\forall (u,v) \in L^2 (\Omega)^2, \quad \tilde{\delta}_1 (u,v)_{L^2 (\Omega)} :=  \sin \measuredangle_{L^2 (\Omega)} \left(u , v \right) = \sqrt{1 - \frac{ \langle u , v \rangle_{L^2 (\Omega)}^2}{|| u ||_{L^2 (\Omega)}^2 || v ||_{L^2 (\Omega)}^2} },
\end{equation} which gives a direct link with the relative projection error:\begin{equation}
\tilde{\delta}_1 (u,v)_{L^2 (\Omega)} = \eta(u, \textrm{span}(\{v\})).
\end{equation} This formula will be used for the computation of the dissimilarity in the applications given at the end of this paper. In this setting, we showed in~\cite{OptimalPartitions-DanielCasenaveAkkariRyckelynck} the following property motivating the choice of k-medoids clustering:}

\begin{proper} [Optimality of k-medoids clustering]
{The partitions $\mathcal{M}_k$ of $\mathcal{M}$ minimizing the k-medoids cost function with dissimilarity $\tilde{\delta}_1 (.,.)_{L^2 (\Omega)}$ are exactly the minimizers of a discretized version of the following cost function:\begin{equation}
\sum_{k=1}^{K} \Prob(u \in \mathcal{M}_k) \   \check{d}_{1} (p_{U | u \in \mathcal{M}_k})^2,
\end{equation} where $\check{d}_{N}$ is a variant of the Kolmogorov width obtained by replacing the worst-case error by the mean squared error as in~\cite{bachmayr2017kolmogorov}:\begin{equation}
\begin{array}{rcl}
\check{d}_{N}(p_U) & :=  & \left( \underset{\mathcal{H}_{N} \in \textrm{Gr}(N, L^2 (\Omega))}{\inf} \  \E_{U \sim p_U} \left[ \eta \left( U, \mathcal{H}_N \right)^2 \right] \right)^{1/2} \\
 & = & \left( \underset{\mathcal{H}_{N} \in \textrm{Gr}(N, L^2 (\Omega))}{\inf} \  \E_{X \sim p_X} \left[ \eta \left( u(X), \mathcal{H}_N \right)^2 \right] \right)^{1/2},
\end{array}
\end{equation} with $p_X$ denoting a probability density function in the parameter space and $p_U$ being the resulting probability density function in the solution space.}
\end{proper}

\subsection{Automatic snapshots selection}

Once clusters have been identified within the dataset, one must select relevant points for which the entire high-fidelity simulation will be run to provide high-fidelity snapshots for the construction of the local ROMs. For each cluster, the high-fidelity snapshots must be well distributed and representative of the cluster's members. When one wants to use only one snapshot per cluster, then the clusters' medoids are good candidates. For more than one snapshot per cluster, a second k-medoids cluster analysis can be conducted within each cluster, with $n_s$ subclusters where $n_s$ is the desired number of high-fidelity snapshots per cluster, using the same dissimilarity measure as for the first clustering. High-fidelity snapshots can then be computed for the subclusters' medoids. This method corresponds to a two-stage hierarchical k-medoids clustering.

\section{ROM-net's efficiency criterion and hyperparameters tuning}
\label{SectionROMnetEffCrit}

\subsection{Gain with respect to a global reduced-order model}

A dictionary-based ROM-net~\cite{ROM-net} is made of a dictionary of $K$ local ROMs and a classifier $\mathcal{C}_K$ which automatically selects the best model from the dictionary for a given point in the parameter space without computing any physics-informed dissimilarity, see Figure~\ref{ROM-net_online}. The real classifier $\mathcal{C}_K$ enables bypassing the simplified simulation that is required to evaluate the perfect classifier $\mathcal{K}_K$, see Figure~\ref{ROM-net_offline}. In this section, it is assumed that all the dictionary's ROMs have the same number of modes, denoted by $N$, and have been built from the same number of high-fidelity snapshots, denoted by $n_s$. A dictionary of $K$ ROBs with $N$ modes and $n_s$ high-fidelity snapshots per basis is denoted by $\{ \Psi_{k}^{(K,N,n_s)} \}_{k\in[\![1;K]\!]}$. The objective of this section is to define a practical method for the calibration of the hyperparameters $K$, $N$ and $n_s$, based on an evaluation criterion quantifying the ROM-net's profitability with respect to a single global ROM. This criterion must be computable very early in the ROM-net's training phase, right after the physics-informed clustering procedure in Figure~\ref{ROM-net_offline} and before the computation of high-fidelity snapshots, the construction of the ROMs, and the classifier's training phase. Therefore, the local ROBs $\{ \Psi_{k}^{(K,N,n_s)} \}_{k\in[\![1;K]\!]}$ used in the evaluation criterion are simply built from $K n_s$ simplified snapshots selected by the two-stage hierarchical k-medoids clustering, instead of the corresponding high-fidelity snapshots that will be computed afterwards. Their performances are compared with the performance of a global ROB $\Psi_{g}^{(1,N, K n_s)}$ containing $N$ modes and inferred from the same $K n_s$ snapshots as $\{ \Psi_{k}^{(K,N,n_s)} \}_{k\in[\![1;K]\!]}$. This global basis thus benefits from the physics-informed clustering procedure for the selection of its snapshots. The ROBs are related to the function {$u^{LF}(X) \in L^{2}(\Omega \times [0; t_f ])$} parametrized by the random variable $X$ representing the current point in the parameter space. The following definition introduces the \textit{gain} used in our evaluation criterion: 

\begin{definition}
[Gain] Given integers $K > 1$, $N>0$ and $n_s >0$ and a classifier $\mathcal{F}_K :\mathcal{X}\rightarrow [\![1;K]\!]$, the \textit{gain} is defined by:\begin{equation}
G(X;K,N,n_s,\mathcal{F}_K) = \frac{\eta\left(u^{LF}(X),\Psi_{g}^{(1,N,K n_s)}\right)}{\eta\left(u^{LF}(X),\Psi_{\mathcal{F}_{K}(X)}^{(K,N,n_s)}\right)} ,
\end{equation} where {$u^{LF}(X)$} results from a simplified simulation. For $K=1$, the gain equals to $1$.
\label{GainDefinition}
\end{definition}

\begin{remark}
In Definition~\ref{GainDefinition}, the primal variable $u$ can be replaced by a quantity of interest, depending on the choice made for the definition of the ROM-oriented dissimilarity.
\end{remark}

As a function of $X$, the gain can be seen as a random variable parametrized by the hyperparameters $K$, $N$ and $n_s$ and the classifier. The notations $G_{\mathcal{C}}(K,N,n_s)$ and $G_{\mathcal{K}}(K,N,n_s)$ denote $G(X;K,N,n_s,\mathcal{C}_K)$ and $G(X;K,N,n_s,\mathcal{K}_K)$ respectively. Right after the physics-informed clustering procedure, the user cannot evaluate the gain $G_{\mathcal{C}}(K,N,n_s)$ since the real classifier $\mathcal{C}_K$ has not been trained yet. However, the clusters implicitly define the perfect classifier $\mathcal{K}_K$ and thus the user has access to values of the gain $G_{\mathcal{K}}(K,N,n_s)$. In the next property, the following assumption is made:
\begin{itemize}
\item[{[A1]}] The gain $G_{\mathcal{K}}(K,N,n_s)$ is assumed to be deterministic, which means that it is no longer a random variable but rather a deterministic function of the hyperparameters $K$, $N$ and $n_s$. In other words, when the right cluster is chosen, the gain does not depend on $X$.
\end{itemize}

\begin{proper}
[Gain decomposition] Under assumption [A1]:\begin{equation}
\E [G_{\mathcal{C}}(K,N,n_s)] = p \ G_{\mathcal{K}}(K,N,n_s) + (1-p) \ E(K,N,n_s) ,
\end{equation} where $p = \Prob(\mathcal{C}_K = \mathcal{K}_K)$ is the classification accuracy and $E(K,N,n_s)$ given by:\begin{equation}
E(K,N,n_s) := \E [G_{\mathcal{C}}(K,N,n_s) \ \vline \ \mathcal{C}_K \neq \mathcal{K}_K] 
\end{equation} is the conditional expectation of the gain $G_{\mathcal{C}}(K,N,n_s)$ when selecting the wrong ROB.
\label{GainDecomp}
\end{proper}

\begin{proof}
The expected gain $\E [G_{\mathcal{C}}(K,N,n_s)]$ satisfies:\begin{equation}
\E [G_{\mathcal{C}}(K,N,n_s)] = p \ \E [G_{\mathcal{C}}(K,N,n_s) \ \vline \ \mathcal{C}_K = \mathcal{K}_K] + (1-p) \ E(K,N,n_s) .
\label{FirstEqGainProp}
\end{equation} If $G_{\mathcal{K}}(K,N,n_s)$ is constant for fixed hyperparameters $(K,N,n_s)$, then:\begin{equation}
G_{\mathcal{K}}(K,N,n_s) = \E [G_{\mathcal{K}}(K,N,n_s) \ \vline \ \mathcal{C}_K = \mathcal{K}_K] = \E [G_{\mathcal{C}}(K,N,n_s) \ \vline \ \mathcal{C}_K = \mathcal{K}_K] ,
\end{equation} because the gains $G_{\mathcal{K}}(K,N,n_s)$ and $G_{\mathcal{C}}(K,N,n_s)$ return the same values when the real classifier $\mathcal{C}_K$ selects the right ROB. Replacing $\E [G_{\mathcal{C}}(K,N,n_s) \ \vline \ \mathcal{C}_K = \mathcal{K}_K]$ by $G_{\mathcal{K}}(K,N,n_s)$ in Equation~\eqref{FirstEqGainProp} ends this proof.
\end{proof}

\noindent Two additional assumptions are made in what follows:\begin{itemize}
\item[{[A2]}] The classification accuracy $p$ is modeled as a  decreasing function of the number of clusters $K$ defined on a finite interval $[\![ 1;K_{\textrm{max}}]\!]$. Indeed, for a fixed number of training examples, increasing the number of classes $K$ makes the classification task more complicated. When the number of classes is too large in comparison with the number of training data, the classifier hardly improves the performance of a random guess classifier.
\item[{[A3]}] The conditional expectation $E(K,N,n_s)$ is constant, meaning that the expected gain when choosing the wrong ROB does not depend on the hyperparameters $K$, $N$ and $n_s$. For all $K,N,n_s$:\begin{equation}
\begin{array}{rcl}
E := E(K,N,n_s) & = & \E [G_{\mathcal{C}}(K,N,n_s) \ \vline \ \mathcal{C}_K \neq \mathcal{K}_K] \\
& \leq & \E [G_{\mathcal{K}}(K,N,n_s) \ \vline \ \mathcal{C}_K \neq \mathcal{K}_K] = G_{\mathcal{K}}(K,N,n_s) ,
\end{array}
\end{equation} so $E \leq G_{\mathcal{K}}(1,N,n_s) = 1$ in particular. In the application presented in the last section of this paper, we take $E = 0.75$. 
\end{itemize}

\noindent The next definition introduces the concept of \textit{real profitability} for a dictionary-based ROM-net:

\begin{definition}
[Real ROM-net profitability] Given integers $K > 1$, $N>0$ and $n_s >0$, a dictionary-based ROM-net with classifier $\mathcal{C}_K$ and ROM dictionary $\{ \Psi_{k}^{(K,N,n_s)} \}_{k\in[\![1;K]\!]}$ is profitable with a real profit $G_{r}^* \in \mathbb{R}_+$ if its expected gain satisfies $\E [G_{\mathcal{C}}(K,N,n_s)] \geq G_{r}^*$.
\end{definition}

This means that, on average, projection errors made by a global ROB are $G_{r}^*$ times larger than those made by the ROM-net, even when classification errors are taken into account. However, the ROM-net profitability cannot be evaluated a priori on $\E [G_{\mathcal{C}}(K,N,n_s)]$, since the real classifier has not been trained yet and the dictionary of ROBs $\{ \Psi_{k}^{(K,N,n_s)} \}_{k\in[\![1;K]\!]}$ inferred from high-fidelity snapshots have not been computed yet, see Figure~\ref{ROM-net_offline}. For these reasons, the following definition introduces the concept of \textit{perfect profitability}:

\begin{definition}
[Perfect ROM-net profitability] Given integers $K > 1$, $N>0$ and $n_s >0$, a dictionary-based ROM-net with perfect classifier $\mathcal{K}_K$ and ROM dictionary $\{ \Psi_{k}^{(K,N,n_s)} \}_{k\in[\![1;K]\!]}$ is perfectly profitable with a perfect profit $G_{p}^* \in \mathbb{R}_+$ if $\E [G_{\mathcal{K}}(K,N,n_s)] \geq G_{p}^*$.
\end{definition}

\begin{proper}
Let $G_{r}^* \in \mathbb{R}_+$. Let us consider a dictionary-based ROM-net with hyperparameters $K$, $N$, $n_s$. Under assumptions [A1], [A2] and [A3], the dictionary-based ROM-net is profitable with real profit $G_{r}^*$ if and only if it is perfectly profitable with the following perfect profit:\begin{equation}
G_{p}^* (G_{r}^*) = \frac{G_{r}^* - (1-p(K))E}{p(K)} .
\end{equation}
\end{proper}

\begin{proof}
It is a direct consequence of the gain decomposition property (Property~\ref{GainDecomp}).
\end{proof}

When the gains are computed with the results of the simplified simulations and with ROBs inferred from simplified snapshots, the dictionary-based ROM-net is said to be \textit{a priori profitable} with real profit $G_{r}^* > 1$ if:\begin{equation}
\E [G_{\mathcal{K}}(K,N,n_s)] \geq \frac{G_{r}^* - (1-p(K))E}{p(K)} .
\label{ProfitabilityCriterion}
\end{equation} The a priori profitability can be assessed early in the ROM-net training phase, right after the physics-informed clustering procedure.

\subsection{Practical method}

The number of clusters $K$, the number of POD modes $N$ and the number of snapshots per cluster $n_s$ are three important hyperparameters when building a dictionary-based ROM-net. Choosing a good number of clusters $K$ may be particularly difficult. The optimal value of $K$ is related to the nonlinearity of the solution manifold: the more curved the solution manifold $\mathcal{M}$ is, the greater $K$ must be to cover $\mathcal{M}$ with several subspaces. It also depends on the number of POD modes $N$: very fast simulations would require $N$ to be small, which would increase the number of local bases required to cover the solution manifold. Last but not least, $K$ also has an influence on the accuracy of the ROM-net's classifier. In a classification problem, increasing the number of classes while keeping the size of the training set constant makes the learning task tougher. Hence, the performance of a dictionary-based ROM-net does not monotonically increase with $K$ since its classifier may choose the wrong model, leading to inaccurate numerical predictions.

\noindent The hyperparameters $K$, $N$ and $n_s$ must satisfy the following requirements:
\begin{itemize}
\item[{[R1]}] \textbf{Limited computational resources:} the total number of high-fidelity snapshots $K n_s$ is limited by the maximum allowable budget in terms of high-fidelity simulations of the entire physics problem.
\item[{[R2]}] \textbf{Speed-up factor requirements:} to effectively reduce the computational cost of high-fidelity simulations, the number $N$ of POD modes per local ROB must not exceed $\mathcal{N}^{1/3}$.
\item[{[R3]}] \textbf{Accuracy requirements:} the mean projection error must be lower than a user-defined threshold $\eta^*$:\begin{equation}
\E [\eta\left(u^{LF}(X),\Psi_{\mathcal{K}_{K}(X)}^{(K,N,n_s)}\right)] \leq \eta^*  .
\label{AccuracyCriterion}
\end{equation}
\item[{[R4]}] \textbf{Gain requirements:} given a user-defined threshold $G_{r}^* > 1$, Equation~\eqref{ProfitabilityCriterion} for the ROM-net a priori profitability must be satisfied to ensure that the local bases give better performances than a single global ROB.
\end{itemize}

\begin{remark}
[Concerning requirement [R2]] After the Galerkin projection of the governing equations onto a ROB made of $N$ modes, the linear system to be solved at each iteration of the Newton-Raphson is full and thus has a complexity of $O(N^\alpha)$ with $2 \leq \alpha \leq 3$, which must be compared with the complexity $O(\mathcal{N}^\beta)$ of the sparse linear system obtained with the finite-element method, with $1 \leq \beta \leq 2$. The worst case is obtained for $\alpha=3$ and $\beta=1$, which gives an upper bound in the order of $\mathcal{N}^{1/3}$ for $N$.
\end{remark}

\noindent Given these constraints, we introduce the definition of \textit{hyperparameters admissible set}:

\begin{definition}
[Hyperparameters admissible set] The \textit{hyperparameters admissible set} is defined by:\begin{equation}
\mathcal{A} = \{ (K,N,n_s) \in (\mathbb{N}^*)^3 \ \vline  \quad \emph{\textrm{[R1], [R2], [R3], [R4] are satisfied.}} \} .
\end{equation}
\end{definition}

\noindent This definition gives a practical method for the ROM-net profitability analysis and hyperparameters tuning. The hyperparameters admissible set can be identified using simplified snapshots right after the clustering step in the training phase, see Figure~\ref{ROM-net_offline}. If the hyperparameters admissible set is empty, then it is not worth continuing the training phase of the dictionary-based ROM-net given the user-defined thresholds $\eta^*$ and $G_{r}^*$ and the maximum number of high-fidelity snapshots $n_{\textrm{snapshots}}^{\textrm{max}}$. The user can either build a global ROB using the physics-informed clustering results to identify snapshots, or weaken some of the requirements [R1] to [R4]. The time spent for simplified simulations is not wasted: the user can justify the choice of using a global ROB, and can benefit from these simulations for high-fidelity snapshots selection. On the contrary, if the hyperparameters admissible set is not empty, then there is a benefit in using a dictionary-based ROM-net. The choice of the best hyperparameters configuration among the admissible ones depends on the user's priorities. However, given the cost of the entire training phase, a ROM-net is generally used for applications where the number of test simulations is very high, \textit{e.g.} for parameter optimization or uncertainty quantification. In this case, once accuracy and gain requirements are met, one should take the smallest number of POD modes $N$ to get the highest possible speed-up factor. Among the admissible configurations with the smallest number of modes, it is recommended to choose the value of $K$ minimizing the mean projection error, to get the most accurate dictionary among the fastest admissible ones. The number of high-fidelity snapshots $n_s$ per cluster must be fixed accordingly so that the total number $K n_s$ of high-fidelity snapshots remains lower than $n_{\textrm{snapshots}}^{\textrm{max}}$.

\begin{remark}
Choosing the smallest possible number of modes $N$ generally implies choosing larger values for $K$, which usually decreases the performance of the ROM-net's classifier for automatic model recommendation. When interesting values for $K$ are rather large (say greater than $8$), one can artificially improve the classifier's accuracy by running several reduced simulations in parallel with the models having the highest membership probabilities. An error estimator could then be used to determine which reduced simulation is the most accurate, as proposed in~\cite{ryckelynck:hal-02559710}. Such a strategy increases the number of simulations to be run in the exploitation phase, but would enable working with large $K$'s and thus small $N$'s, lessening the computational complexity of online reduced simulations. In addition, when the number of training examples is not large enough compared to the number of clusters $K$ for the classification task, the data augmentation algorithm presented in~\cite{mca26010017} for the classification of numerical simulations can be applied to reduce the risk of overfitting.
\end{remark}

\section{Numerical applications}
\label{SectionApplication}

\subsection{1D steady heat equation}

\subsubsection{Problem description}

\noindent Let us consider the following ordinary differential equation:\begin{equation}
\left\lbrace \begin{array}{rclc}
-\left( \lambda u' \right)'(\xi) & = & s(\xi) & \quad \forall \xi \in [0;L] \\
u(0) & = & u_0 & \\
u(L) & = & u_0 &
\end{array}\right.
\label{PDEexample1}
\end{equation} where $\lambda\in L^2 ([0;L])$, $s\in L^2 ([0;L])$, $u_0 \in \mathbb{R}$ and $u-u_0\in H^{1}_{0} ([0;L])$. This equation describes the thermal behavior of an heterogeneous continuous medium of length $L$ with thermal conductivity $\lambda(\xi)$ and temperature $u(\xi)$, in the presence of a heat source $s(\xi)$. We are interested in the behavior of the solution $u$ under variying source terms and conductivity functions. The conductivity function $\lambda$ is defined by:\begin{equation}
\lambda (\xi) = \lambda_1 \mathds{1}_{\{\epsilon(\zeta) L \leq \xi \leq (\epsilon(\zeta) + \zeta)L \}} + \lambda_2 (\mathds{1}_{\{\xi < \epsilon(\zeta) L \}} + \mathds{1}_{\{\xi > (\epsilon(\zeta) + \zeta)L\}}) ,
\label{RandomConductivity}
\end{equation} with $\lambda_2 = 1000 \lambda_1 \in \mathbb{R}_{+}^{*}$ and with $\zeta$ being a random variable following the uniform distribution $\mathcal{U}(0.1,0.5)$. The random variable $\epsilon(\zeta)$ follows the uniform distribution $\mathcal{U}(0,1-\zeta)$. The source term~$s$ is modeled by a zero-mean Gaussian process with an exponential covariance function. The problem described by Equation~\eqref{PDEexample1} is therefore parametrized by the heat source distribution $s$ and the microstructural parameters $\zeta$ and $\epsilon(\zeta)$. The weak formulation of Equation~\eqref{PDEexample1} reads:\begin{equation}
\int_{0}^{L} \lambda(\xi) v'(\xi) u'(\xi) d\xi = \int_{0}^{L} s(\xi)v(\xi)d\xi \quad \forall v \in H^{1}_{0} ([0;L]) .
\end{equation} The interval $[0;L]$ is discretized into $\mathcal{N}-1 = 1999$ subdivisions of length $h=L/(\mathcal{N}-1)$. The vertices $\{\xi_i = ih \}_{0\leq i \leq \mathcal{N}-1}$ define a finite-element mesh whose P1 shape functions are denoted by $\{\phi_i \}_{1\leq i \leq \mathcal{N}-2}$. The shape functions $\phi_0$ and $\phi_{\mathcal{N}-1}$ are not used because of the Dirichlet boundary conditions. The finite-element method computes a high-fidelity approximate solution $u-u_0$ in the space $\textrm{span}\left( \{\phi_i \}_{1\leq i \leq \mathcal{N}-2} \right)$, whose coordinates are stored in a vector $\g{q}\in\mathbb{R}^{\mathcal{N}-2}$. This vector is the solution of the following linear system:\begin{equation}
\g{K}\g{q} = \g{f},
\label{Kqf}
\end{equation} with $\g{K}\in\mathbb{R}^{(\mathcal{N}-2)\times (\mathcal{N}-2)}$ given by:\begin{equation}
K_{ij} = \int_{0}^{L} \lambda(\xi) \phi_i '(\xi) \phi_j '(\xi) d\xi \quad \forall (i,j) \in [\![1; \mathcal{N}-2 ]\!] ,
\end{equation} and $\g{f}\in\mathbb{R}^{\mathcal{N}-2}$ given by: \begin{equation}
f_i = \int_{0}^{L} s(\xi)\phi_i (\xi)d\xi  \quad \forall (i,j) \in [\![1; \mathcal{N}-2 ]\!] .
\end{equation} 

\begin{figure}[!h]
\centering
\includegraphics[scale=0.5]{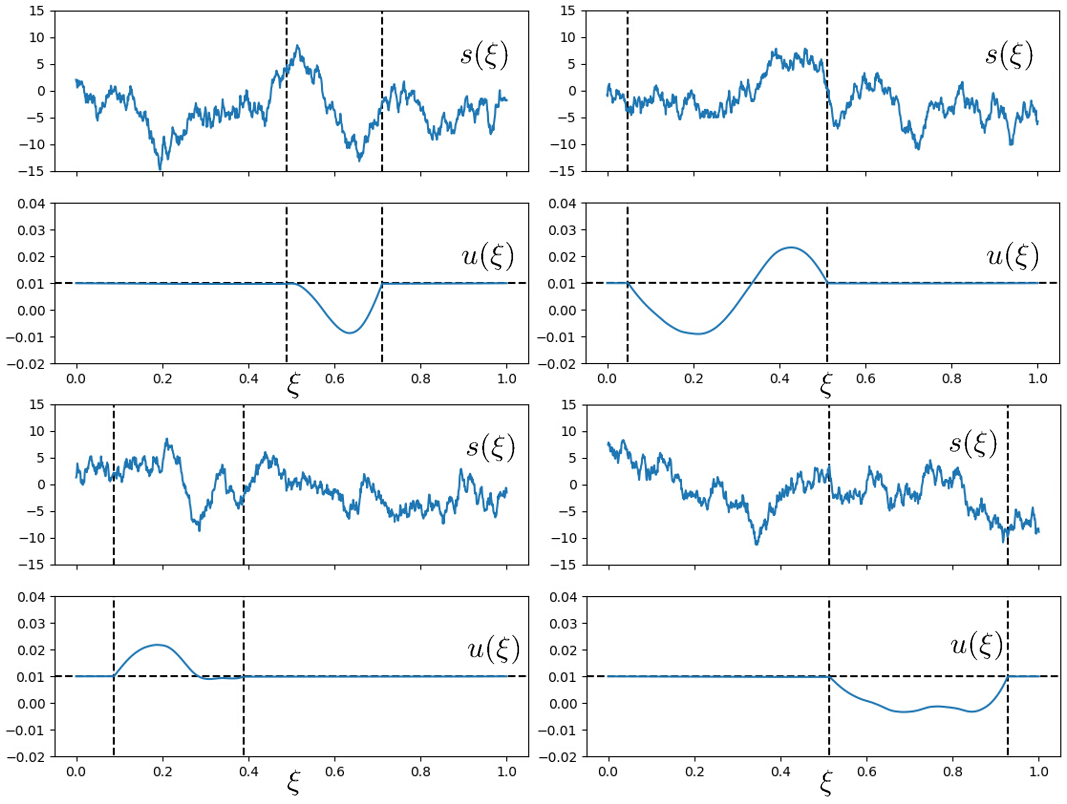}
\caption{Examples of solutions $u(\xi)$ for different source terms $s(\xi)$. The vertical dashed lines indicate the locations of the interfaces between the constituents of the bimaterial. Between the two interfaces, the thermal conductivity is $\lambda_1$. Outside of this interval, the thermal conductivity is $\lambda_2 = 1000 \lambda_1$.}
\label{solutions}
\end{figure}

\begin{figure}[!h]
\centering
\includegraphics[scale=0.39]{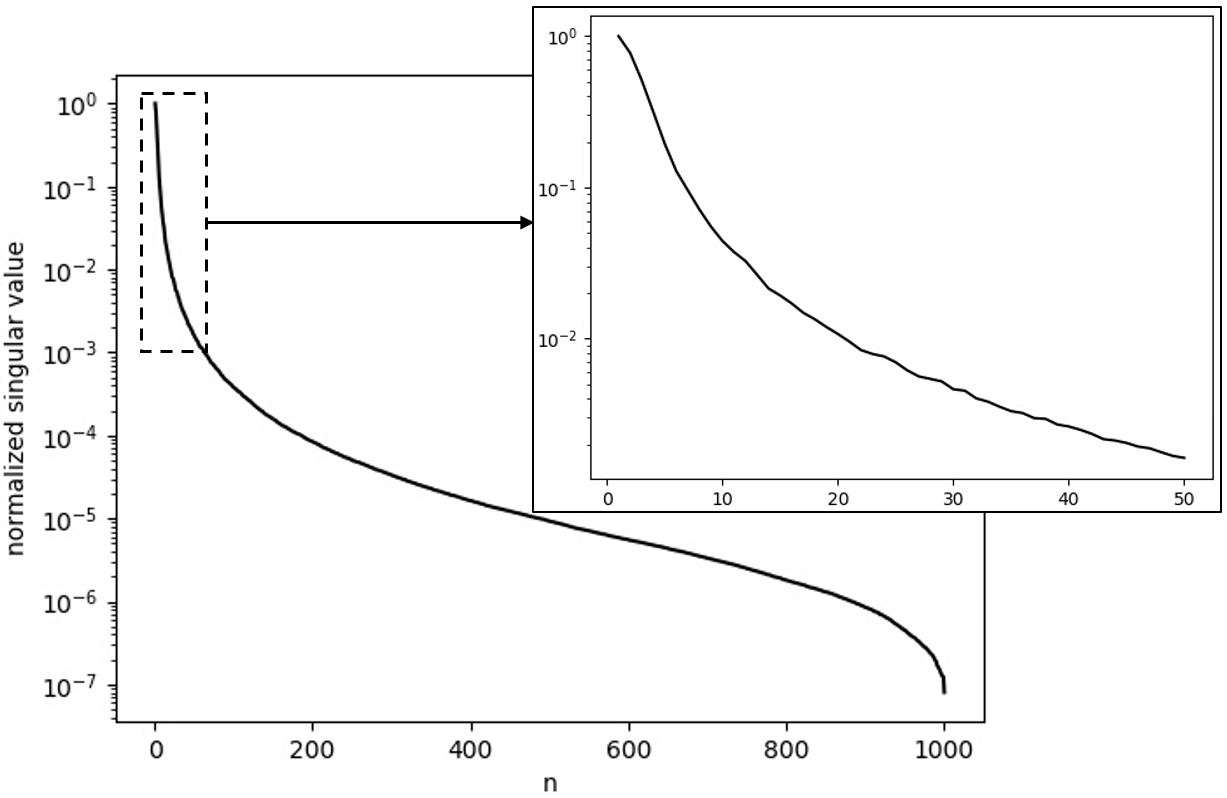}
\caption{Decay of the singular values obtained by singular value decomposition on the matrix containing all the examples in the dataset.}
\label{singValues}
\end{figure}

A dataset of $1000$ realizations of the random source term and microstructural parameters is generated. For each example in the dataset, the finite-element solution $\g{q}$ is computed with a Python routine. Figure~\ref{solutions} shows the solution's behavior for different configurations. One can observe that the solution is not affected by the source term in high-conductivity regions. Figure~\ref{singValues} gives the singular values of the matrix containing the $1000$ solutions. It can be observed that the decay of the singular values is rather slow for a 1D problem, meaning that this problem is non-reducible and that a dictionary of local ROBs may be required. The database is splitted into two subsets: a training set and a test set, both containing $500$ examples. The training set is used to identify clusters and build the ROBs, while the test set is used for evaluation purposes.

\begin{remark}
In the training phase of a dictionary-based ROM-net for time-dependent physics problems, the simplified problem that is simulated to provide data for the clustering procedure generally corresponds to a few time steps of the target problem. In this example, Equation~\eqref{PDEexample1} does not define a time-dependent problem. In this case, the simplified problem can be defined as the target problem solved on a coarse finite-element mesh.
\end{remark}

\subsubsection{Hyperparameters calibration}

This section deals with the calibration of the hyperparameters $(K,N,n_s)$. K-medoids is applied on the training set with the ROM-oriented dissimilarity $\delta_n$ introduced in Equation~\eqref{DefROMdissimilarity}. Since Equation~\eqref{PDEexample1} is time-independent, one must take $n=1$ for the ROM-oriented dissimilarity. We simply use the notation $\delta$ instead of $\delta_1$ for the ROM-oriented dissimilarity obtained by computing the sine dissimilarity in the solution space. The ROBs $\Psi_{1}(u)$ and $\Psi_{1}(v)$ used to calculate $\delta (u,v)$ are obtained by normalizing the solutions $u$ and $v$. For clustering, we use our own implementation of PAM~\cite{kMedoidsPAM, kaufman1990finding} k-medoids algorithm, with multiple random initializations for the medoids.

The physics problem considered in this section gives only one field $u$ per set of parameters. Therefore, the number of POD modes $N$ is necessarily lower than or equal to the number of high-fidelity snapshots per cluster $n_s$. For simplification purposes, we take $N=n_s$. Given that $\mathcal{N}=2000$, the number of POD modes $N$ must be lower than $\lfloor \mathcal{N}^{1/3} \rfloor = 12$. To effectively reduce the computational cost of high-fidelity simulations, the maximum number of modes considered in this paper is $N=5$.

Let us say that we are given a budget of $20$ high-fidelity simulations. The hyperparameters must satisfy the inequality $K n_s = K N \leq 20$. Our thresholds for the mean projection error and the mean gain are  $\eta^* = 0.35$ and $G_{r}^* = 2$. A polynomial of degree $2$ is considered for the model $p(K)$ for the classification accuracy, and its coefficients are determined by imposing $p(1)=1$, $p(6)=0.8$ (value taken from~\cite{ROM-net}), $p(K_{\textrm{max}}) = 1/K_{\textrm{max}}$ (accuracy of a random guess for balanced classes) and $K_{\textrm{max}}=20$.

\begin{figure}[!h]
\centering
\includegraphics[scale=0.39]{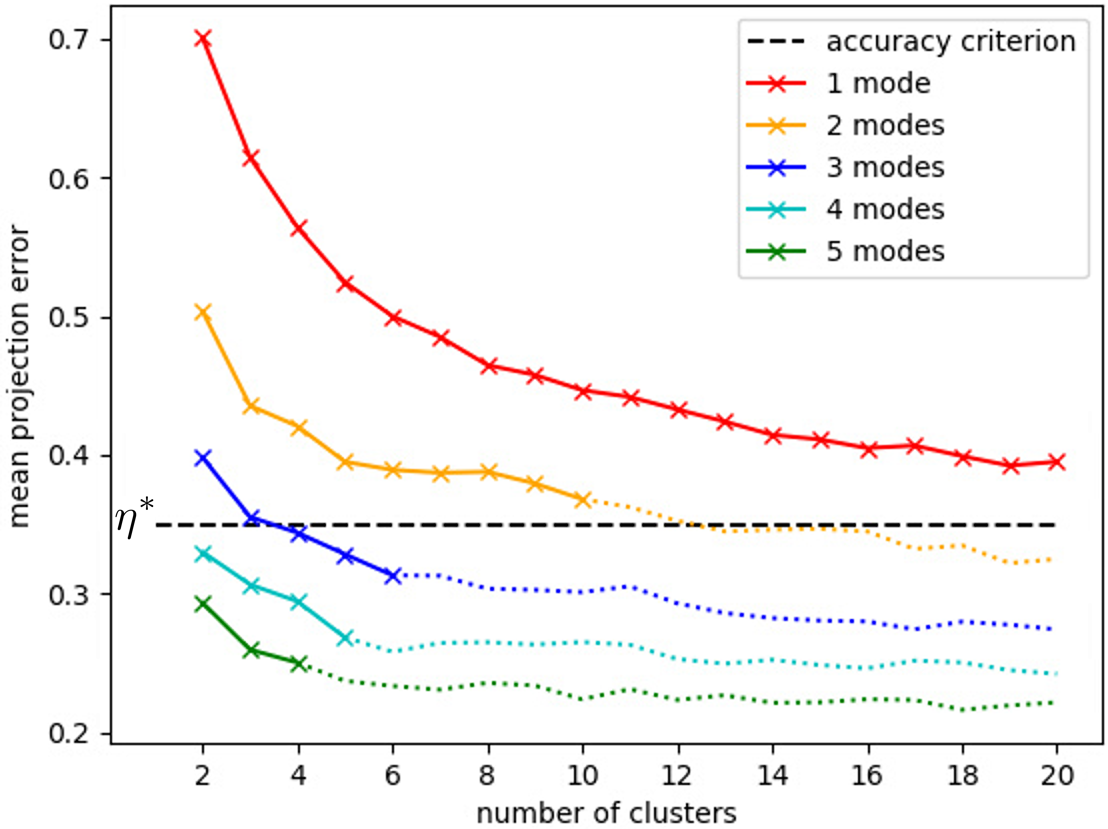}
\caption{Mean projection error as a function of the number of clusters $K$ for different number of modes $N$. Dotted lines indicate the configurations which do not comply with the allocated number of high-fidelity snapshots.}
\label{errorCurves}
\end{figure}

\begin{figure}[!h]
\centering
\includegraphics[scale=0.78]{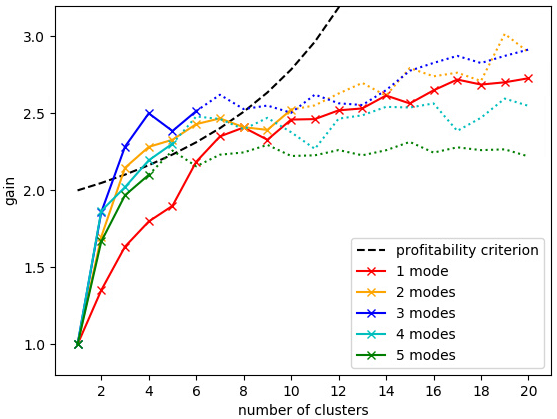}
\caption{Gain as a function of the number of clusters $K$ for different number of modes $N$. Dotted lines indicate the configurations which do not comply with the allocated number of high-fidelity snapshots.}
\label{gainCurves}
\end{figure}

Figure~\ref{errorCurves} gives the mean projection error as a function of $K$ and $N$. For $N=4$ and $N=5$ modes, the mean projection errors are below the tolerance $\eta^*$ for all $K$. For $N=3$, the accuracy criterion is satisfied for $K\geq 4$. The mean projection error for $N=2$ modes falls below the tolerance for $K \geq 13$, which does not conform to the constraint $K \leq 10$ imposed by the allocated number of high-fidelity snapshots. With $N=1$ mode, the mean projection error remains too large, which rejects configurations with $N=1$. The configurations $(K,N)$ satisfying the accuracy criterion and respecting the budget for high-fidelity snapshots are $K=4,5,6$ for $N=3$, $K=2,3,4,5$ for $N=4$, and $K=2,3,4$ for $N=5$. 

\begin{figure}[!h]
\centering
\includegraphics[scale=0.78]{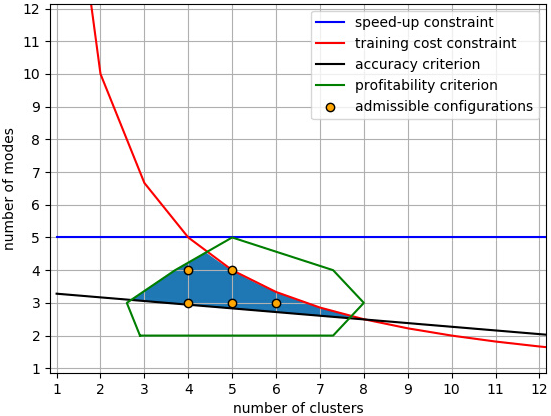}
\caption{Hyperparameters admissible set.}
\label{admSet}
\end{figure}

The gain curves are given in Figure~\ref{gainCurves}. The dashed line in black delimits the ROM-net's profitability domain: configurations under this curve are irrelevant, either because the corresponding expected gain is too low, or because misclassification errors would be too frequent because too many classes are considered. The configuration $(K,N)=(5,5)$ meets both gain and accuracy requirements, but violates the constraint $K \leq 4$ for $N=5$ and thus requires too many high-fidelity snapshots. For $(K,N)=(3,3)$, the gain is large enough but the mean projection error is larger than the tolerance, as seen in Figure~\ref{errorCurves}. Finally, the admissible configurations are $K=4,5,6$ for $N=3$, and $K=4,5$ for $N=4$. The hyperparameters admissible set is represented in Figure~\ref{admSet}. Among the admissible configurations, those with $N=3$ are more interesting in terms of speed of online reduced simulations. The lowest mean projection error is obtained for $K=6$ when $N=3$, see Figure~\ref{errorCurves}. Therefore, we choose the hyperparameters $(K,N,n_s) = (6,3,3)$, corresponding to the lower right dot in Figure~\ref{admSet}.

\begin{remark}
It has been decided to take the configuration with the best accuracy among the admissible configurations with the smallest value for $N$, in order to have a simple and systematic approach for hyperparameters calibration. However, in the present example, one could also use the \textit{elbow method}. The elbow method is commonly used for selecting the number of clusters for k-means clustering or the number of components for a PCA. It consists in choosing the elbow or knee point of the curve of an evaluation criterion. In spite of the difficulties of defining clearly the elbow point in some situations, this method raises interesting questions. In our example, if one uses the elbow method with the error curve, the best number of clusters is still $K=6$: for $N=3$, $K=6$ is the elbow point. When using this method with the gain curve, the best number of clusters turns out to be $K=4$, even when considering a smoothed version of the blue curve in Figure~\ref{gainCurves} to avoid undesirable fluctuations due to sampling and medoids initializations. Indeed, taking $K=5$ or $K=6$ does not significantly improve the gain when $N=3$, whereas the number of high-fidelity snapshots and the complexity of the classification problem would be increased. The practical method presented in this paper can be adapted according to the user's priorities between training cost, online speed, accuracy, and gain.
\end{remark}

\subsubsection{Comparison of different model order reduction strategies}

Let $x = (s,\zeta, \epsilon(\zeta))$ denote the parameter of the problem. After projection of the source term in the finite-element basis, the parameter $x$ is represented by a $\mathcal{N}+2$-dimensional vector $\g{x}$ whose coordinates are centered and scaled to unit variance. This way, distances in the parameter space can be computed with the Euclidean distance:\begin{equation}
\delta_{\mathcal{X}}(x,x') = || \g{x}-\g{x'} ||_2 .
\end{equation} Introducing the notation $\g{q}(x)\in\mathbb{R}^{\mathcal{N}-2}$ for the solution of Equation~\eqref{Kqf} for a given parameter $x$, one can define a physics-informed dissimilarity measure $\delta_{\mathcal{H}}$ using the Euclidean distance in the solution space:\begin{equation}
\delta_{\mathcal{H}}(x,x') = || \g{q}(x)-\g{q}(x') ||_2 .
\end{equation} These dissimilarity measures are compared with the ROM-oriented dissimilarity measure $\delta$ introduced in Equation~\eqref{DefROMdissimilarity}, obtained by computing the sine dissimilarity in the solution space. K-medoids clustering is used for both snapshots selection and manifold partitioning in conjunction with one of these three dissimilarity measures. Different model order reduction strategies are compared in terms of projection errors under the following setting:
\begin{itemize}
\item \textbf{Equivalent number of snapshots:} all the strategies use the same total number of snapshots, which ensures equal budgets for high-fidelity simulations in the training phase. It is recalled that the high-fidelity snapshots are given by high-fidelity simulations that are more expensive than the simplified simulations used to generate the database, find clusters and evaluate their quality.
\item \textbf{Equivalent number of POD modes:} all the ROBs use the same number of modes, which ensures equivalent speed-ups when exploiting the ROMs.
\end{itemize}
If a dictionary of $K$ local ROBs is compared with a global ROB made of $N$ modes, then each local ROB must have $N$ modes. For the construction of these local cluster-specific ROBs, $n_s = N$ snapshots are selected in each cluster using the two-stage hierarchical k-medoids clustering procedure. Hence, the total number of snapshots is $K n_s = K N$. Snapshots for the construction of the global ROB are therefore selected by taking the medoids of a single k-medoids clustering with $KN$ clusters.

\noindent Six model order reduction strategies are considered, namely:
\begin{itemize}
\item Three global ROBs containing $N$ modes computed from $KN$ snapshots. The snapshots are selected thanks to a k-medoids cluster analysis with $KN$ clusters, using different dissimilarities:
\begin{itemize}
\item \textbf{Global ROM 1} uses the dissimilarity $\delta_{\mathcal{X}}$ (Euclidean distance in the parameter space).
\item \textbf{Global ROM 2} uses the dissimilarity $\delta_{\mathcal{H}}$ (Euclidean distance in the solution space).
\item \textbf{Global ROM 3} uses the ROM-oriented dissimilarity $\delta$ (sine dissimilarity in the solution space).
\end{itemize}
\item Three ROM dictionaries consisting of $K$ local ROBs with $N$ modes each. Each local ROB is inferred from $N$ snapshots. Again, k-medoids is applied with different dissimilarity measures:
\begin{itemize}
\item \textbf{ROM dictionary 1} uses the dissimilarity $\delta_{\mathcal{X}}$ (Euclidean distance in the parameter space). This strategy is the most natural and simple one among ROM dictionaries.
\item \textbf{ROM dictionary 2} uses the dissimilarity $\delta_{\mathcal{H}}$ (Euclidean distance in the solution space) like in~\cite{localROB, localROB2, Amsallem2015localHROM, RyckelynckComputerVision, Grimberg2020}. This strategy belongs to physics-informed strategies.
\item \textbf{ROM dictionary 3} uses the ROM-oriented dissimilarity $\delta$ (sine dissimilarity in the solution space). This is the strategy we have introduced in this paper for dictionary-based ROM-nets. Like ROM dictionary 2, it relies on a physics-informed cluster analysis, but with another dissimilarity.
\end{itemize}
\end{itemize}

\begin{figure}[!h]
\centering
\includegraphics[scale=0.8]{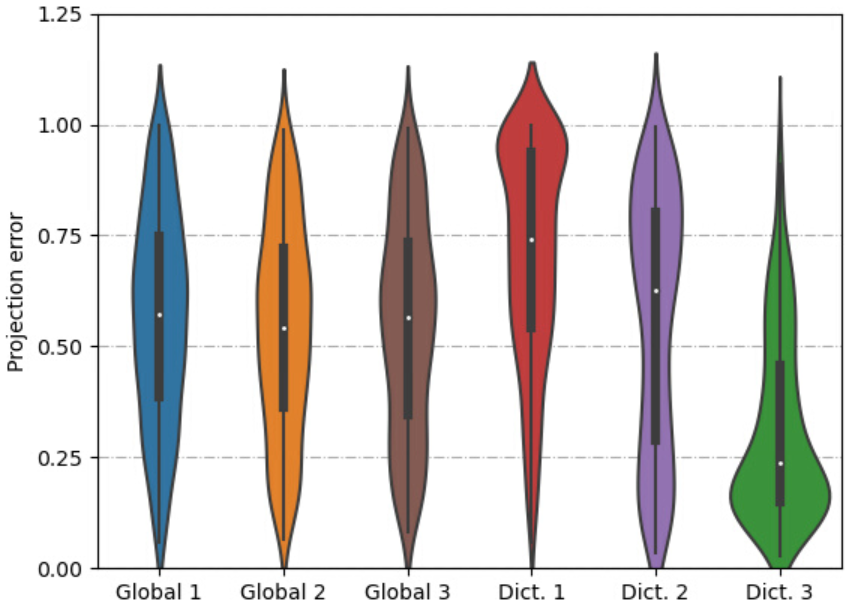}
\caption{Violin plots of the projection errors for different model order reduction strategies, with $(K,N,n_s)=(6,3,3)$.}
\label{errDistrib}
\end{figure}

In this section, the comparison is presented for $K=6$ and $N= n_s =3$, the configuration identified in the previous section thanks to the gain curves and the projection error curves. Projection errors as defined in Equation~\eqref{ProjectionError} are computed for the $500$ test examples for each strategy, which enables estimating their probability density functions using Gaussian kernel density estimation (see section~6.6.1. of~\cite{Hastie2005TheEO}). The violin plots of the projection errors are given in Figure~\ref{errDistrib}, and the values of the quartiles and expectations are given in Table~\ref{StatsProjErrors}. The third ROM dictionary using the ROM-oriented dissimilarity clearly outperforms the other strategies. Although using a physics-informed clustering procedure, ROM dictionary 2 fails to improve the performances of global ROMs on this specific example. This result illustrates the fact that the Euclidean distance is not always appropriate for model order reduction purposes. ROM dictionary 1 gives the worst results, showing that integrating physics in cluster analyses is crucial when the final objective is to build local approximation spaces. Interestingly, these results also show that using local ROBs can deteriorate the performances of a global ROB when choosing an improper dissimilarity measure for clustering. In this example, the three global ROMs give approximately the same projection errors. These errors are lower than those obtained with ROM dictionary 1 and ROM dictionary 2 because the global ROMs have more relevant snapshots, since they use $KN$ well-distributed snapshots instead of $N$ badly-distributed snapshots. Hence, the dissimilarities $\delta_{\mathcal{X}}$ and $\delta_{\mathcal{H}}$ both define inefficient notions of locality in this example.

\begin{table}[h!]
\centering
\caption{Quartiles and means of the projection errors for different model order reduction strategies, with $(K,N,n_s)=(6,3,3)$.}
      \begin{tabular}{|c|c|c|c|c|c|}
        \hline & & & & & \\[-1em]
        Strategy & Dissimilarity & $Q_1$ & Median & $Q_3$ & Mean  \\
         & & & & & \\[-1em] \hline
        Global ROM 1 & $\delta_{\mathcal{X}}$ & $0.3863$ & $0.5735$ & $0.7477$ & $0.5636$ \\
        Global ROM 2 & $\delta_{\mathcal{H}}$ & $0.3611$ & $0.5427$ & $0.7208$ & $0.5397$ \\
        Global ROM 3 & $\delta$ & $0.3460$ & $0.5666$ & $0.7346$ & $0.5480$ \\
        ROM dictionary 1 & $\delta_{\mathcal{X}}$ & $0.5434$ & $0.7412$ & $0.9379$ & $0.7071$ \\
        ROM dictionary 2 & $\delta_{\mathcal{H}}$ & $0.2874$ & $0.6280$ & $0.8038$ & $0.5586$ \\
        ROM dictionary 3 & $\delta$ & $0.1482$ & $0.2369$ & $0.4584$ & $0.3132$ \\
\hline
      \end{tabular}
      \label{StatsProjErrors}
\end{table}

\begin{remark}
Figure~\ref{errDistrib} gives projection errors obtained when choosing the correct cluster and thus the most suitable local ROB. The ROM-net's classification errors would have the effect of moving the distribution of ROM dictionary 3 towards larger errors, reducing the gap between the errors made by the different model order reduction strategies. Therefore, particular attention must be paid to the training of the ROM-net's classifier.
\end{remark}

\begin{figure}[!h]
\centering
\includegraphics[scale=0.32]{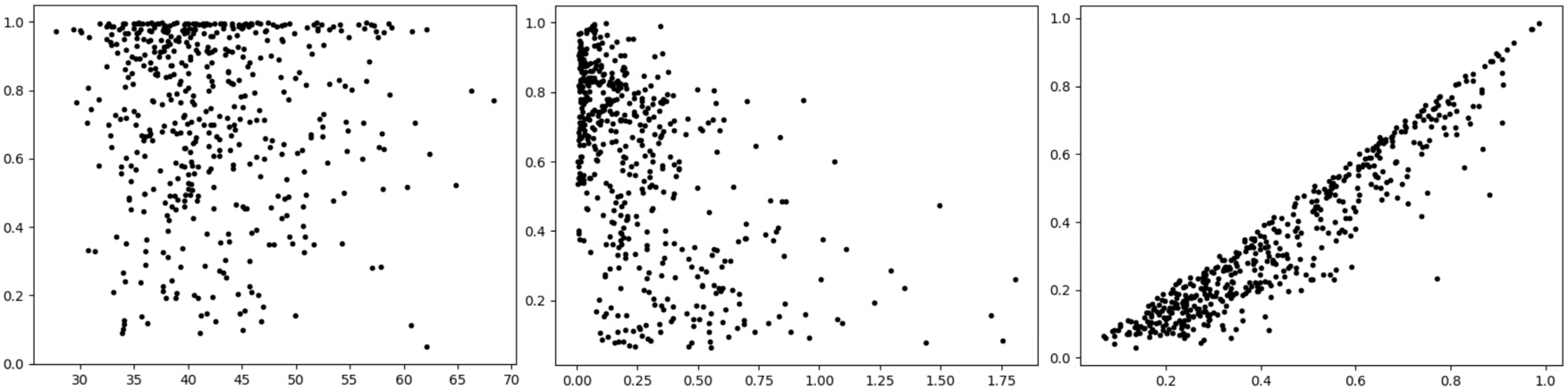}
\caption{Scatter plots giving the projection error (y-axis) for test data against the dissimilarity with the closest snapshot (x-axis). From the left to the right: Euclidean distances in the parameter space, Euclidean distances in the solution space, and ROM-oriented dissimilarity.}
\label{correlations}
\end{figure}

Figure~\ref{correlations} plots the projection error against the dissimilarity measure $\delta_{\mathcal{X}}$ (left), $\delta_{\mathcal{H}}$ (middle) and $\delta$ (right) separating a test example from its closest snapshot. One can clearly see the correlation between the projection error and our ROM-oriented dissimilarity $\delta$, contrasting with the absence of correlations between the projection error and the other dissimilarities. {This figure also shows that the dissimilarity with the closest snapshot is generally larger than the projection error onto the POD basis. Indeed, for this time-independent problem, the sine dissimilarity corresponds to the relative projection error; the closest snapshot for this dissimilarity is therefore a better approximation of the solution than the orthogonal projection onto the POD basis, because the POD basis does not perfectly approximate each of its snapshots.}

{In the following Sections~\ref{sec:2DTransportSnapshots}-\ref{sec:3DMecaApplication}, we compare the performances of dictionaries of local ROMs obtained by clustering using the $L^2$ and ROM-oriented dissimlarity measures for various physics problems.}

\subsection{2D advection equation}
\label{sec:2DTransportSnapshots}

\begin{figure}[!h]
\centering
\includegraphics[width=\textwidth]{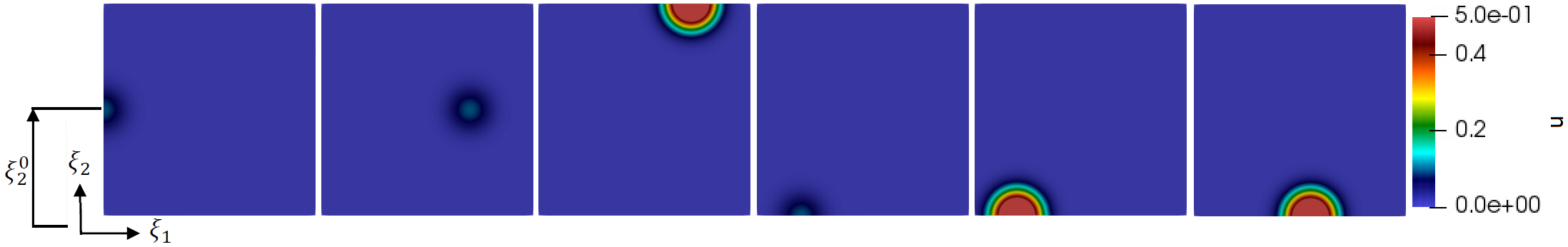}
\caption{Snapshots of the solution $u$ of Equation~\eqref{eq:transportEq} for various values of $U^0$, $\xi_2^0$ and $t$.}
\label{fig:2DTransportSnapshots}
\end{figure}

We consider the following advection equation:
\begin{equation}
\left\{
\begin{aligned}
\frac{\partial u}{\partial t}+c\frac{\partial u}{\partial \xi_1}&=0\\
u(\xi_1,\xi_2,t=0)&=U^0\exp\left(-\frac{\xi_1^2+\left(\xi_2-\xi_2^0\right)}{l^2}\right),
\end{aligned}
\right.
\label{eq:transportEq}
\end{equation}
where $(\xi_1,\xi_2)\in [0,1]^2$ and $t\in [0,100]$. The analytical solution is known: $\displaystyle u(\xi_1,\xi_2,t)=U^0\exp\left(-\frac{\left(\xi_1-ct\right)^2+\left(\xi_2-\xi_2^0\right)}{l^2}\right)$. The quantities $l=0.1$ and $c=1$ are constants, while $U^0\in \{0.1,1\}$ and $\xi_2^0\in \{0,0.5,1\}$ are the parameters of the problem. Some snapshots are illustrated on a mesh with 10201 vertices in Figure~\ref{fig:2DTransportSnapshots}, for various values of $U^0$, $\xi_2^0$ and $t$. A total of $600$ snapshots are generated.

\begin{figure}[!h]
\centering
\includegraphics[width=\textwidth]{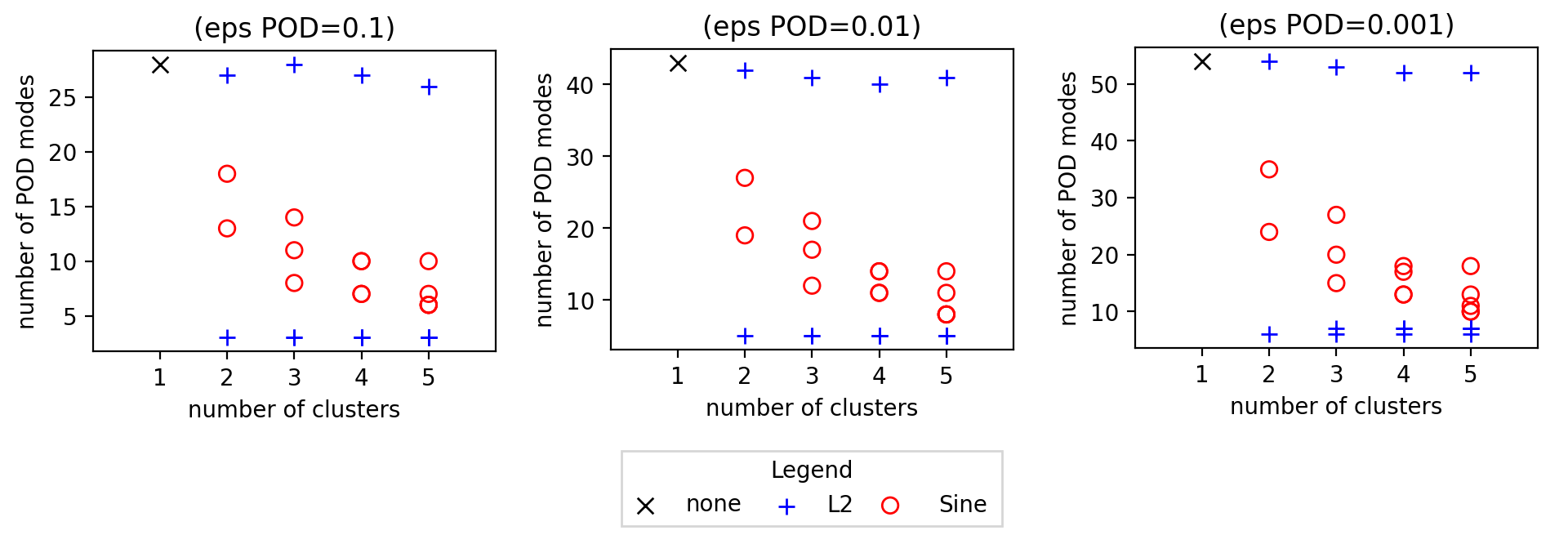}
\caption{Number of POD modes for each local basis with respect to the number of clusters, for various accuracy criterions of the POD truncature, applied to the advection problem.}
\label{fig:2DTransportNbePODModes}
\end{figure}

Figure~\ref{fig:2DTransportNbePODModes} shows the number of POD modes for each local basis with respect to the number of clusters, for various accuracy criterions of the POD truncature. For all the considered levels of truncature, the clustering carried out using the sine dissimilarity measure leads to the smallest maximal size of local reduced-order basis.

\begin{figure}[!h]
\centering
\includegraphics[width=\textwidth]{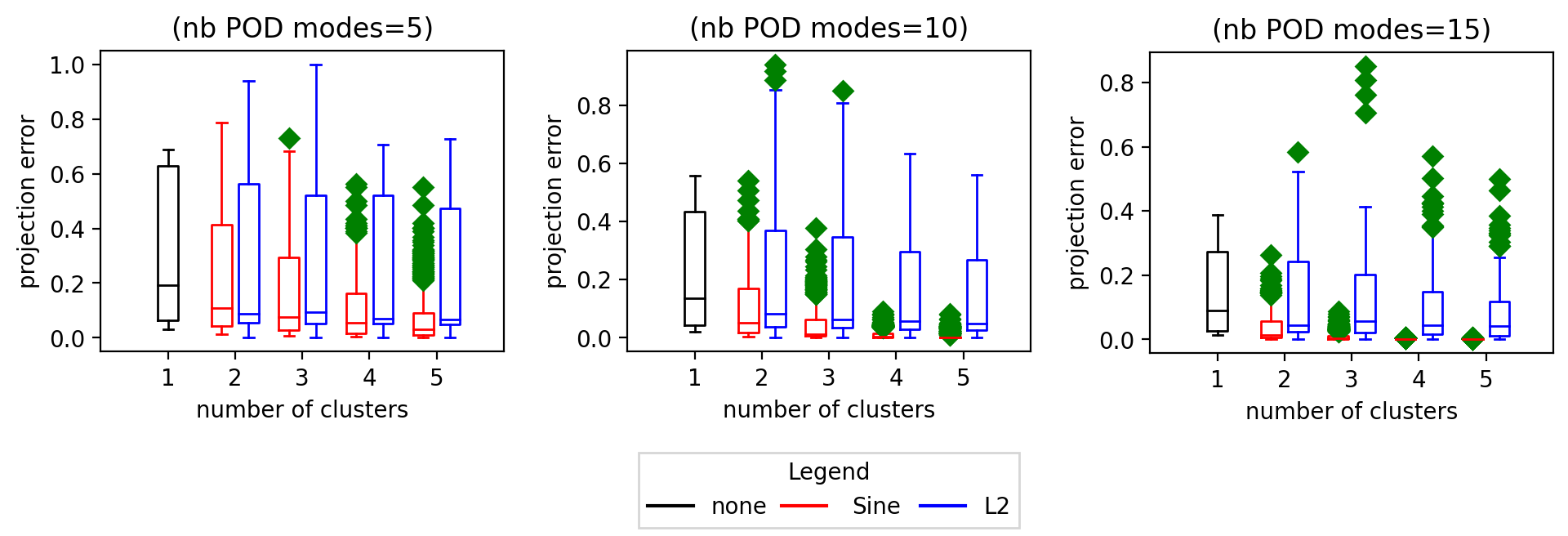}
\caption{Box plots for the projection errors with respect to the number of clusters, for different size of the local POD basis, applied to the advection problem. The box plots show the median, first and third quartiles, as well as extremal values (in the form of green dots if outliers).}
\label{fig:2DTransportProjError}
\end{figure}

Figure~\ref{fig:2DTransportProjError} shows the projection errors with respect to the number of clusters, for different size of the local POD basis. For all the considered size of local reduced-order basis, the clustering carried out using the sine dissimilarity measure leads to the smallest projection errors.

\begin{figure}[!h]
\centering
\includegraphics[width=\textwidth]{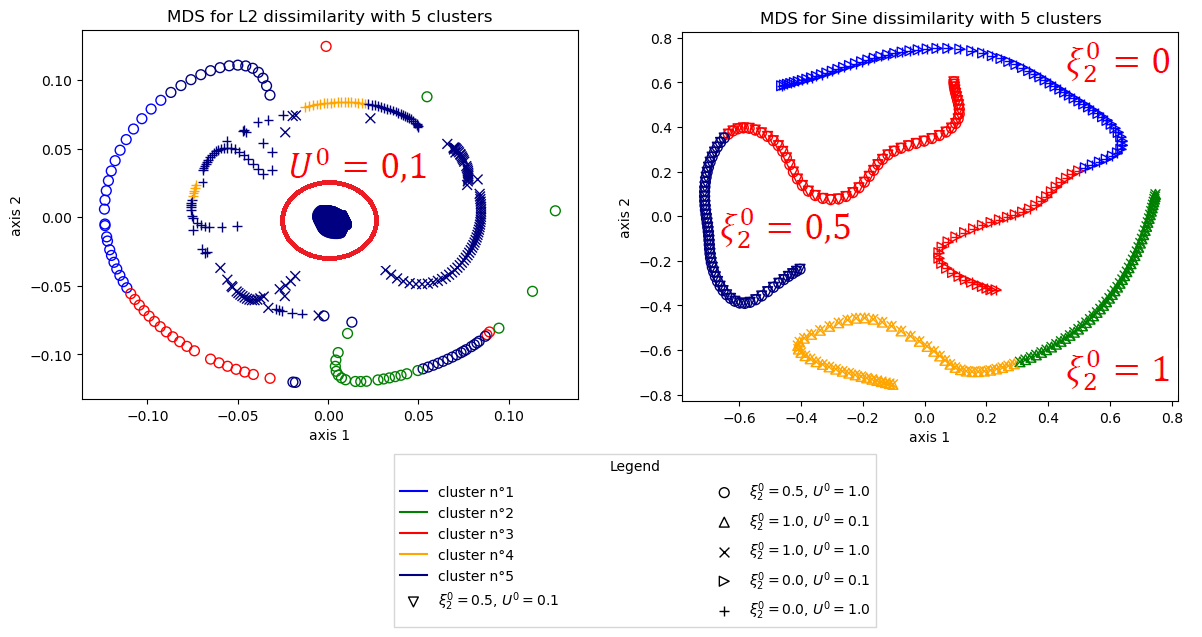}
\caption{MDS representations of $L^2$ and sine dissimilary measures, with coloring depending on cluster affectation for 5 clusters, applied to the advection problem.}
\label{fig:2DTransportMDS}
\end{figure}

Figure~\ref{fig:2DTransportMDS} shows MultiDimensional Scaling (MDS) representations of $L^2$ and sine dissimilary measures, with coloring depending on cluster affectation for 5 clusters. A 2-dimensional MDS representation aims to locate points, each representing a solution of the considered physical problem, in such a fashion that the pairwise 2D Euclidean distances between each points is as close as possible to the corresponding dissimilarity. In the $L^2$ dissimilarity case, all the snapshots corresponding to $U^0=0.1$ in the MDS representation are very close to each other, which means that their corresponding $L^2$ pairwise distances are small compared to the rest of the pairwise distances. It is explained by the fact that the $L^2$ norm quantifies magnitudes. Hence, when applying a k-medoid clustering algorithm, all the snapshots corresponding to $U^0=0.1$ are affected to the same cluster. However, these small-magnitude snapshots contain all the independant directions of the solution function space described by the whole snapshot set. As a consequence, the local reduced-order model corresponding to the cluster containing these small-magnitude snapshots has the same reducibility as the complete set, for any accuracy level and even when increasing the number of clusters. This is illutrated in Figure~\ref{fig:2DTransportNbePODModes}, where the $L^2$ dissimilary case exhibits one local reduced-order model having a number of mode very close to the global reduced-order model. On the contrary, the MDS for the sine dissimilarity in Figure~\ref{fig:2DTransportMDS} shows three trajectories corresponding to the three different values of $\xi_2^0$. Actually, each pair of snapshots corresponding of same values of $\xi_2^0$ and $t$, for $U^0=0.1$ and $1$, are at the same location of the MDS representation, which means that their pairwise sine dissimilarities are zero. Hence, a new snapshot collinear to an existing snapshot do not increase the number of independant directions in the snapshot set. The corresponding clustering produces balanced clusters, having local small-sized reduced-order basis as seen in Figure~\ref{fig:2DTransportNbePODModes}.

Since the same analysis can be done in the following three additional numerical experiments, we do not repeat it and simply explain the new physical settings.

\subsection{2D incompressible Navier-Stokes}
\label{sec:2DNSApplication}

\begin{figure}[!h]
\centering
\includegraphics[width=0.7\textwidth]{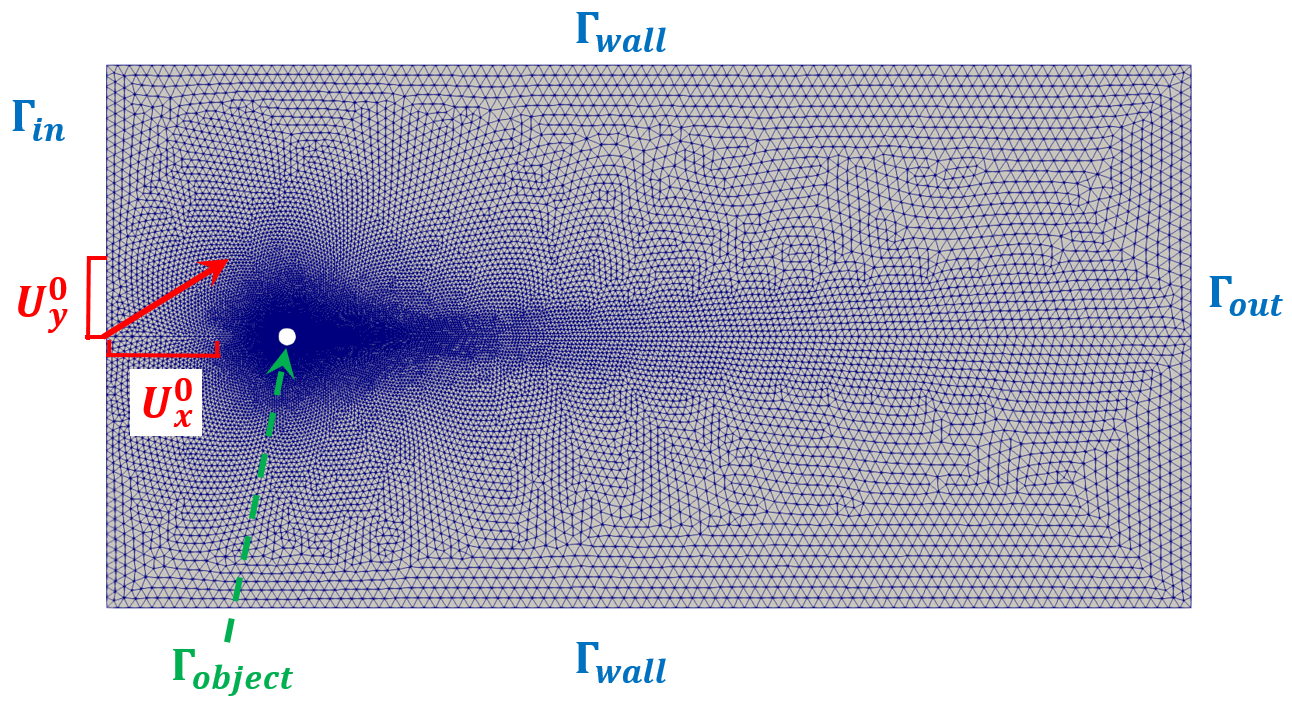}
\caption{Setting and mesh for the Navier-Stokes problem.}
\label{fig:2DCFDMesh}
\end{figure}

\begin{figure}[!h]
\centering
\includegraphics[width=\textwidth]{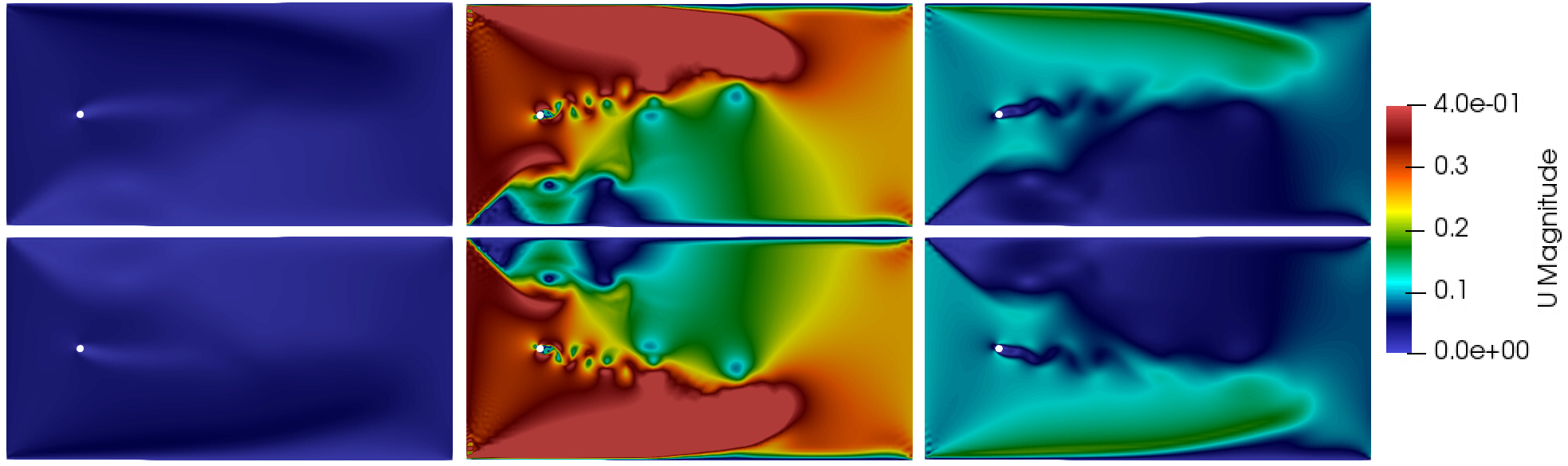}
\caption{Snapshots of the solution velocity field $U$ of the Navier-Stokes problem: last time step for $\left(U^0_x, U^0_y\right)\in \left\{(0.025, 0.025)\right.$, $(0.075, 0.075)$, $(0.25, 0.25)$, $(0.025, -0.025)$, $(0.075, -0.075)$, $\left.(0.25, -0.25)\right\}$.}
\label{fig:2DCFDSnapshots}
\end{figure}

We consider the 2D incompressible Navier-Stokes equations in the setting illustrated in Figure~\ref{fig:2DCFDMesh}: the air flows from the left to the right in the rectangular domain, with a uniform Dirichlet boundary condition for the velocity $U$ on $\Gamma_{in}$, outflow boundary condition on $\Gamma_{in}$, and no-slip boundary condition on the walls $\Gamma_{wall}$ and on the circular object $\Gamma_{object}$. The mesh is constituted of 19818 vertices; the low-Mach number solver YALES2~\cite{yales2} for unstructured grids is used. The parameters of the problem are components $U^0_x$ and $U^0_y$ of the uniform incoming velocity boundary condition, and we consider 6 temporal simulations for $(U^0_x, U^0_y)\in \left\{(0.025, 0.025)\right.$, $(0.075, 0.075)$, $(0.25, 0.25)$, $(0.025, -0.025)$, $(0.075, -0.075)$, $\left.(0.25, -0.25)\right\}$ leading to a total of 600 snapshots. The last time step of each of these simulations is illustrated in Figure~\ref{fig:2DCFDSnapshots}.

\begin{figure}[!h]
\centering
\includegraphics[width=\textwidth]{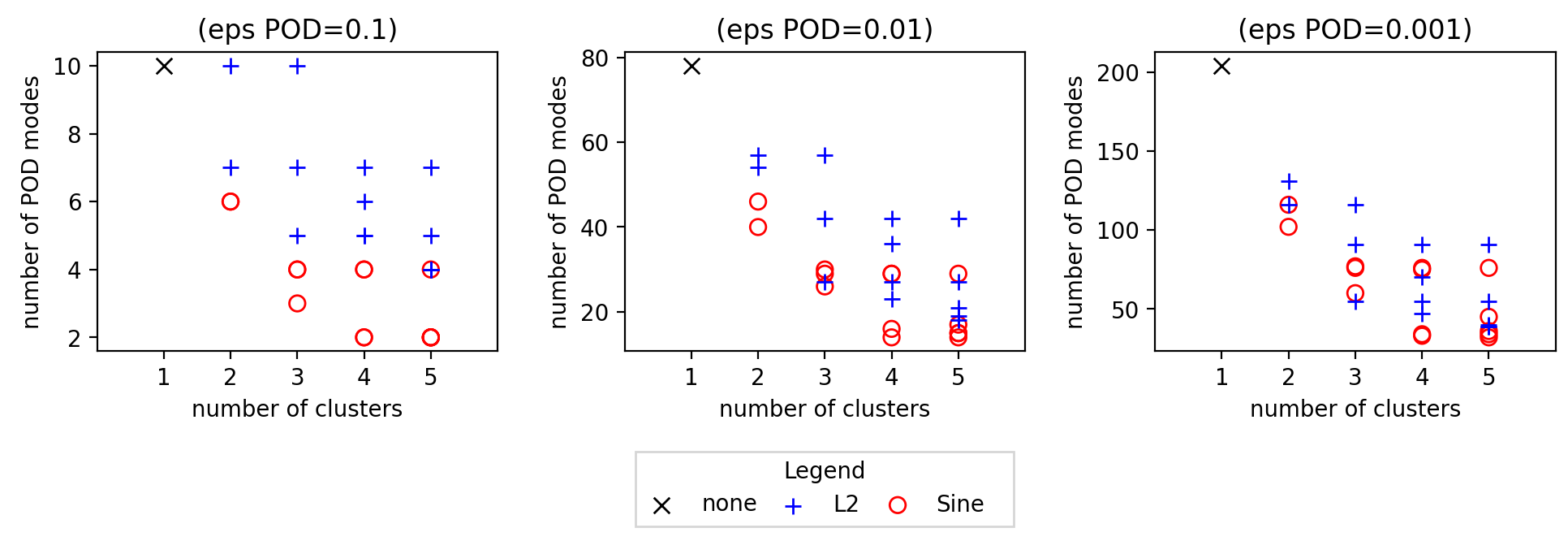}
\caption{Number of POD modes for each local basis with respect to the number of clusters, for various accuracy criterions of the POD truncature, applied to the Navier-Stokes problem.}
\label{fig:2DCFDNbePODModes}
\end{figure}

\begin{figure}[!h]
\centering
\includegraphics[width=\textwidth]{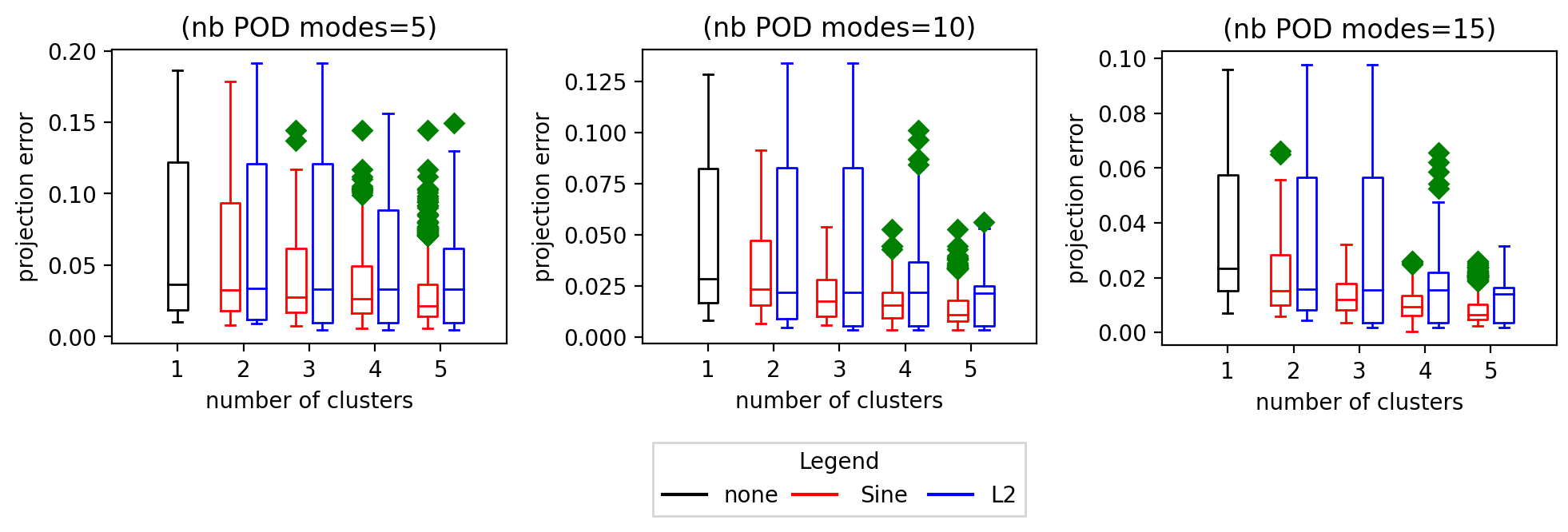}
\caption{Box plots for the projection errors with respect to the number of clusters, for different size of the local POD basis, applied to the Navier-Stokes problem. The box plots show the median, first and third quartiles, as well as extremal values (in the form of green dots if outliers).}
\label{fig:2DCFDProjError}
\end{figure}

The improved performance of the sine dissimilarity based clustering, with respect to the $L^2$ one, is illustrated in Figures~\ref{fig:2DCFDNbePODModes}-\ref{fig:2DCFDProjError}.

\begin{figure}[!h]
\centering
\includegraphics[width=\textwidth]{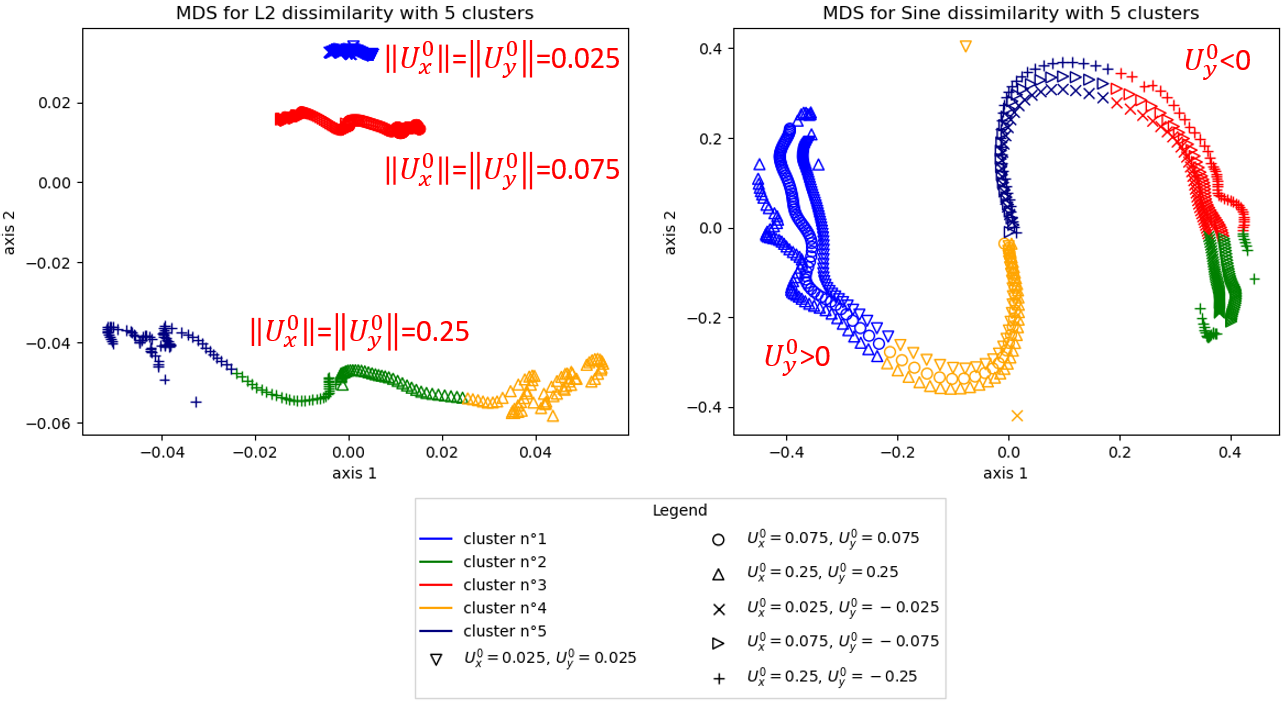}
\caption{MDS representations of $L^2$ and sine dissimilary measures, with coloring depending on cluster affectation for 5 clusters, applied to the Navier-Stokes problem.}
\label{fig:2DCFDMDS}
\end{figure}

The MDS representations in Figure~\ref{fig:2DCFDMDS} shows the advantages of the sine dissimilarity in the same fashion as the previous section: with $L^2$-based clustering, snapshots are grouped by magnitude, where all the small-magnitude ones are tightly packed in a single cluster, whereas with sine-based clustering, snapshots are grouped by direction of the initial condition.

\subsection{2D heat equation}
\label{sec:2DHeatApplication}

\begin{figure}[!h]
\centering
\includegraphics[width=\textwidth]{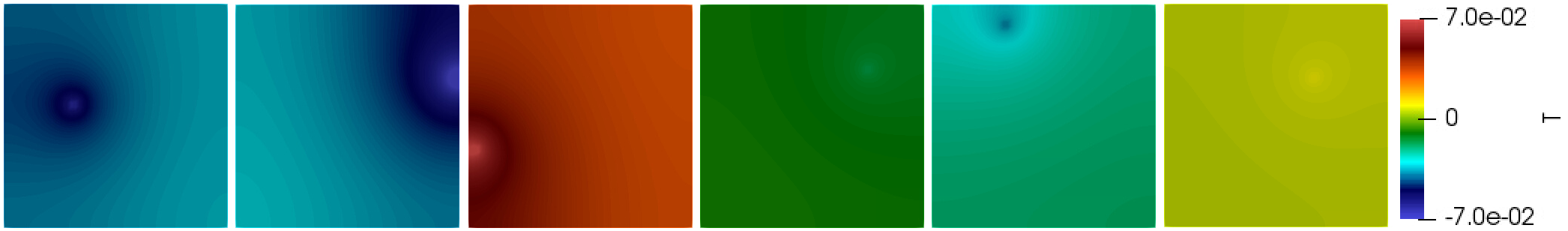}
\caption{Snapshots of the solution temperature field $T$ of the 2D heat problem for various values of the parameters.}
\label{fig:2DThermalSnapshots}
\end{figure}

We consider the 2D transient linear heat equation on a square domain, with a localized volumetric heat source term. The parameters of the problem are the location of the source term and its magnitude; only the final time step is kept, see Figure~\ref{fig:2DThermalSnapshots} for examples of snapshots. The mesh is composed of 2601 vertices and 500 snapshots for 500 random values of the parameters are computed using the finite-element software \textit{Z-set}~\cite{zset}.

\begin{figure}[!h]
\centering
\includegraphics[width=\textwidth]{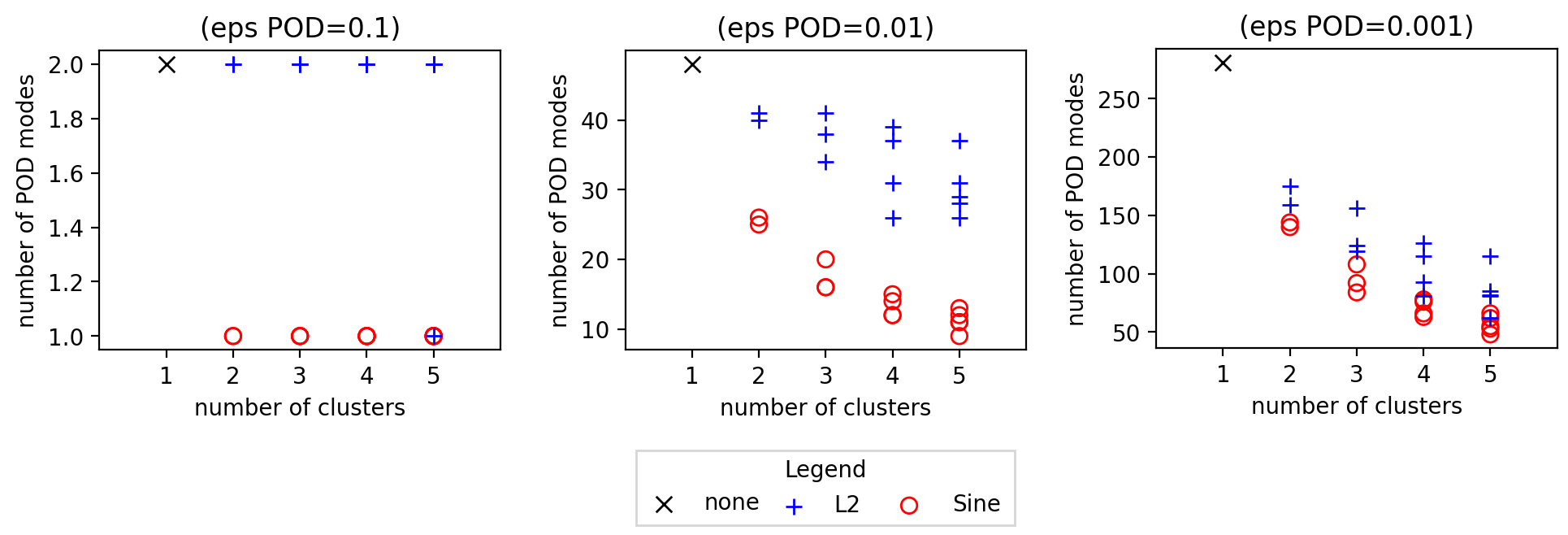}
\caption{Number of POD modes for each local basis with respect to the number of clusters, for various accuracy criterions of the POD truncature, applied to the 2D heat problem.}
\label{fig:2DThermalNbePODModes}
\end{figure}

\begin{figure}[!h]
\centering
\includegraphics[width=\textwidth]{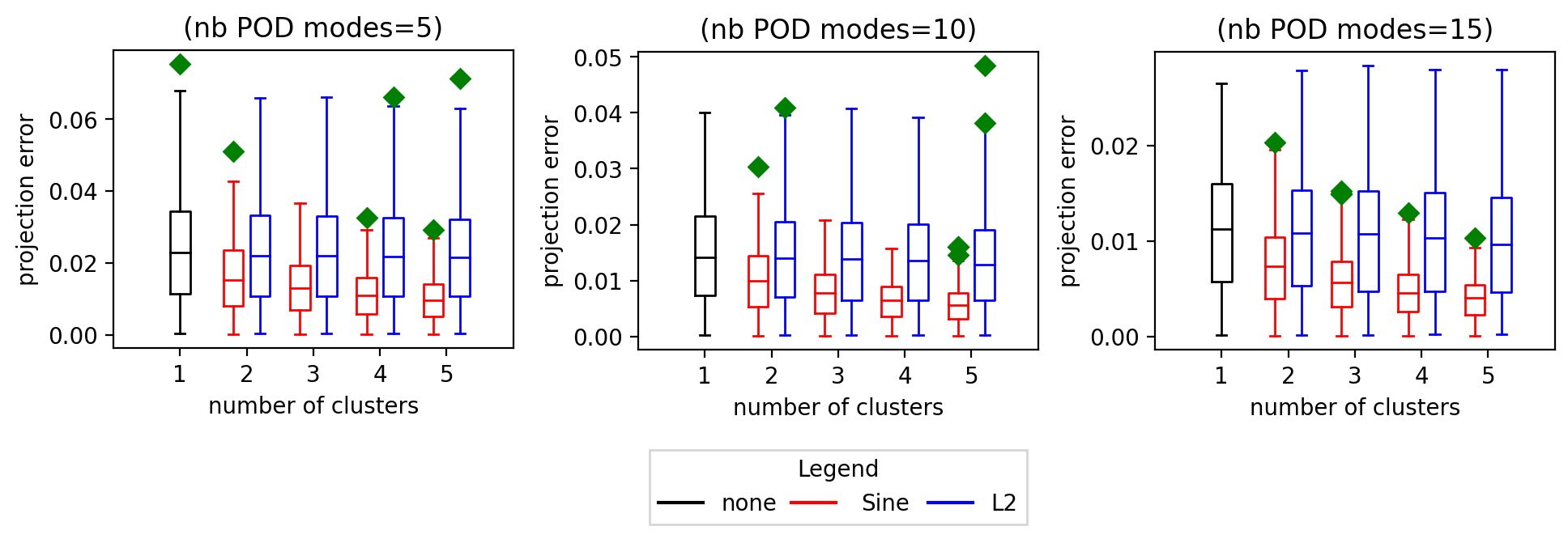}
\caption{Box plots for the projection errors with respect to the number of clusters, for different size of the local POD basis, applied to the 2D heat problem. The box plots show the median, first and third quartiles, as well as extremal values (in the form of green dots if outliers).}
\label{fig:2DThermalProjError}
\end{figure}

The improved performance of the sine dissimilarity based clustering, with respect to the $L^2$ one, is illustrated in Figures~\ref{fig:2DThermalNbePODModes}-\ref{fig:2DThermalProjError}.

\subsection{3D nonlinear structural mechanics}
\label{sec:3DMecaApplication}

\begin{figure}[!h]
\centering
\includegraphics[width=\textwidth]{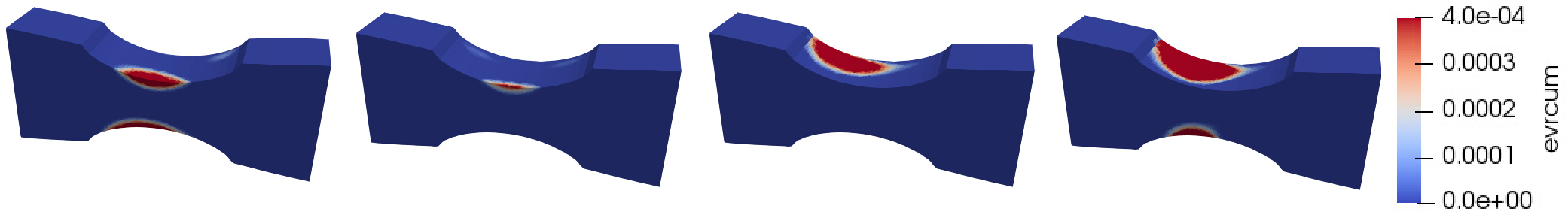}
\caption{Snapshots of the accumulated plasticity field $evrcum$ of the mechanical problem for various values of the parameters.}
\label{fig:3DMecaSnapshots}
\end{figure}

We consider a 3D quasistatic nonlinear structural mechanics problem: an object is rotated along an axis intersecting the center of gravity of this object. The material is modeled by a viscoplastic constitutive law with a Von Mises criterion and a Norton flow. The orientation of the axis is the parameter of the problem; only the final time step in kept and the quantity of interest is the accumulated plasticity field, see Figure~\ref{fig:3DMecaSnapshots} for examples of snapshots. The mesh is composed of 61741 vertices and 100 snapshots for 100 random values of the parameters are computed using the finite-element software \textit{Z-set}~\cite{zset}.

\begin{figure}[!h]
\centering
\includegraphics[width=\textwidth]{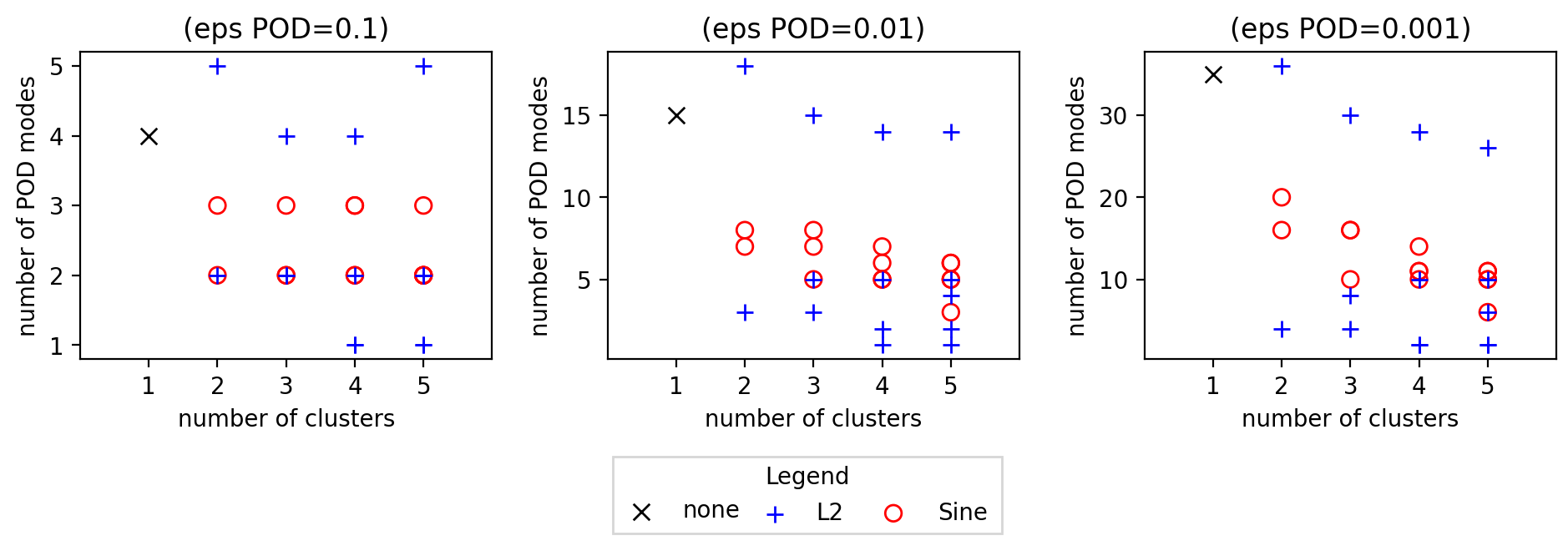}
\caption{Number of POD modes for each local basis with respect to the number of clusters, for various accuracy criterions of the POD truncature, applied to the mechanical problem.}
\label{fig:3DMecaNbePODModes}
\end{figure}

\begin{figure}[!h]
\centering
\includegraphics[width=\textwidth]{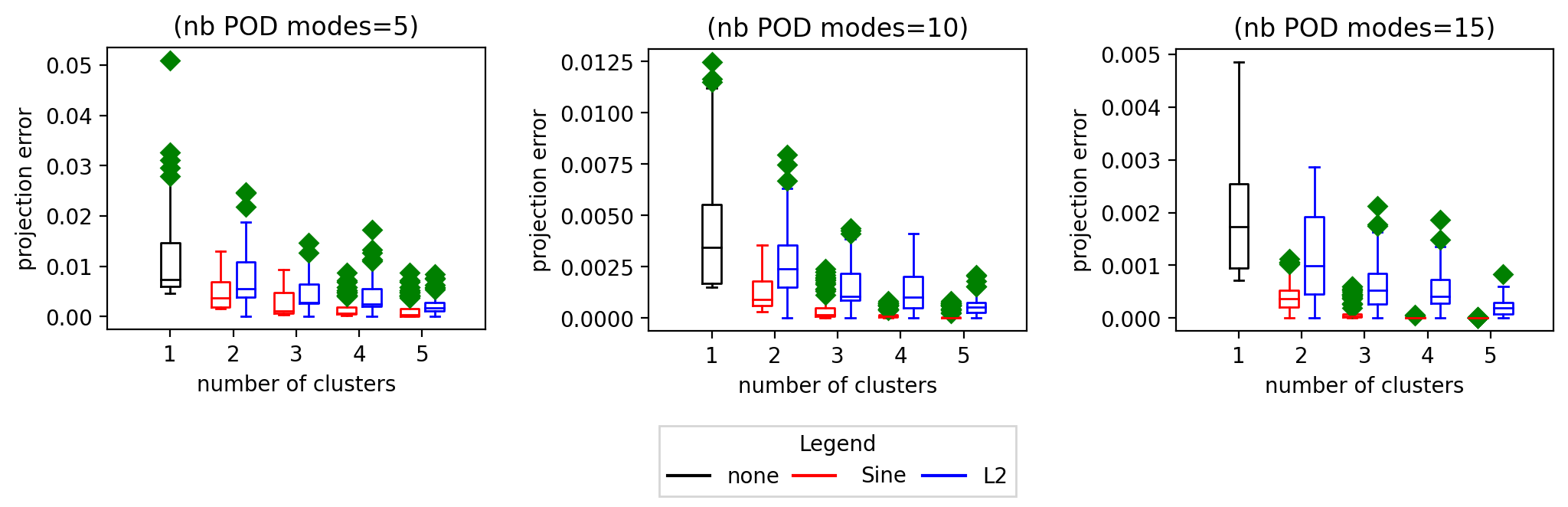}
\caption{Box plots for the projection errors with respect to the number of clusters, for different size of the local POD basis, applied to themechanical problem. The box plots show the median, first and third quartiles, as well as extremal values (in the form of green dots if outliers).}
\label{fig:3DMecaProjError}
\end{figure}

The improved performance of the sine dissimilarity based clustering, with respect to the $L^2$ one, is illustrated in Figure~\ref{fig:3DMecaProjError}.

\section{Conclusion}

Dictionaries of local ROBs are commonly used for nonlinear model order reduction. A natural way of building such dictionaries is to partition the parameter space or the solution space with a clustering algorithm, and then define one local basis per cluster. This article shows that the choice of the dissimilarity measure for clustering is crucial, as it highly affects the performances of the local ROBs. In particular, it is shown that using Euclidean distances in the parameter space or in the solution space can lead to local bases whose performances are worse than those of a global ROB with the same number of modes. To remedy this problem, a ROM-oriented dissimilarity measure involving the principal angles between subspaces spanned by simulation results is introduced. It enables focusing on the shape of simulation results rather than their magnitudes. The strength of this dissimilarity comes from its link with the projection error appearing in the definition of the Kolmogorov $N$-width. The resulting dictionary of ROBs can be integrated in a ROM-net, where a classifier is used in the exploitation phase for fast and automatic model recommendation. The present paper gives an a priori efficiency criterion enabling hyperparameters calibration before time-consuming steps of the ROM-net's training phase. Future work will consider the application of this methodology to complex three-dimensional problems and their simulations with local ROMs.

\section*{Abbreviations}
\noindent DEIM: Discrete Empirical Interpolation Method; LDEIM: Localized Discrete Empirical Interpolation Method; PCA: Principal Component Analysis; PDE: Partial Differential Equation; POD: Proper Orthogonal Decomposition; ROB: Reduced-Order Basis; ROM: Reduced-Order Model; SVD: Singular Value Decomposition.

\section*{Funding}
\noindent Study funded by Safran and ANRT (Association Nationale de la Recherche et de la Technologie, France).

\section*{Appendix A: Properties of the normalized Kolmogorov width}

\begin{proper}
[Inequalities on Kolmogorov widths] If $\mathcal{M}$ is bounded and contains at least one nonzero element, then:\begin{equation}
\forall N \in \mathbb{N}^{*}, \quad d_{N}(\mathcal{M}) \leq \underset{v \in \mathcal{M}}{\sup} \ || v ||_{\mathcal{H}} \ \tilde{d}_{N}(\mathcal{M}) .
\label{InequalityKolmogorovWidths}
\end{equation}
\end{proper}

\begin{proof}
The boundedness of $\mathcal{M}$ implies the existence of $\underset{v \in \mathcal{M}}{\sup} \ || v ||_{\mathcal{H}}$, thus for all $u \in \mathcal{M}\setminus\{0\}$: \begin{equation}
|| u - \pi_{\mathcal{H}_N}(u) ||_{\mathcal{H}} \leq \underset{v \in \mathcal{M}}{\sup} \ || v ||_{\mathcal{H}} \frac{|| u - \pi_{\mathcal{H}_N}(u) ||_{\mathcal{H}}}{|| u ||_{\mathcal{H}}} ,
\end{equation} which implies Equation~\eqref{InequalityKolmogorovWidths}.
\end{proof}

\begin{proof} [\textbf{Proof of Property~\ref{LinkKolmogorovSineAngle}}]
Given that $\measuredangle_{\mathcal{H}} \left(u,\mathcal{H}_N\right) \in [0;\pi /2]$, the sine of the angle $\measuredangle_{\mathcal{H}} \left(u,\mathcal{H}_N\right)$ satisfies:\begin{equation}
\begin{array}{rcl}
\sin \measuredangle_{\mathcal{H}} \left(u,\mathcal{H}_N \right) & = & | \sin \measuredangle_{\mathcal{H}} \left(u,\mathcal{H}_N \right) | \\
& = & \sqrt{1 - \cos^2 \measuredangle_{\mathcal{H}} \left(u,\mathcal{H}_N \right)} \\
& = & \sqrt{1 - \cos^2 \underset{v \in \mathcal{H}_N}{\inf} \ \arccos \left( \displaystyle\frac{ | \langle u,v \rangle_{\mathcal{H}} |}{|| u ||_{\mathcal{H}} || v ||_{\mathcal{H}}} \right)} ,
\end{array}
\end{equation} and, since the function $\alpha \rightarrow \sqrt{1-\cos^2 \alpha}$ is increasing on the interval $[0;\pi/2]$:\begin{equation}
\sin \measuredangle_{\mathcal{H}} \left(u,\mathcal{H}_N \right) = \displaystyle \underset{v \in \mathcal{H}_N}{\inf} \ \sqrt{1 - \cos^2  \arccos \left( \frac{ | \langle u,v \rangle_{\mathcal{H}} |}{|| u ||_{\mathcal{H}} || v ||_{\mathcal{H}}} \right)}
= \displaystyle \underset{v \in \mathcal{H}_N}{\inf} \ \sqrt{1 - \frac{ \langle u,v \rangle_{\mathcal{H}}^2}{|| u ||_{\mathcal{H}}^2 || v ||_{\mathcal{H}}^2} } .
\label{eqProofA}
\end{equation} Furthermore:\begin{equation}
\frac{|| u - \pi_{\textrm{span}(\{v\})}(u) ||_{\mathcal{H}}^2}{|| u ||_{\mathcal{H}}^2} = \frac{|| u - \langle u, \frac{v}{|| v ||_{\mathcal{H}}} \rangle_{\mathcal{H}} \frac{v}{|| v ||_{\mathcal{H}}}||_{\mathcal{H}}^2}{|| u ||_{\mathcal{H}}^2} = 1 - \frac{ \langle u,v \rangle_{\mathcal{H}}^2}{|| u ||_{\mathcal{H}}^2 || v ||_{\mathcal{H}}^2} ,
\label{eqProofB}
\end{equation} where $\pi_{\textrm{span}(\{v\})}(u)$ is the orthogonal projection of $u$ on $\textrm{span}(\{v\})$. Using Equation~\eqref{eqProofB} in Equation~\eqref{eqProofA} yields:\begin{equation}
\sin \measuredangle_{\mathcal{H}} \left(u,\mathcal{H}_N \right) = \underset{v \in \mathcal{H}_N}{\inf} \frac{|| u - \pi_{\textrm{span}(\{v\})}(u) ||_{\mathcal{H}}}{|| u ||_{\mathcal{H}}} = \underset{v \in \mathcal{H}_N}{\inf} \ \underset{w \in \textrm{span}(\{v\})}{\inf} \frac{|| u - w ||_{\mathcal{H}}}{|| u ||_{\mathcal{H}}} .
\end{equation} Finally, given that $\underset{v \in \mathcal{H}_N}{\cup} \textrm{span}(\{v\})= \mathcal{H}_N$:\begin{equation}
\sin \measuredangle_{\mathcal{H}} \left(u,\mathcal{H}_N \right) = \underset{v \in \mathcal{H}_N}{\inf} \frac{|| u - v ||_{\mathcal{H}}}{|| u ||_{\mathcal{H}}} ,
\end{equation} which ends the proof.
\end{proof}

\section*{Appendix B: Properties of the dissimilarity measure}

\begin{proof} [\textbf{Proof of Property~\ref{SineDissAndProjErrorProp}}]
Let us first develop the square of the right-hand side of Equation~\eqref{firstFormulaDeltaTilde}, denoted by $f_{n}(u,v)^2$, using Equation~\eqref{DefOrthogProj}, the bilinearity of the $L^2$ inner product and the orthonormality of the bases $\Psi_{n}(u)$ and $\Psi_{n}(v)$:\begin{equation}
\begin{array}{ccl}
f_{n}(u,v)^2 & = & \displaystyle\sum_{i=1}^n || \psi_{i}(u) ||_{L^{2}(\Omega)}^2 - 2 \sum_{i=1}^n \sum_{j=1}^n \langle \psi_{i}(u),\psi_{j}(v) \rangle_{L^{2}(\Omega)}^2 \\
& & + \displaystyle\sum_{i=1}^n \sum_{j=1}^n \sum_{k=1}^n \langle \psi_{i}(u),\psi_{j}(v) \rangle_{L^{2}(\Omega)} \langle \psi_{i}(u),\psi_{k}(v) \rangle_{L^{2}(\Omega)} \langle \psi_{j}(v),\psi_{k}(v) \rangle_{L^{2}(\Omega)} \\
& = & n - 2 \displaystyle\sum_{i=1}^n \sum_{j=1}^n \langle \psi_{i}(u),\psi_{j}(v) \rangle_{L^{2}(\Omega)}^2 \\
& & + \displaystyle\sum_{i=1}^n \sum_{j=1}^n \sum_{k=1}^n \langle \psi_{i}(u),\psi_{j}(v) \rangle_{L^{2}(\Omega)} \langle \psi_{i}(u),\psi_{k}(v) \rangle_{L^{2}(\Omega)} \delta_{jk} \\
& = & n - \displaystyle\sum_{i=1}^n \sum_{j=1}^n \langle \psi_{i}(u),\psi_{j}(v) \rangle_{L^{2}(\Omega)}^2 ,
\end{array}
\end{equation} where $\delta_{jk}$ is the Kronecker delta function. Let $\g{C}\in\mathbb{R}^{n \times n}$ be the matrix whose entries are the inner products $\langle \psi_{i}(u),\psi_{j}(v) \rangle_{L^{2}(\Omega)}$. Its SVD reads $\g{C} = \g{V} \cos \bs{\Theta} \ \g{W}^T$ where $\bs{\Theta}$ is a diagonal matrix containing the principal angles $\theta_{k}(\mathcal{V}_{n}(u), \mathcal{V}_{n}(v))$, and where $\g{V}$ and $\g{W}$ are orthogonal matrices. Then, we obtain:\begin{equation}
\begin{array}{ccl}
f_{n}(u,v)^2 & = & n - \textrm{tr}(\g{C}^T \g{C}) \\
& = & n - \textrm{tr}\left( \g{W} \cos \left( \bs{\Theta} \right)^T \g{V}^T \g{V} \cos \left( \bs{\Theta} \right) \g{W}^T \right) \\
& = & n - \textrm{tr}\left( \g{W}^T \g{W} \cos \left( \bs{\Theta} \right)^T \g{V}^T \g{V} \cos \left( \bs{\Theta} \right) \right) \\
& = & n - \textrm{tr}\left( \cos \left( \bs{\Theta} \right)^T \cos \left( \bs{\Theta} \right) \right) \\
& = & n - \displaystyle\sum_{k=1}^n \cos^2 \theta_k (\mathcal{V}_{n}(u), \mathcal{V}_{n}(v)) \\
& = & \displaystyle\sum_{k=1}^n 1 - \cos^2 \theta_k (\mathcal{V}_{n}(u), \mathcal{V}_{n}(v)) \\
& = & \displaystyle\sum_{k=1}^n \sin^2 \theta_k (\mathcal{V}_{n}(u), \mathcal{V}_{n}(v)) \\
& = & \tilde{\delta}_{n}(u,v)^2 .
\end{array}
\end{equation} These equations remain true when exchanging $u$ and $v$, which ends the proof.
\end{proof}

\begin{proof} [\textbf{Proof of Property~\ref{PropHSdistance}}]
Since the Hilbert-Schmidt inner product on $HS(L^2 (\Omega))$ does not depend on the choice of the orthonormal basis of $L^2 (\Omega)$, let us choose a basis that is relevant for calculations. For $u \in L^{2}(\Omega \times [0;t_f])$, the $n$-dimensional elementary basis $\Psi_{n}(u)$ is completed with an orthonormal basis of the orthogonal complement of $\mathcal{V}_{n}(u)$ in $L^2 (\Omega)$. The resulting orthonormal basis of $L^2 (\Omega)$ is denoted by $\{ \psi_{k}(u) \}_{k \in \mathbb{N}^*}$, where the $n$ first basis vectors are those of the basis $\Psi_{n}(u)$. Let us now expand the term $|| \pi_{\mathcal{V}_{n}(u)} - \pi_{\mathcal{V}_{n}(v)} ||_{HS(L^{2}(\Omega))}^2$:\begin{equation}
|| \pi_{\mathcal{V}_{n}(u)} - \pi_{\mathcal{V}_{n}(v)} ||_{HS(L^{2}(\Omega))}^2 = || \pi_{\mathcal{V}_{n}(u)} ||_{HS(L^{2}(\Omega))}^2 + || \pi_{\mathcal{V}_{n}(v)} ||_{HS(L^{2}(\Omega))}^2 - 2 \langle \pi_{\mathcal{V}_{n}(u)} , \pi_{\mathcal{V}_{n}(v)} \rangle_{HS(L^{2}(\Omega))} .
\end{equation} Using the definition of the Hilbert-Schmidt inner product given by Equation~\eqref{HSinnerProd}, one has:\begin{equation}
\begin{array}{ccl}
\langle \pi_{\mathcal{V}_{n}(u)} , \pi_{\mathcal{V}_{n}(v)} \rangle_{HS(L^{2}(\Omega))} & = & \displaystyle\sum_{i=1}^{\infty} \langle \pi_{\mathcal{V}_{n}(u)} (\psi_{i}(u)) , \pi_{\mathcal{V}_{n}(v)} (\psi_{i}(u)) \rangle_{L^{2}(\Omega)} \\
& = & \displaystyle\sum_{i=1}^{n} \langle \psi_{i}(u) , \pi_{\mathcal{V}_{n}(v)} (\psi_{i}(u)) \rangle_{L^{2}(\Omega)} \\
& = & \displaystyle\sum_{i=1}^{n} \sum_{j=1}^{n} \langle \psi_{i}(u) , \psi_{j}(v) \rangle_{L^{2}(\Omega)}^2 ,
\end{array}
\end{equation} where the last equality results from the expression of $\pi_{\mathcal{V}_{n}(v)} (\psi_{i}(u))$ given by Equation~\eqref{DefOrthogProj}. Furthermore:\begin{equation}
\begin{array}{ccl}
|| \pi_{\mathcal{V}_{n}(u)} ||_{HS(L^{2}(\Omega))}^2 & = & \displaystyle\sum_{i=1}^{\infty} \langle \pi_{\mathcal{V}_{n}(u)} (\psi_{i}(u)) , \pi_{\mathcal{V}_{n}(u)} (\psi_{i}(u)) \rangle_{L^{2}(\Omega)} \\ 
& = & \displaystyle\sum_{i=1}^{n} \langle \psi_{i}(u) , \psi_{i}(u) \rangle_{L^{2}(\Omega)} \\ 
& = & n .
\end{array}
\end{equation} Similarly, one can prove that $|| \pi_{\mathcal{V}_{n}(v)} ||_{HS(L^{2}(\Omega))}^2 = n$. Finally:\begin{equation}
|| \pi_{\mathcal{V}_{n}(u)} - \pi_{\mathcal{V}_{n}(v)} ||_{HS(L^{2}(\Omega))}^2 = 2n - 2 \sum_{i=1}^{n} \sum_{j=1}^{n} \langle \psi_{i}(u) , \psi_{j}(v) \rangle_{L^{2}(\Omega)}^2 = 2 f_{n}(u,v)^2 = 2\tilde{\delta}_{n}(u,v)^2 ,
\end{equation} where $f_{n}(u,v)$ was introduced in the proof of Property~\ref{SineDissAndProjErrorProp}.
\end{proof}

\begin{proper}
For all $n \in [\![ 1 ; \mathcal{N} ]\!]$, the sine dissimilarity is a pseudometric on $L^{2}(\Omega \times [0;t_f])$.
\label{pseudometric}
\end{proper}

\begin{proof}
The sine dissimilarity $\tilde{\delta}_n$ is nonnegative and symmetric. Equation~\eqref{symFormulaDeltaTilde} implies that for all $u \in L^{2}(\Omega \times [0;t_f])$, $\tilde{\delta}_{n}(u,u)=0$. Equation~\eqref{symFormulaDeltaTilde} also yields:\begin{equation}
\tilde{\delta}_{n}(u,v) = \displaystyle\frac{1}{\sqrt{2}} || \pi_{\mathcal{V}_{n}(u)} - \pi_{\mathcal{V}_{n}(v)} ||_{HS(L^{2}(\Omega))} = \frac{1}{\sqrt{2}} || \pi_{\mathcal{V}_{n}(u)} - \pi_{\mathcal{V}_{n}(w)} + \pi_{\mathcal{V}_{n}(w)} - \pi_{\mathcal{V}_{n}(v)} ||_{HS(L^{2}(\Omega))} ,
\end{equation} so the triangle inequality on the Hilbert-Schmidt norm gives the triangle inequality $\tilde{\delta}_{n}(u,v) \leq \tilde{\delta}_{n}(u,w) + \tilde{\delta}_{n}(w,v)$.
\end{proof}

Note that $\tilde{\delta}_n(u,v)=0$ does not imply that $u=v$, which is the reason why the sine dissimilarity is not a metric on $L^{2}(\Omega \times [0;t_f])$. This is not a problem since we want this dissimilarity measure to be zero for all pairs of functions $(u,v) \in L^{2}(\Omega \times [0;t_f])^2$ whose trajectories over time in $L^{2}(\Omega)$ give the same POD approximation space. The next property shows the link between the sine dissimilarity and the Grassmann dissimilarity used in~\cite{ROM-net} for dictionary-based ROM-nets:

\begin{proper}
[Equivalence with the Grassmann dissimilarity for small angles] Given $n \in [\![ 1 ; \mathcal{N} ]\!]$, let $\bs{\theta}_n$ denote the vector of principal angles between $\mathcal{V}_{n}(u)$ and $\mathcal{V}_{n}(v)$ for two square-integrable functions $u$ and $v$. As $|| \bs{\theta}_n ||_2$ tends towards zero, the sine dissimilarity is asymptotically equivalent to the Grassmann dissimilarity $|| \bs{\theta}_n ||_2$, that is, using Bachmann-Landau notations:\begin{equation}
\tilde{\delta}_{n}(u,v) \underset{|| \bs{\theta}_n ||_2 \rightarrow 0}{\sim} || \bs{\theta}_n ||_2 .
\end{equation}
\end{proper}

\begin{proof}
One must show that:\begin{equation}
|| \sin \bs{\theta}_n ||_2 = || \bs{\theta}_n ||_2 + o(|| \bs{\theta}_n ||_2) .
\end{equation} As $|| \bs{\theta}_n ||_2$ tends towards zero:
\begin{equation}
|| \sin \bs{\theta}_n ||_{2} = \left( \sum_{i=1}^n \sin^2 \theta_{n,k} \right)^{1/2} = \left( \sum_{i=1}^n (\theta_{n,k} + o(\theta_{n,k}^2))^2 \right)^{1/2} = \left( || \bs{\theta}_n ||_{2}^2 + o(|| \bs{\theta}_n ||_{2}^2) \right)^{1/2} ,
\end{equation} which gives:\begin{equation}
|| \sin \bs{\theta}_n ||_{2} = || \bs{\theta}_n ||_{2} \sqrt{1+o(1)} = || \bs{\theta}_n ||_{2} + o(|| \bs{\theta}_n ||_{2}) .
\end{equation}
\end{proof}

\bibliographystyle{unsrt}
\bibliography{Biblio}

\begin{thebibliography}{10}

\bibitem{10.5555/2568435}
A.~Quarteroni and G.~Rozza.
\newblock {\em Reduced Order Methods for Modeling and Computational Reduction}.
\newblock Springer Publishing Company, Incorporated, 2013.

\bibitem{keiper2018reduced}
W.~Keiper, A.~Milde, and S.~Volkwein.
\newblock {\em Reduced-Order Modeling ({ROM}) for Simulation and Optimization:
  Powerful Algorithms as Key Enablers for Scientific Computing}.
\newblock Springer International Publishing, 2018.

\bibitem{PGDreview}
F.~Chinesta, P.~Ladeveze, and E.~Cueto.
\newblock A short review on model order reduction based on {P}roper
  {G}eneralized {D}ecomposition.
\newblock {\em Archives of Computational Methods in Engineering}, 18:395--404,
  11 2011.

\bibitem{PGDbook}
F.~Chinesta and E.~Cueto.
\newblock {\em {PGD}-Based Modeling of Materials, Structures and Processes}.
\newblock 01 2014.

\bibitem{RBmethodPrudhomme}
C.~Prud'homme, D.~Rovas, K.~Veroy, L.~Machiels, Y.~Maday, A.~Patera, and
  G.~Turinici.
\newblock Reliable real-time solution of parametrized partial differential
  equations: Reduced-basis output bound methods.
\newblock {\em Journal of Fluids Engineering}, 124:70, 03 2002.

\bibitem{RozzaReducedBasis}
G.~Rozza, D.~Huynh, and A.~Patera.
\newblock Reduced basis approximation and a posteriori error estimation for
  affinely parametrized elliptic coercive partial differential equations.
\newblock {\em Archives of Computational Methods in Engineering}, 15:1--47, 09
  2007.

\bibitem{cordier:hal-00417819}
L.~Cordier and M.~Bergmann.
\newblock {Proper Orthogonal Decomposition: an overview}.
\newblock In {\em {Lecture series 2002-04, 2003-03 and 2008-01 on
  post-processing of experimental and numerical data, Von Karman Institute for
  Fluid Dynamics, 2008.}}, page 46 pages. {VKI}, 2008.

\bibitem{RowleyPODGalerkin}
C.~Rowley, T.~Colonius, and R.~Murray.
\newblock Model reduction for compressible flow using {POD} and {G}alerkin
  projection.
\newblock {\em Physica D: Nonlinear Phenomena}, 189:115--129, 01 2003.

\bibitem{10.1093/imanum/dru066}
A.~Cohen and R.~DeVore.
\newblock {Kolmogorov widths under holomorphic mappings}.
\newblock {\em IMA Journal of Numerical Analysis}, 36(1):1--12, 03 2015.

\bibitem{GREIF2019216}
C.~Greif and K.~Urban.
\newblock Decay of the kolmogorov n-width for wave problems.
\newblock {\em Applied Mathematics Letters}, 96:216 -- 222, 2019.

\bibitem{CagniartBook}
N.~Cagniart, Y.~Maday, and B.~Stamm.
\newblock {\em Model Order Reduction for Problems with Large Convection
  Effects}, pages 131--150.
\newblock 01 2019.

\bibitem{nonino2019overcoming}
M.~Nonino, F.~Ballarin, G.~Rozza, and Y.~Maday.
\newblock Overcoming slowly decaying kolmogorov n-width by transport maps:
  application to model order reduction of fluid dynamics and fluid--structure
  interaction problems, 2019.

\bibitem{doi:10.1137/110823158}
W.~Dahmen, C.~Huang, C.~Schwab, and G.~Welper.
\newblock Adaptive {P}etrov--{G}alerkin methods for first order transport
  equations.
\newblock {\em SIAM Journal on Numerical Analysis}, 50(5):2420--2445, 2012.

\bibitem{dahmen_plesken_welper_2014}
W.~Dahmen, C.~Plesken, and G.~Welper.
\newblock Double greedy algorithms: Reduced basis methods for transport
  dominated problems.
\newblock {\em ESAIM: Mathematical Modelling and Numerical Analysis},
  48(3):623–663, 2014.

\bibitem{rim2020manifold}
D.~Rim, B.~Peherstorfer, and K.T. Mandli.
\newblock Manifold approximations via transported subspaces: Model reduction
  for transport-dominated problems, 2020.

\bibitem{taddei2020spacetime}
T.~Taddei and L.~Zhang.
\newblock Space-time registration-based model reduction of parameterized
  one-dimensional hyperbolic pdes, 2020.

\bibitem{localROB}
D.~Amsallem, M.~Zahr, and C.~Farhat.
\newblock Nonlinear model order reduction based on local reduced-order bases.
\newblock {\em International Journal for Numerical Methods in Engineering},
  pages 1--31, 2012.

\bibitem{localROB2}
K.~Washabaugh, D.~Amsallem, M.~Zahr, and C.~Farhat.
\newblock Nonlinear model reduction for {CFD} problems using local reduced
  order bases.
\newblock {\em 42nd AIAA Fluid Dynamics Conference}, 2012.

\bibitem{Lumley}
J.~Lumley.
\newblock The structure of inhomogeneous turbulent flows.
\newblock {\em Atm. Turb. and Radio Wave. Prop.}, pages 166--178, 1967.

\bibitem{Interpolation0}
D.~Amsallem and C.~Farhat.
\newblock An online method for interpolating linear parametric reduced-order
  models.
\newblock {\em SIAM Journal on Scientific Computing}, 33(5):2169--2198, 2011.

\bibitem{Interpolation1}
T.~Lieu and M.~Lesoinne.
\newblock Parameter adaptation of reduced order models for three-dimensional
  flutter analysis.
\newblock {\em AIAA Paper 2004-0888}, 2004.

\bibitem{Interpolation2}
T.~Lieu, C.~Farhat, and M.~Lesoinne.
\newblock {POD}-based aeroelastic analysis of a complete {F}-16 configuration:
  {ROM} adaptation and demonstration.
\newblock {\em AIAA Paper 2005-2295}, 2005.

\bibitem{Interpolation3}
T.~Lieu and C.~Farhat.
\newblock Adaptation of {POD}-based aeroelastic {ROM}s for varying {M}ach
  number and angle of attack: {A}pplication to a complete {F}-16 configuration.
\newblock {\em AIAA Paper 2005-7666}, 2005.

\bibitem{Interpolation4}
T.~Lieu, C.~Farhat, and M.~Lesoinne.
\newblock Reduced-order fluid/structure modeling of a complete aircraft
  configuration.
\newblock {\em Computer Methods in Applied Mechanics and Engineering},
  195:5730--5742, 2006.

\bibitem{Interpolation5}
T.~Lieu and C.~Farhat.
\newblock Adaptation of aeroelastic reduced-order models and application to an
  {F}-16 configuration.
\newblock {\em AIAA Journal}, 45:1244--1257, 2007.

\bibitem{Interpolation6}
D.~Amsallem and C.~Farhat.
\newblock Interpolation method for adapting reduced-order models and
  application to aeroelasticity.
\newblock {\em AIAA Journal}, 46(7):1803--1813, 2008.

\bibitem{doi:10.2514/1.J050233}
D.~Amsallem, J.~Cortial, and C.~Farhat.
\newblock Towards real-time computational-fluid-dynamics-based aeroelastic
  computations using a database of reduced-order information.
\newblock {\em AIAA Journal}, 48(9):2029--2037, 2010.

\bibitem{AMSALLEM2016373}
D.~Amsallem, R.~Tezaur, and C.~Farhat.
\newblock Real-time solution of linear computational problems using databases
  of parametric reduced-order models with arbitrary underlying meshes.
\newblock {\em Journal of Computational Physics}, 326:373 -- 397, 2016.

\bibitem{Interpolation7}
R.~Mosquera, A.~Hamdouni, A.~El~Hamidi, and C.~Allery.
\newblock {POD} basis interpolation via {I}nverse {D}istance {W}eighting on
  {G}rassmann manifolds.
\newblock {\em Discrete and Continuous Dynamical Systems, Series S.},
  12(6):1743--1759, 2018.

\bibitem{Interpolation8}
R.~Mosquera, A.~El~Hamidi, A.~Hamdouni, and A.~Falaize.
\newblock Generalization of the {N}eville-{A}itken {I}nterpolation {A}lgorithm
  on {G}rassmann {M}anifolds : {A}pplications to {R}educed {O}rder {M}odel.
\newblock \url{https://arxiv.org/pdf/1907.02831.pdf}, 2019.

\bibitem{CHOI2020109787}
Y.~Choi, G.~Boncoraglio, S.~Anderson, D.~Amsallem, and C.~Farhat.
\newblock Gradient-based constrained optimization using a database of linear
  reduced-order models.
\newblock {\em Journal of Computational Physics}, 423:109787, 2020.

\bibitem{LEE2020108973}
K.~Lee and K.T. Carlberg.
\newblock Model reduction of dynamical systems on nonlinear manifolds using
  deep convolutional autoencoders.
\newblock {\em Journal of Computational Physics}, 404:108973, 2020.

\bibitem{KimChoi2020}
Y.~Kim, Y.~Choi, D.~Widemann, and T.~Zohdi.
\newblock A fast and accurate physics-informed neural network reduced order
  model with shallow masked autoencoder, 2020.

\bibitem{EIM}
M.~Barrault, Y.~Maday, N.C. Nguyen, and A.T. Patera.
\newblock An empirical interpolation method: application to efficient
  reduced-basis discretization of partial differential equations.
\newblock {\em Comptes Rendus Mathematiques}, 339(9):666--672, 2004.

\bibitem{Ryckelynck2005}
D.~Ryckelynck.
\newblock A priori hyperreduction method: an adaptive approach.
\newblock {\em Journal of Computational Physics, Elsevier}, 202(1):346--366,
  2005.

\bibitem{ECSW}
C.~Farhat, P.~Avery, T.~Chapman, and J.~Cortial.
\newblock Dimensional reduction of nonlinear finite element dynamic models with
  finite rotations and energy-based mesh sampling and weighting for
  computational efficiency.
\newblock {\em International Journal for Numerical Methods in Engineering},
  98(9):625--662, 2014.

\bibitem{ECM}
J.A. Hernandez, M.A. Caicedo, and A.~Ferrer.
\newblock Dimensional hyper-reduction of nonlinear finite element models via
  empirical cubature.
\newblock {\em Computer methods in applied mechanics and engineering},
  313:687--722, 2017.

\bibitem{Grimberg2020}
S.~Grimberg, C.~Farhat, R.~Tezaur, and C.~Bou-Mosleh.
\newblock Mesh sampling and weighting for the hyperreduction of nonlinear
  {P}etrov-{G}alerkin reduced-order models with local reduced-order bases.
\newblock 08 2020.

\bibitem{kmeans}
J.B. MacQueen.
\newblock Some methods for classification and analysis of multivariate
  observations.
\newblock {\em Proceedings of 5-th Berkeley Symposium on Mathematical
  Statistics and Probability}, 1:281--297, 1967.

\bibitem{aggarwal2015data}
C.C. Aggarwal.
\newblock {\em Data Mining: The Textbook}.
\newblock Springer International Publishing, 2015.

\bibitem{kMedoidsPAM}
L.~Kaufmann and P.~Rousseeuw.
\newblock Clustering by means of medoids.
\newblock {\em Data Analysis based on the L1-Norm and Related Methods}, pages
  405--416, 01 1987.

\bibitem{aggarwal2013data}
C.C. Aggarwal and C.K. Reddy.
\newblock {\em Data Clustering: Algorithms and Applications}.
\newblock Chapman \& Hall/CRC Data Mining and Knowledge Discovery Series.
  Taylor \& Francis, 2013.

\bibitem{Hastie2005TheEO}
T.~Hastie, R.~Tibshirani, and J.H. Friedman.
\newblock {\em The Elements of Statistical Learning: Data Mining, Inference,
  and Prediction, 2nd Edition}.
\newblock Springer series in statistics. Springer, 2009.

\bibitem{10.1145/331499.331504}
A.~K. Jain, M.~N. Murty, and P.~J. Flynn.
\newblock Data clustering: A review.
\newblock {\em ACM Comput. Surv.}, 31(3):264--323, September 1999.

\bibitem{Saxena2017}
A.~Saxena, M.~Prasad, A.~Gupta, N.~Bharill, O.~Patel, A.~Tiwari, M.~Er,
  W.~Ding, and C.~Lin.
\newblock A review of clustering techniques and developments.
\newblock {\em Neurocomputing}, 267, 07 2017.

\bibitem{balabanov2021randomized}
O.~Balabanov and A.~Nouy.
\newblock Randomized linear algebra for model reduction—part ii: minimal
  residual methods and dictionary-based approximation.
\newblock {\em Advances in Computational Mathematics}, 47(2):1--54, 2021.

\bibitem{AmsallemHaasdonk}
D.~Amsallem and B.~Haasdonk.
\newblock {PEBL-ROM}: {P}rojection-error based local reduced-order models.
\newblock {\em Advanced Modeling and Simulation in Engineering Sciences}, 3, 12
  2016.

\bibitem{ROM-net}
T.~Daniel, F.~Casenave, N.~Akkari, and D.~Ryckelynck.
\newblock Model order reduction assisted by deep neural networks ({ROM}-net).
\newblock {\em Advanced Modeling and Simulation in Engineering Sciences},
  7(16), 2020.

\bibitem{HenriPOD}
T.~Henri and J.P. Yvon.
\newblock Convergence estimates of {POD}-{G}alerkin methods for parabolic
  problems.
\newblock volume 166, pages 295--306, 01 2006.

\bibitem{gohberg1990classes}
I.~Gohberg, S.~Goldberg, and M.A. Kaashoek.
\newblock {\em Classes of Linear Operators}.
\newblock Number vol. 1 in Classes of Linear Operators. Springer, 1990.

\bibitem{cheverry:cel-01587623}
C.~Cheverry and N.~Raymond.
\newblock {Handbook of Spectral Theory}.
\newblock Lecture, September 2019.

\bibitem{Djouadi2009}
S.~Djouadi.
\newblock On the optimality of the proper orthogonal decomposition and balanced
  truncation.
\newblock pages 4221 -- 4226, 01 2009.

\bibitem{Djouadi2012}
S.~Djouadi and S.~Sahyoun.
\newblock On a generalization of the proper orthogonal decomposition and
  optimal construction of reduced order models.
\newblock In {\em 2012 American Control Conference (ACC)}, pages 1436--1441,
  2012.

\bibitem{Sirovich}
L.~Sirovich.
\newblock Turbulence and the dynamics of coherent structures, {P}arts {I}, {II}
  and {III}.
\newblock {\em Quarterly of Applied Mathematics}, XLV:561 -- 590, 1987.

\bibitem{Chatterjee}
A.~Chatterjee.
\newblock An introduction to the proper orthogonal decomposition.
\newblock {\em Current Science}, 78:808 -- 817, 2000.

\bibitem{Meyer2003}
M.~Meyer and H.G. Matthies.
\newblock Efficient model reduction in non-linear dynamics using the
  {K}arhunen-{L}o\`{e}ve expansion and dual-weighted-residual methods.
\newblock {\em Computational Mechanics}, 31:179--191, 05 2003.

\bibitem{PhysRevE.89.022923}
A.~Iollo and D.~Lombardi.
\newblock Advection modes by optimal mass transfer.
\newblock {\em Phys. Rev. E}, 89:022923, Feb 2014.

\bibitem{doi:10.1137/17M1140571}
J.~Reiss, P.~Schulze, J.~Sesterhenn, and V.~Mehrmann.
\newblock The {S}hifted {P}roper {O}rthogonal {D}ecomposition: A mode
  decomposition for multiple transport phenomena.
\newblock {\em SIAM Journal on Scientific Computing}, 40(3):A1322--A1344, 2018.

\bibitem{MadayConvection}
N.~Cagniart, Y.~Maday, and B.~Stamm.
\newblock Model order reduction for problems with large convection effects.
\newblock {\em In: Chetverushkin B., Fitzgibbon W., Kuznetsov Y., Neittaanmäki
  P., Periaux J., Pironneau O. (eds) Contributions to Partial Differential
  Equations and Applications. Computational Methods in Applied Sciences}, 47,
  2019.

\bibitem{Zimmermann2017}
R.~Zimmermann, B.~Peherstorfer, and K.~Willcox.
\newblock Geometric subspace updates with applications to online adaptive
  nonlinear model reduction.
\newblock {\em SIAM Journal on Matrix Analysis and Applications}, 39, 11 2017.

\bibitem{10.1145/1618452.1618469}
T.~Kim and D.L. James.
\newblock Skipping steps in deformable simulation with online model reduction.
\newblock {\em ACM Trans. Graph.}, 28(5):1–9, December 2009.

\bibitem{doi:10.1137/151003660}
M.~Ohlberger and F.~Schindler.
\newblock Error control for the localized reduced basis multiscale method with
  adaptive on-line enrichment.
\newblock {\em SIAM Journal on Scientific Computing}, 37(6):A2865--A2895, 2015.

\bibitem{CasenaveAkkari}
F.~Casenave and N.~Akkari.
\newblock An error indicator-based adaptive reduced order model for nonlinear
  structural mechanics - {A}pplication to high-pressure turbine blades.
\newblock {\em Math. Comput. Appl.}, 24(2), 2019.

\bibitem{he2020insitu}
W.~He, P.~Avery, and C.~Farhat.
\newblock In-situ adaptive reduction of nonlinear multiscale structural
  dynamics models, 2020.

\bibitem{doi:10.1137/140989169}
B.~Peherstorfer and K.~Willcox.
\newblock Online adaptive model reduction for nonlinear systems via low-rank
  updates.
\newblock {\em SIAM Journal on Scientific Computing}, 37(4):A2123--A2150, 2015.

\bibitem{doi:10.1137/19M1257275}
B.~Peherstorfer.
\newblock Model reduction for transport-dominated problems via online adaptive
  bases and adaptive sampling.
\newblock {\em SIAM Journal on Scientific Computing}, 42(5):A2803--A2836, 2020.

\bibitem{ETTER2020112931}
P.A. Etter and K.T. Carlberg.
\newblock Online adaptive basis refinement and compression for reduced-order
  models via vector-space sieving.
\newblock {\em Computer Methods in Applied Mechanics and Engineering}, 364,
  2020.

\bibitem{doi:10.1137/120873868}
Y.~Maday and B.~Stamm.
\newblock Locally adaptive greedy approximations for anisotropic parameter
  reduced basis spaces.
\newblock {\em SIAM Journal on Scientific Computing}, 35(6):A2417--A2441, 2013.

\bibitem{Kaulmann2012ONLINEGR}
S.~Kaulmann and B.~Haasdonk.
\newblock Online greedy reduced basis construction using dictionaries.
\newblock {\em VI International Conference on Adaptive Modeling and Simulation
  (ADMOS)}, 2012.

\bibitem{Drohmann2010}
M.~Drohmann, B.~Haasdonk, and M.~Ohlberger.
\newblock Adaptive reduced basis methods for nonlinear convection–diffusion
  equations.
\newblock volume~4, pages 369--377, 12 2010.

\bibitem{Dihlmann2011}
M.~Dihlmann, M.~Drohmann, and B.~Haasdonk.
\newblock Model reduction of parametrized evolution problems using the reduced
  basis method with adaptive time partitioning.
\newblock 01 2011.

\bibitem{Eftang2010}
J.~Eftang, A.~Patera, and E.~Ronquist.
\newblock An “hp” certified reduced basis method for parametrized elliptic
  partial differential equations.
\newblock {\em SIAM J. Scientific Computing}, 32:3170--3200, 09 2010.

\bibitem{Haasdonk2011}
B.~Haasdonk, M.~Dihlmann, and M.~Ohlberger.
\newblock A training set and multiple bases generation approach for
  parametrized model reduction based on adaptive grids in parameter space.
\newblock {\em Mathematical and Computer Modelling of Dynamical Systems},
  17:423--442, 08 2011.

\bibitem{LDEIM2014}
B.~Peherstorfer, D.~Butnaru, K.~Willcox, and H.J. Bungartz.
\newblock {L}ocalized {D}iscrete {E}mpirical {I}nterpolation {M}ethod.
\newblock {\em SIAM Journal on Scientific Computing}, 36, 01 2014.

\bibitem{doi:10.2514/6.2020-0418}
M.G. Kapteyn, D.J. Knezevic, and K.E. Willcox.
\newblock Toward predictive digital twins via component-based reduced-order
  models and interpretable machine learning.
\newblock {\em AIAA Scitech 2020 Forum}, 2020.

\bibitem{kapteyn2020physicsbased}
M.G. Kapteyn and K.E. Willcox.
\newblock From physics-based models to predictive digital twins via
  interpretable machine learning, 2020.

\bibitem{Amsallem2015localHROM}
D.~Amsallem, M.~Zahr, and K.~Washabaugh.
\newblock Fast local reduced basis updates for the efficient reduction of
  nonlinear systems with hyper-reduction.
\newblock {\em Advances in Computational Mathematics}, 41, 02 2015.

\bibitem{RyckelynckComputerVision}
F.~Nguyen, S.M. Barhli, D.P. Mu\~{n}oz, and D.~Ryckelynck.
\newblock Computer vision with error estimation for reduced order modeling of
  macroscopic mechanical tests.
\newblock {\em Complexity}, 2018.

\bibitem{BuhrIapichinoSmetana+2020+245+306}
A.~Buhr, L.~Iapichino, and K.~Smetana.
\newblock {\em 6 Localized model reduction for parameterized problems}, pages
  245--306.
\newblock De Gruyter, 2020.

\bibitem{lee2020DeepConservation}
K.~Lee and K.~Carlberg.
\newblock {D}eep {C}onservation: A latent-dynamics model for exact satisfaction
  of physical conservation laws, 2020.

\bibitem{DEIM}
S.~Chaturantabut and D.~Sorensen.
\newblock Discrete empirical interpolation for nonlinear model reduction.
\newblock {\em Decision and Control, 2009 held jointly with the 2009 28th
  Chinese Control Conference, CDC/CCC 2009, Proceedings of the 48th IEEE
  Conference}, pages 4316--4321, 2010.

\bibitem{LocalPODGPR}
R.~Dupuis, J.-C. Jouhaud, and P.~Sagaut.
\newblock Surrogate modeling of aerodynamic simulations for multiple operating
  conditions using machine learning.
\newblock 12 2019.

\bibitem{mainini2015surrogate}
Laura Mainini and Karen Willcox.
\newblock Surrogate modeling approach to support real-time structural
  assessment and decision making.
\newblock {\em AIAA Journal}, 53(6):1612--1626, 2015.

\bibitem{mca26010017}
T.~Daniel, F.~Casenave, N.~Akkari, and D.~Ryckelynck.
\newblock Data augmentation and feature selection for automatic model
  recommendation in computational physics.
\newblock {\em Mathematical and Computational Applications}, 26(1), 2021.

\bibitem{10.2307/2529785}
A.~D. Gordon and J.~T. Henderson.
\newblock An algorithm for euclidean sum of squares classification.
\newblock {\em Biometrics}, 33(2):355--362, 1977.

\bibitem{kaufman1990finding}
L.~Kaufman, P.J.R. Leonard~Kaufman, and P.J. Rousseeuw.
\newblock {\em Finding Groups in Data: An Introduction to Cluster Analysis}.
\newblock A Wiley-Interscience publication. Wiley, 1990.

\bibitem{kMedoidsCLARA}
L.~Kaufman and P.~Rousseeuw.
\newblock {\em Clustering Large Data Sets}, pages 425--437.
\newblock 12 1986.

\bibitem{kMedoidsCLARANS2}
R.T. Ng and J.~Han.
\newblock Efficient and effective clustering methods for spatial data mining.
\newblock Technical report, CAN, 1994.

\bibitem{kMedoidsCLARANS}
R.~Ng and J.~Han.
\newblock {CLARANS}: A method for clustering objects for spatial data mining.
\newblock {\em Knowledge and Data Engineering, IEEE Transactions on}, 14:1003--
  1016, 10 2002.

\bibitem{kMedoids2019}
E.~Schubert and P.J. Rousseeuw.
\newblock Faster k-medoids clustering: Improving the {PAM}, {CLARA}, and
  {CLARANS} algorithms.
\newblock In G.~Amato, C.~Gennaro, V.~Oria, and M.~Radovanovi{\'{c}}, editors,
  {\em Similarity Search and Applications}, pages 171--187, Cham, 2019.
  Springer International Publishing.

\bibitem{kMedoidsParkJun2009}
H.S. Park and C.H. Jun.
\newblock A simple and fast algorithm for k-medoids clustering.
\newblock {\em Expert Systems with Applications}, 36:3336--3341, 2009.

\bibitem{Bjorck73}
A.~Bjorck and G.~Golub.
\newblock Numerical methods for computing angles between linear subspaces.
\newblock {\em Mathematics of Computation}, 27:123, 07 1973.

\bibitem{chordalDist}
J.~Conway, R.~Hardin, and N.~Sloane.
\newblock Packing lines, planes, etc.: Packings in grassmannian space.
\newblock {\em Experimental Mathematics}, 5:139--159, 01 1996.

\bibitem{GappyPOD}
R.~Everson and L.~Sirovich.
\newblock Karhunen-{L}oeve procedure for gappy data.
\newblock {\em JOSA A}, 12, 08 1995.

\bibitem{OptimalPartitions-DanielCasenaveAkkariRyckelynck}
T.~Daniel, F.~Casenave, N.~Akkari, and D.~Ryckelynck.
\newblock Optimal piecewise linear data compression for solutions of
  parametrized partial differential equations.
\newblock {\em arXiv preprint: 2108.12291}, 2021.

\bibitem{bachmayr2017kolmogorov}
M.~Bachmayr and A.~Cohen.
\newblock Kolmogorov widths and low-rank approximations of parametric elliptic
  {PDE}s.
\newblock {\em Mathematics of Computation}, 86(304):701--724, 2017.

\bibitem{ryckelynck:hal-02559710}
D.~Ryckelynck, T.~Goessel, and F.~Nguyen.
\newblock {Mechanical dissimilarity of defects in welded joints via {G}rassmann
  manifold and machine learning}.
\newblock Preprint, July 2020.

\bibitem{yales2}
Vincent Moureau, Pascale Domingo, and Luc Vervisch.
\newblock Design of a massively parallel {CFD} code for complex geometries.
\newblock {\em Comptes Rendus M{\'e}canique}, 339(2-3):141 -- 148, 2011.

\bibitem{zset}
{M}ines {P}aris{T}ech and {ONERA} the {F}rench aerospace lab. {Z}-set:
  nonlinear material \& structure analysis suite.
\newblock \url{http://www.zset-software.com}, 1981-present.

\end{thebibliography}

\end{document}